\documentclass[12pt]{article}

\usepackage{epic}
\usepackage{eepic}
\usepackage{epsf} 
\usepackage{ amssymb, amscd}
\usepackage{amsmath}
\usepackage{color}
\numberwithin{equation}{section}
\usepackage{epsfig}
\usepackage{graphicx}
\usepackage{tikz}
\usetikzlibrary{matrix}

\def\cA            {{\mathcal{A}}}
\def\cB            {{\mathcal{B}}}
\def\cC            {{\mathcal{C}}}
\def\cD            {{\mathcal{D}}}
\def\cE            {{\mathcal{E}}}
\def\cF            {{\mathcal{F}}}
\def\cG            {{\mathcal{G}}}

\def\cL            {{\mathcal{L}}}
\def\cM            {{\mathcal{M}}}
\def\cN            {{\mathcal{N}}}
\def\cO            {{\mathcal{O}}}

\def\cV            {{\mathcal{V}}}

\def\bbC           {\mathbb{C}}

\def\bbH           {\mathbb{H}}

\def\bbX           {\mathbb{X}}

\def\bbQ           {\mathbb{Q}}
\def\bbR           {\mathbb{R}}
\def\bbT           {\mathbb{T}}
\def\bbZ           {\mathbb{Z}}

\def\NXN           {{}_N {\mathcal{X}}_N}

\def\i{{\rm i}}
\def\sdprod{{\times\!\vrule height5pt depth0pt width0.4pt\,}}

\title{Non-unitary fusion categories and their doubles via endomorphisms}
\author{
{\sc David E.\ Evans}\\
 {\footnotesize School of Mathematics, Cardiff University,}\\
 {\footnotesize Senghennydd Road, Cardiff CF24 4AG, Wales, U.K.}\\
 {\footnotesize e-mail: {\tt EvansDE@cf.ac.uk}}\\ \\
 {\sc Terry  Gannon }\\
 {\footnotesize Department of Mathematics, University of Alberta,}\\
{\footnotesize Edmonton, Alberta, Canada T6G 2G1}\\
{\footnotesize e-mail: {\tt tgannon@math.ualberta.ca}} }

\begin{document}
\maketitle

\begin{abstract} 
We realise non-unitary fusion categories using subfactor-like methods,
and compute their quantum doubles and modular data. For concreteness we focus on {}{generalising}
the Haagerup--Izumi family of {}{Q-systems.} For example, we construct endomorphism
realisations of  the (non-unitary)
Yang--Lee model, and  non-unitary analogues of one of the even subsystems of the Haagerup subfactor and of the Grossman--Snyder system. We supplement Izumi's equations for identifying the half-braidings, which were incomplete even in {}{his Q-system} setting.  We conjecture a remarkably simple form
for the modular $S$ and $T$ matrices of the doubles of these fusion categories. We would expect all of these
doubles to be realised as the category of modules of a rational VOA {}{and conformal net of factors}. We expect our approach will also suffice to realise the non-semisimple tensor categories arising in logarithmic conformal field theories.

\end{abstract}

{\footnotesize
\tableofcontents
}

\section{Introduction}

{}{The chiral part of a  \textit{unitary} rational conformal field theory (CFT) can be represented as either a 
completely rational conformal net of factors
on a circle or a rational vertex operator algebra (VOA). Whilst conformal nets and subfactors  theory focus and exploit the analytic aspects, %and elliptic genus and Chiral de Rham complexes 
%exploit the geometry, 
vertex operator algebras focus on the algebraic aspects. %is constrained by a purely algebraic setting.
The relation between these approaches is studied in \cite{CKLW}; at the simplest level, they both must give rise
to the same modular tensor category (MTC)  if they are to correspond to the same CFT. }

{}{Conformal nets of factors are a particularly rich framework, % within algebraic quantum field  theory, 
with connections with twisted equivariant $K$-theory and non-commutative geometry. 
Subfactor methods have proved to be much more effective {than VOA methods in many ways.}
For example, 
structure theorems such as rationality of orbifolds or cosets is much
easier in the conformal nets of subfactors picture (see e.g. \cite{Xu}) than in the VOA picture.
{Also, the factor %, and indeed the subfactor 
setting captures in a natural way} the \textit{full} CFT as an inclusion of (local) nets \cite{BE4}, \cite{Rehren}}.

{}{However, the VOA setting for the chiral CFT is apparently more flexible in allowing non-unitary examples. For example,
the Virasoro minimal models are parametrised by pairs $p>q$ of coprime numbers;
they  are unitary if and only if $p=q+1$. The simplest of these is the Yang--Lee model $V(2,5)$ (see e.g. section 7.4.1 of \cite{DMS}), which Cardy \cite{Cardy} showed arises as the Yang--Lee edge singularity in the Ising model in an imaginary magnetic field. {Other non-unitary statistical mechanical examples are the scaling limit
of critical dense polymers, and critical perculation, both with central charge $c=-2$. An unrelated non-unitary example crucial to string theory is the (super-)ghost CFT; what must be unitary is space where the physical states lie, namely the 
BRST cohomology of the ghosts coupled to  a matter CFT. Wess--Zumino--Witten models on Lie
supergroups provide other non-unitary examples important to string theory.}} {In the VOA setting, realising non-unitary CFTs presents no special problems, whereas subfactors and nets of factors have unitarity built in.}

{A fundamental question is
whether there are any rational CFTs beyond those constructed from loop groups or quantum groups, 
using standard methods such as orbifolds and cosets (see e.g. \cite{Xu} for a discussion on this point). 
It is known that all unitary fusion categories, hence all unitary MTCs, can be realised
by endomorphisms on a factor. These methods have produced countless 
 `exotic' examples of \textit{unitary} MTCs  \cite{EG2,EG4}. Indeed, the relative abundance of these examples suggests that most modular tensor categories may be
 `exotic'.  Finding conformal net and VOA realisations of these `exotic' MTCs is an important but 
difficult challenge --- we expect most or all of them to have such realisations. The situation for the  (double of the) Haagerup subfactor is discussed in detail in \cite{EG2}. In any case, the effectiveness of these subfactor methods in
constructing new \textit{unitary} MTCs provides another compelling reason for extending
these methods to the non-unitary setting.}

%In this paper we explain how to extend the subfactor methods to the non-unitary setting.
The main purpose of this paper is to provide a broader context, {dropping the requirement of unitarity,} in which the subfactor methods
can be applied. {After all, most rational CFT are non-unitary, and one would like to exploit the powerful methods of subfactors and nets of factors in the general case.}
%We focus on the Haagerup--Izumi series for concreteness. We realise the (non-unitary) Yang--Lee model, and ... 

{In the remainder of the Introduction we sketch in more detail some of the terms used earlier, as
well as the content of the paper.}

The sectors of a rational CFT, or modules of a rational VOA, give rise to a tensor category
of a very special type, namely an MTC. More generally, we are interested
in \textit{fusion categories}, which {roughly speaking} are MTCs without the braiding (we review their definition  in section 3). {Given a fusion category, the \textit{double} or \textit{centre construction} canonically associates
an MTC.}  Unitarity {in a category} can be defined
as follows.  A $*$-operation on a $\bbC$-linear category $\cC$ is a {conjugate-linear involution  $\mathrm{Hom}(X,Y)\rightarrow\mathrm{Hom}(Y,X)$ satisfying $(fg)^*=g^*f^*$ for all
$f\in\mathrm{Hom}(X,Y)$, $g\in\mathrm{Hom}(Z,Y)$. If the category is tensor (and strict),  we also 
require $(f\otimes g)^*=f^*\otimes g^*$ for all $f\in\mathrm{Hom}(X,Y)$, $g\in\mathrm{Hom}(Z,W)$.} 
A $*$-operation is called positive if $f^* f = 0$ implies $f = 0$. 
A category equipped with a (positive) $*$-operation is called \textit{hermitian} (resp. \textit{unitary}).

{Associated to an MTC is a representation of SL$_2(\bbZ)$ called
\textit{modular data}. It is generated by a symmetric unitary matrix $S$ which gives the fusion coefficients
({structure constants of} the Grothendieck ring of the category)
through Verlinde's formula, together with a diagonal matrix $T$ of finite order. Some column of $S$ must
be strictly positive --- e.g. in a \textit{unitary} MTC  that {Perron--Frobenius} column corresponds to the
unit. In a rational CFT, the characters {$\chi_M(\tau)=q^{h_M-c/24}\sum_{n=0}^\infty
\mathrm{dim}\,M_n\,q^n$} of the irreducible modules {$M=\coprod_nM_n$ form a vector-valued modular function for SL$_2(\bbZ)$ with modular data
as its multiplier. The minimal conformal weight $h_M$ corresponds to the positive column of $S$.
The conformal weights and central charge $c$ must be rational, but in a unitary theory they will also be non-negative.
For more comparisons between the modular data of non-unitary versus unitary theories, see \cite{Ganun}.}

{A very convenient realisation of tensor categories is through endomorphisms on an
algebra, where objects are algebra endomorphisms and morphisms are intertwiners. The tensor product
of objects corresponds to composition and of morphisms to the (twisted) product in the underlying
algebra. However, it is awkward to realise other properties in the category, such as additivity
or rigidity, without assuming special structures on the algebra. When the underlying algebra is
a $C^*$-algebra such as the Cuntz algebra, these other properties arise naturally. Indeed,
any \textit{unitary} fusion category can be realised as a category of endomorphisms on} {a hyperfinite von Neumann algebra} {(see section 7 of \cite{HaYa})}.

{A natural question is, can we find systematic realisations by endomorphisms
of \textit{non-unitary} fusion categories? We will see that the answer is yes.}

{Our approach was influenced by} recent work of Phillips \cite{Phil}, who studies
non-unitary analogues of the Cuntz algebra. But {all of our calculations are within 
a polynomial algebra (the \textit{Leavitt algebra}).} Rather than completing {that
 algebra} as {studied} by Phillips, we have found it {sufficient} to work exclusively within the
Leavitt algebra itself.

For concreteness we focus on the \textit{Haagerup--Izumi {family of} fusion rings}, but our method works more
generally. Let $G$ be any finite abelian group. The (isomorphism classes of) simple objects in
these fusion rings are $[\alpha_g]$ and $[\alpha_g\rho]$ as $g$ ranges over $G$. The fusions are given by
\begin{align}[\alpha_g][\alpha_h]
=[\alpha_{g+h}]\,,& \ [\alpha_g][\alpha_h\rho]=[\alpha_{g+h}\rho]=[\alpha_h\rho][\alpha_{-g}]\,, \nonumber\\
&[\alpha_g\rho][\alpha_h\rho]=[\alpha_{g-h}]+\sum{}_k\,[\alpha_k\rho]\,.\label{HIfus}\end{align}
{In the following sections} we explain explicitly how to construct, using endomorphisms on {the Leavitt} algebra, fusion categories (not necessarily unitary) which realise
the Haagerup--Izumi fusions when $G$ has odd order. We compute the corresponding tube algebras and from that
obtain the modular data $S,T$ of the double of the system. We give several examples {and
{explicitly} classify these systems for small $G$}.

The {(unitary)} Haagerup--Izumi {fusions} \eqref{HIfus} for $|G|$ odd was introduced  by Izumi  in \cite{iz3}. His motivation was to construct the Haagerup subfactor \cite{Haag,ahaag}, so he focussed on the special class of
systems of  {Cuntz algebra endomorphisms, called \textit{Q-systems}}, 
 which arise as the even subsystem of a
subfactor with canonical endomorphism $1+\rho$. {Q-systems correspond to especially {constrained
$\rho$; their} fusion categories are always unitary.} He showed that there was a unique {Q-}system {satisfying \eqref{HIfus}} for the group 
$G=\bbZ_3$, and {comparing indices observed that  it must} correspond to the Haagerup subfactor. {Likewise, he showed} that there is a unique {Q-}system  for $G=\bbZ_5$.
He also computed the modular data for the doubles of his systems {(modulo a 
technicality discussed shortly)}. {Evans--Gannon} \cite{EG2} pushed
this further, finding {Q-}systems in {this} class for {all $G$ with} $|G|\le 19$ (including the complete lists for $|G|\le 9$),
and simplifying considerably Izumi's expressions for the modular data. {Thanks to this
work,} it is now expected that there
are subfactors (usually several) for each odd order, {and they are all expected to correspond
through their doubles to rational VOAs etc}. Grossman--Snyder \cite{GrSn}
found new systems of endomorphisms realising \eqref{HIfus} ({unitary but} not Q-systems), for
$G=\bbZ_3$ and $\bbZ_5$, which are Morita equivalent to Izumi's systems (and thus have the 
same doubles).
This treatment has been extended to even order $G$, and to all unitary systems {(not only Q-systems)} {realising \eqref{HIfus}}, by 
{Evans--Gannon \cite{EGm},\cite{EGn} and independently Izumi  \cite{iz4}}.  

{In this paper, as an illustration of our method, we characterise all realisations
by endomorphisms (not necessarily Q-systems {nor unitary}) of the Haagerup--Izumi fusions \eqref{HIfus} for $|G|$ odd. We show they all yield fusion
categories.  Like \cite{iz3}, our systems correspond to solutions
of finitely many equations in finitely many variables, but unlike \cite{iz3} our equations are all
polynomials (those of \cite{iz3} involve complex conjugates).
%, and ours include quartic equations (all of those of \cite{iz3} are cubic or less). 
In broad strokes the method we use is analogous
to that of \cite{iz3}, but the absence of unitarity introduces several complications and our
argument is required to be much more subtle. We find the doubles and modular data of our systems.}

For example, there are
precisely 2,4,4 inequivalent fusion categories realised by endomorphisms, of Haagerup--Izumi type for $G=\bbZ_1,\bbZ_3,\bbZ_5$
respectively (of course we recover all of them). Precisely 1,2,2 of these, respectively, are unitary: 1,1,1 are {Q-systems,} %in the Izumi class, 
and  0,1,1 are the aforementioned Grossman--Snyder systems. The Yang--Lee system is the {unique} non-unitary
one corresponding to $G=\bbZ_1$. 

Every fusion category $\cC$ is defined over some number field \cite{ENO}. An automorphism
$\sigma$ of that field acts on the quantities of that category in the natural way, defining a 
new fusion category $\cC^\sigma$. These categories may or may not be equivalent ---
e.g. a Galois associate of a unitary fusion category may not be unitary.  In general, $\cC$ and 
$\cC^\sigma$ will
have identical fusion rings, but their modular data for example will be Galois associates.
Our construction, {unlike that of e.g. Izumi,} is closed under this Galois action.

{It turns out that} all 5 non-unitary fusion categories for $G=\bbZ_1,\bbZ_3,\bbZ_5$ are
Galois associates of unitary categories. We expect though that this is an accident of small $G$. Our system of 
equations involve twice as many  variables as in the unitary case, and approximately the same
number of equations. %, and most of ours are quartic rather than cubic. 
{For these reasons,  we would} expect typically
many more non-unitary categories than unitary ones.

{In any case, it is easy to construct non-unitary fusion categories, all of whose Galois associates
are also non-unitary. A simple example is the tensor product of affine $G_2$ at level 1 (a unitary MTC) with the Yang--Lee model {(a non-unitary one)}.}

Actually, the equations in \cite{iz3} are not sufficient to determine the half-braidings, and hence
the modular data, for most odd abelian $G$, even in the {Q-system} case.
In section 6 below we supply additional equations which are {both necessary and} sufficient.

Incidentally, another interesting class of CFTs and  VOAs are the so-called \textit{logarithmic} {or \textit{$C_2$-cofinite non-rational}} ones \cite{CrRi}, {for example the symplectic fermions \cite{DaRu}}. Unlike the
rational CFTs, their category of modules will not be semisimple {and so direct (sub)factor realisations
of them wouldn't be possible}.   Logarithmic theories
appear to be intimately connected with non-unitarity: all known ones are conformal embeddings
of non-unitary rational VOAs ({with states of negative conformal weight}). In any case, although we address in this paper only fusion categories (which are semisimple by definition),
 modelling non-semisimple systems is also possible by our methods and
we would expect we could realise with endomorphisms these logarithmic theories.

\section{The Yang-Lee model}

{This section} illustrates the ideas developed in the following sections, with {the simplest} non-unitary example:
 the Yang-Lee model  ({this CFT is described e.g. in} section 7.4.1 of \cite{DMS}). It consists of
two simple objects 1 and $\rho$, which obey the fusion rule
\begin{equation}
[\rho][\rho]=[1]+[\rho]\,.\label{YangLee}\end{equation}

Let us try to realise \eqref{YangLee} as a system of algebra endomorphisms on some algebra {$\cA$. To motivate our solution though, let's reverse the logic and derive the consequences of such a realisation. It would}  require the relation
\begin{equation}\rho(\rho(x))=sxs'+t\rho(x)t'\,,\label{ylrho2}\end{equation}
where $s,s',t,t'\in\cA$ satisfy the Leavitt--Cuntz relations
\begin{equation} ss'+tt'=1\,,\ \ s's=t't=1\,,\ \ s't=t's=0\,.\label{LC2}\end{equation}
{More precisely, these relations say that \eqref{ylrho2} expresses $\rho\circ\rho$ as a direct sum
of objects id and $\rho$ in the category ${\cE\cN\cD}(\cA)$ (we describe this category in detail
next section).} 
These {elements $s,s',t,t'$ generate {by definition} a copy of the Leavitt algebra $\cL_2$ inside $\cA$; we will see shortly that $\rho$ restricts to an endomorphism of $\cL_2$}.
In order  to identify {the restriction of $\rho$ to $\cL_2$},  it is necessary and sufficient to determine
the values $\rho(s)$, $\rho(s')$, $\rho(t)$ and $\rho(t')$ {of $\rho$ on the generators}. For $*$-maps, we would have $\rho(s')=\rho(s)'$ etc,
but we cannot require that here  if we
hope to realise the Yang--Lee model.

{We {require that} both endomorphisms id and $\rho$ be simple, equivalently that the intertwiner spaces Hom(id,id) and Hom$(\rho,\rho)$ in the algebra $\cA$ be $\bbC 1$, 
and that Hom$(\rho,\mathrm{id})=\mathrm{Hom}(\mathrm{id},\rho)=0$}. {(The definition of intertwiners is given next section.)} From \eqref{ylrho2} we obtain
 \begin{equation}
 \rho^2(x)s=sx\,,\ \  \rho^2(x)t=t\rho(x)\,,\ \  s'\rho^2(x)=xs'\,,\ \ 
 t'\rho^2(x)=\rho(x)t'\,.\label{YLintert}\end{equation}
{The first means $s\in\mathrm{Hom}(\mathrm{id},\rho^2)$. Conversely, suppose
$r\in\mathrm{Hom}(\mathrm{id},\rho^2)$, i.e.  $rx=\rho^2(x)r$ for all $x$. Then $s'rx=s'\rho^2(x)r
=xs'r$ and $t'rx=t'\rho^2(x)r=\rho(x)t'r$. Thus by simplicity of id and $\rho$ we have $s'r\in\bbC$ 
and $t'r=0$, so $r=(ss'+tt')r=ss'r\in\bbC s$ using the Leavitt--Cuntz relation $ss'+tt'=1$. We have shown Hom$(\mathrm{id},\rho^2)=\bbC s$.
 In the same way (see Lemma 3 below for details and the generalisation),} we can identify the intertwiner spaces 
 $\bbC t=\mathrm{Hom}(\rho,\rho^2)$,  $\bbC s'=\mathrm{Hom}(\rho^2,\mathrm{id})$ and 
  $\bbC t'=\mathrm{Hom}(\rho^2,\rho)$. {These observations are} crucial for what follows.

Note, using \eqref{YLintert}, that
\begin{equation}s'\rho(s)\rho(x)=s'\rho(sx)=s'\rho(\rho^2(x)s)
=s'\rho^2(\rho(x))\rho(s)=\rho(x)s'\rho(s)\,.\end{equation}
In other words, $s'\rho(s)\in\mathrm{Hom}(\rho,\rho)=\bbC$, so $s'\rho(s)$ equals some complex
number $a$. Likewise, $t'\rho(s)\in\mathrm{Hom}(\rho,\rho^2)$ so $t'\rho(s)=bt$ for some $b\in\bbC$.
{The point is that}
 \begin{equation}\rho(s)=(ss'+tt')\rho(s)=s(s'\rho(s))+t(t'\rho(s))=as+btt\,.\label{YLrhos}\end{equation}
Similar calculations (see section 4 for details and the generalisation) give \begin{equation}
\rho(s')=a's'+b't't'\,,\ \ \rho(t)=cst'+dtss'+ettt'\,,\ \ \rho(t')=c'ts'+d'ss't'+e'tt't'\,,\label{YLrho}\end{equation}
for some $a',b',c,c',d,d',e,e'\in\bbC$. {Because $\rho$ sends the generators of $\cL_2$ into $\cL_2$,} 
this means $\rho$ is actually an endomorphism of $\cL_2$. 
If we required $\rho$ to be a $*$-map, then {we would have} $a'=\overline{a}$ etc, but again we shouldn't do that
{if we are to recover Yang-Lee}.

{We can now use the constraints on $\rho$ to solve for those 10 parameters. First, $\rho$ is required
to be an algebra endomorphism, so it must respect the}  Leavitt--Cuntz relations \eqref{LC2}.
One relation requires $1=\rho(s')\rho(s)$, i.e.
\begin{equation}\label{aabb1}
1=(a's'+b't't')(as+btt)=a'a+{bb'}\,.\end{equation}
Similarly,  $\rho(s)\rho(s')+\rho(t)\rho(t')=1$ gives the identities {$1=aa'+cc'$
(hence $b'b=c'c$), $aa'+cc'=dd'$ (hence $d'd=1$), and  $ab'=-ce'$, amongst others.
More precisely, Lemma 1 below gives a unique form for any element of a Leavitt algebra, so once we expand out  $\rho(s)\rho(s')+\rho(t)\rho(t')=1$ and put it into reduced form (e.g. replacing
$ss'$ by $1-tt'$), the identities fall out by comparing corresponding coefficients.} 

We also require that $\rho$ satisfy \eqref{ylrho2}.  It {implies for instance} that $s'\rho^2(s)=ss'$. We can compute 
$s'\rho(\rho(s))$ directly from \eqref{YLrhos},\eqref{YLrho}, and we find
\begin{equation}
s'\rho^2(s)=as'\rho(s)+bs'\rho(t)\rho(t)=a^2+bct'(cst'+dtss'+ettt')=a^2+bcdss'+bcett'\,.\nonumber\end{equation}
{This must equal $ss'$, which {(using $1=ss'+tt'$)} gives} $1=a^2+bcd$ and $1=bcd-bce$ (hence $a^2=-bce$). Likewise, $\rho^2(s')s=ss'$ gives  {$1=a^{\prime \,2}+b'c'd'$ and $a^{\prime\,2}=-b'c'e'$. Similarly, \eqref{ylrho2} implies 
$t'\rho^2(s)=\rho(s)t'\,$; its $t$ and $st'$ coefficients give $ab=-bde$ and $a=bcd$ respectively. Likewise,
$\rho^2(s')t=t\rho(s')$ gives $a'b'=-b'd'e'$ and $a'=b'c'd'$.}

{Plugging  $a=bcd$ into $1=a^2+bcd$ (and likewise for the {primed quantities}) gives $1=a^2+a=a'{}^2+a'$, 
which means $a,a'\in\{(-1\pm\sqrt{5})/2\}$. {Note} that if $a\ne a'$ then $aa'=1$ --- we
will use this shortly.
Since $a=bcd,a'=b'c'd'$ are both  non-zero, so are all $b,c,d,b',c',d'$.}
 Note that we are free to rescale $s$ by $\lambda\in\bbC^\times$
(hence $s'$ by $1/\lambda$) without affecting {\eqref{ylrho2} nor the Leavitt--Cuntz relations. Choosing $\lambda$
appropriately we can simultaneously force $b=c$ and also $0\le \mathrm{Arg}(b)<\pi$, and then $bb'=cc'$ also gives $b'=c'$.
Comparing $a^2=-bce$, $a=bcd$, and $ab=-bde$ give $e=-da$ and $d\in\{\pm 1\}$ (and likewise
$e'=-d'a'$ and $d'\in\{\pm 1\}$). {But we knew $dd'=1$, so we have} $d'=d$. Putting $aa'=(bcd)(b'c'd')=b^2b'{}^2$ into \eqref{aabb1}
gives $bb'\in\{(-1\pm\sqrt{5})/2\}$. In particular, $aa'$ cannot be 1, so we must have $a=a'$ and thus $b=\sqrt{da}$
and $b'=db$.}

{We eliminate} the possibility that $d=-1$ by considering the $st't'$ coefficient of $t'\rho^2(t)=\rho(t)t'$, which
gives $c=ab'd^2+cdee'$. So we have determined that $a=a'=-e=-e'=(-1\pm\sqrt{5})/2$, $b=b'=c=c'=\sqrt{a}$, and $d=d'=1$, where we can take the square-root for $b$ so that $b\in \bbR_{>0}\cup i\bbR_{>0}$.
So we have 2 possible solutions, corresponding to the choice of signs in
$a=(-1\pm\sqrt{5})/2$. {In section 4 we generalise this argument to arbitrary odd order abelian $G$
in \eqref{HIfus}.}

{Conversely, given either solution $a=(-1\pm\sqrt{5})/2$, we can define $\rho$ on the generators
$s,s',t,t'$ of $\cL_2$ by \eqref{YLrhos}-\eqref{YLrho}. Using
Corollary 1 below, this choice extends to} an algebra endomorphism $\rho$ on $\cL_2$ iff it respects the Leavitt--Cuntz relations: 
{i.e. $1=\rho(s')\rho(s)=\rho(t')\rho(t)=\rho(s)\rho(s')+\rho(t)\rho(t')$ and $0=\rho(s')\rho(t)=
\rho(t')\rho(s)$. It is straightforward to verify this (this is done in full generality in section 5). To show $\rho$
satisfies \eqref{ylrho2}, note that  both sides of \eqref{ylrho2} are
manifestly endomorphisms, so it suffices to verify it {for each of the four generators $x\in\{s,s',t,t'\}$}. If we can} show
$s'\rho^2(x)=xs'$ and $t'\rho^2(x)=\rho(x)t'$ for $x=s,t$ {(these must hold if \eqref{ylrho2} is to hold),
then $\rho^2(x)=(ss'+tt')\rho^2(x)$ shows \eqref{ylrho2} holds for $x=s,t$. Likewise, if} 
$\rho^2(y)s=sy$ and $\rho^2(y)t=t\rho(y)$ for $y=s',t'$, then \eqref{ylrho2} holds for $x=s',t'$.
Again, the details are given in full generality in section 5. Thus $\rho$ {defined by \eqref{YLrho}}
obeys the Yang--Lee fusions \eqref{YangLee}.
Finally, we can confirm that the endomorphism $\rho$
we have just constructed is indeed {simple, i.e. Hom$(\rho,\rho)=\bbC$ as well as Hom$(\mathrm{id},
\mathrm{id})=\bbC$ and Hom$(\rho,\mathrm{id})=\mathrm{Hom}(\mathrm{id},\rho)=0$} (this is done in full generality in Proposition 1
below).

Much more delicate is to associate a (strict) fusion category to {both of these} $\rho$. The biggest
challenge {here for arbitrary $G$} is to define arbitrary (but finite) sums of endomorphisms using the Leavitt algebra, {in the sense of the right-side
of \eqref{ylrho2}}.  We are lucky
here with the Yang--Lee: because its Leavitt algebra has $2\times 2$ generators, we can capture arbitrary sums ---
e.g. $\rho\oplus\rho^3\oplus\rho^5$ can be written $s\rho s'+ts\rho^3s't'+tt\rho^5t't'$, {to choose a random example}. {The resulting fusion category for the solution with  $a=(-1+\sqrt{5})/2$ is the unitary
category associated to e.g. {the integrable modules of the affine $G_2$ algebra} at level 1, whilst for  $a=(-1-\sqrt{5})/2$, we obtain the Yang-Lee fusion
category. These two fusion categories are inequivalent {even though they share the same fusions
\eqref{YangLee}} ---  indeed, {it can be shown that} the categorical dimension of $\rho$ (defined 
next section) is $1/a=(1\pm\sqrt{5})/2$,
so is positive in one and negative in the other. Nevertheless they are clearly related by the Galois automorphism
interchanging $a=(-1\pm\sqrt{5})/2$.}

To
realise the fusions \eqref{HIfus} {for general $G$}, we will need a Leavitt algebra $\cL$ with $(1+|G|)\times 2$ generators {(one pair for each term on the right of \eqref{HIfus}), but for such an algebra 
only direct sums with $n\equiv 1$ (mod $|G|$) terms can be realised}. When $\rho$ is a $*$-map
({e.g. the case studied in \cite{iz3},\cite{EG2}}), {we can extend} $\rho$ to an endomorphism of an infinite von
Neumann factor $N$ \cite{iz3}; semisimplicity  is then automatic, since $N$ contains copies of
the Leavitt algebras of arbitrary rank, {so arbitrary sums of  endomorphisms can be made. On the 
other hand, when $\rho$ is not a $*$-map, we obtain semisimplicity by first forming the idempotent completion.}   
 In section 5 we show that any solution to the various consistency
equations yields a (usually non-unitary) fusion category. 

 {MTC
structures can be placed on both of the $G=1$ fusion categories constructed in this section, though in more than one way --- e.g. 
{the $a=(-1+\sqrt{5})/2$ category is realised by both affine} $G_{2}$ level 1 
and affine $F_{4}$ level 1,  which {which are inequivalent as MTC since they} have  different central charges mod 8.}
%,  and their doubles both yield the modular tensor category for affine $G_2\oplus F_4$ at level (1,1).
{This behaviour too is special to $G=1$: the fusion categories for larger $G$ never come
with a braiding {(this is clear from \eqref{HIfus}, as $[\alpha_g][\rho]\ne[\rho][\alpha_g]$ when
$|G|$ is odd and $>1$)}. For these other $G$, we realise  in section 6}  the 
associated MTC through the centre of the tube algebra. {Incidentally, this construction applied
to e.g. the fusion category of affine $G_2$ at level 1, would yield the MTC of affine $G_2\oplus F_4$ at level
$(1,1)$.} 

{Although the fusion (or modular tensor) categories of Yang--Lee and affine $G_2$ or $F_4$ at level 1
are merely related by a Galois automorphism, the corresponding VOAs do not seem
related in any simple way. For example, the characters of Yang--Lee are
$$q^{11/60}(1+q^{2}+q^3+q^4+q^5+2q^6+\cdots)\,,\ \ 
q^{-1/60}(1+q+q^2+q^{3}+2q^4+2q^5+3q^6+\cdots)$$
with modular data
$$S=\frac{1}{\sqrt{5}}\left(\begin{matrix}-\sin(2\pi/5)&\sin(\pi/5)\\ \sin(\pi/5)&\sin(2\pi/5)\end{matrix}\right)\,,\
T=\left(\begin{matrix}e^{2\pi i 11/60}&0\\ 0&e^{-2\pi i /60}\end{matrix}\right)\,,$$
while those for affine $G_2$ at level 1 are
$$q^{-7/60}(1+14q+42q^{2}+140q^3+350q^4+840q^5%+1827q^6
+\cdots)\,,\ \ 
q^{17/60}(7+34q+119q^2+322q^{3}+819q^4+1862q^5%+4025q^6
+\cdots)$$
with modular data
$$S=\frac{1}{\sqrt{5}}\left(\begin{matrix}\sin(\pi/5)&\sin(2\pi/5)\\ \sin(2\pi/5)&-\sin(\pi/5)\end{matrix}\right)\,,\
T=\left(\begin{matrix}e^{-2\pi i 7/60}&0\\ 0&e^{2\pi i 17/60}\end{matrix}\right)$$
and those for affine $F_4$ at level 1 are
$$q^{-13/60}(1+52q+377q^{2}+1976q^3+7852q^4%+27404q^5+84981q^6
+\cdots)\,,\ \ 
q^{23/60}(26+299q+1702q^2+7475q^{3}+27300q^4%+88452q^5+260650q^6
+\cdots)$$
with modular data
$$S=\frac{1}{\sqrt{5}}\left(\begin{matrix}\sin(\pi/5)&\sin(2\pi/5)\\ \sin(2\pi/5)&-\sin(\pi/5)\end{matrix}\right)\,,\
T=\left(\begin{matrix}e^{-2\pi i 13/60}&0\\ 0&e^{2\pi i 23/60}\end{matrix}\right)\,.$$
In these cases, the first character given is that of the VOA {$\cV=\coprod_{n=0}^\infty\cV_n$} itself, and so lists the dimensions of its graded spaces $\cV_n$,
so we see that there appears little relation between the Yang--Lee VOA and that of say the $G_2$ one.
On the other hand, the naive inner product of the $G_2$ and $F_4$ character vectors is $j(\tau)^{1/3}$,
reflecting the fact that the VOA $\cV(G_2,1)\otimes \cV(F_4,1)$ is a conformal subalgebra of 
the $E_8$ lattice VOA. Note also that the first column of the matrix $S$ is strictly positive
for the VOAs $\cV(G_2,1)$ and $\cV(F_4,1)$ (as it must be for unitary VOAs), and isn't for the Yang--Lee (as is typical
for non-unitary VOAs).}

\section{Leavitt algebras and categories of endomorphisms}

For each $n>1$ define the \textit{Leavitt algebra} $\cL_n$ to be the associative $*$-algebra
freely generated over $\bbC$ by $x_1,\ldots,x_n,x_1',\ldots,x_n'$, modulo the \textit{Leavitt--Cuntz relations} 
\begin{equation} x'_ix_j=\delta_{i,j}\,,\ \ 
\sum_{i=1}^nx_ix'_i=1\,.\label{leavitt}\end{equation}
 The elements of $\cL_n$ are polynomials in the non-commuting
variables $x_i,x_j'$. The $*$-operation 
 sends  $x_i\mapsto x_i'$, $x_i'\mapsto x_i$, and obeys  $(cyz)'=\overline{c}z'y'$ for all
$c\in\bbC$ and $x,y\in\cL_n$. It has an obvious grading by $\bbZ^n$.
The Leavitt algebra $\cL_n$ can be regarded as the polynomial part of the Cuntz algebra $\cO_n$, {its $C^*$-algebra
completion}.

The Leavitt algebras $\cL_n$ are all non-isomorphic for $n=2,3,4, \ldots $, since the {inclusion of} $\cL_n $ in $\cO_n$ induces an isomorphism on $K$-theory with the cyclic group $\mathbb Z_{n-1}$ \cite{ABC}. The only obstruction to embedding $\cL_m$ unitally in $\cL_n$ is given by the $K$-theory
\cite{Phil1}. More precisely $\cL_m$ embeds  unitally in $\cL_n$ if and only if $m-1$ divides $n-1$. 
In the Cuntz framework of Izumi \cite{iz0, iz3, iz4} and Evans--Gannon \cite{EG2,EG4,EGn}, one constructs endomorphisms on a fixed
Cuntz algebra $\cO_n$, with prescribed fusion rules  and then extends these to a completion
as an infinite von Neumann factor $N$. Any Cuntz algebra $\cO_m$ can be unitally embedded in  the  factor $N$ 
for any $m$, even though usually it cannot be unitally  embedded in $\cO_n$. The fusion category will then be realised as 
a system of endomorphisms of $N$, since addition of any number $m$ of endomorphisms can be
expressed {in $N$}.   

We will realise fusion categories through endomorphisms of $\cL_n$. But we do not require that
our endomorphisms be $*$-maps, so they need not extend to the completion, the Cuntz algebra
 $\cO_n$ or the Banach algebras of Phillips \cite{Phil}. 
 
 Note that if $\rho$ is any algebra endomorphism on $\cL_n$,
then so is $\tilde{\rho}$ defined by 
{\begin{equation}\label{tilde}\tilde{\rho}(y)=\rho(y')'\,.\end{equation}}
\textit{Throughout this paper we distinguish an algebra endomorphism from a
$*$-algebra endomorphism.} The latter must obey $f(y)'=f(y')$ {(equivalently $\tilde{\rho}=\rho$)} while the former may not. 

%$\cL_n$ is graded: For any $n$-tuple $(k_1,\ldots,k_n)\in\bbZ^n$, let $\cL_n(k_1,\ldots,k_n)$
%be the span of all monomials where, for each $i$, the number of $x_i$ appearing in the monomial
%minus the number of $x_i'$, equals $k_i$. Then $\cL_n=\coprod_{(k_i)\in\bbZ^n}\cL_n(k_1,\ldots,k_n)$ and $\cL_n(k_1,\ldots,k_n)\cL_n(k_1',\ldots,k_n')\subset\cL_n(k_1+k_1',\ldots,k_n+k'_n)$. 

There is a canonical way to write any element of $\cL_n$. Call any monomial in the generators
$x_i,x_j'$ \textit{reduced} if no primed variable appears to the left of any unprimed variable,
and $x_1$ is not adjacent to $x_1'$ in the monomial. Call any  linear combination over $\bbC$ of finitely many distinct reduced monomials, a \textit{reduced sum}.

\medskip\noindent\textbf{Lemma 1.} \cite{Leav2} Any $y\in\cL_n$ can be written in one and only one way
as a reduced sum.\medskip

This simple observation has several easy consequences, as we'll see. {It easily implies
the centre of $\cL_n$ is trivial \cite{Leav2}. Moreover:}

\medskip\noindent\textbf{{Corollary 1.}}  An algebra  endomorphism $\rho$ on $\cL_n$ is uniquely defined by
its values $\rho(x_i),\rho(x_j')$ on the generators, and these can be assigned arbitrarily provided they
respect the Leavitt--Cuntz relations \eqref{leavitt}.

%(a)} If $y\in\cL_n$ obeys $yx_1=0$, then 
%$y$ is a reduced sum of  the form $\sum_{j>1}y_jx_j'$, where the $y_j$ are themselves reduced.

\medskip There are several {complications} caused by avoiding the completion and working
exclusively with $\cL_n$. {In particular}, two serious challenges are how to add {endomorphisms}, and how to
get rigidity.
We accomplish the former through the idempotent completion (described below), and the latter by hand. 

Recall that because $\cL_n$ is a {unital algebra over $\bbC$}, by general nonsense its algebra endomorphisms  define a $\bbC$-linear {preadditive strict} tensor category ${\cE\cN\cD}(\cL_n)$.
More precisely, the objects {in ${\cE\cN\cD}(\cL_n)$} are algebra (but not $*$-algebra) endomorphisms
of $\cL_n$. The morphisms  ${r}\in\mathrm{Hom}
(\beta,\gamma)$  are intertwiners,  i.e. $r\in\cL_n$ for which $r\beta(x)=\gamma(x)r$ for all $x\in\cL_n$; composition of morphisms is multiplication in $\cL_n$.
{$\cE\cN\cD(\cL_n)$ is $\bbC$-linear, i.e. each Hom$(\beta,\gamma)$ is a vector space over $\bbC$; it is also preadditive, i.e. composition of morphisms is bilinear.}
The tensor product of objects is composition: $\beta\otimes\gamma=\beta\circ\gamma$, whilst
of morphisms is: $r\otimes s=r\beta(s)=\gamma(s)r\in\mathrm{Hom}(\beta\circ\rho,\gamma\circ\sigma)$ when $r\in\mathrm{Hom}(\beta,\gamma)$, $s\in\mathrm{Hom}(\rho,\sigma)$.

A \textit{fusion category} {\cite{ENO}} is a $\bbC$-linear semisimple rigid {tensor} category with finitely many 
isomorphism classes of simple objects and finite dimensional spaces of morphisms, such that {the unit object 1 is simple}. 
{A simple object $X$ is one with End$(X)=\bbC \,\mathrm{id}_X$; amongst other things, every object in a
semisimple category is a direct sum of simple ones.}
 We say object $X$ has a {right-dual}  $X^\vee$ iff there is a pair of morphisms \textit{evaluation} $e_X\in\mathrm{Hom}(X^\vee\otimes X,1)$ and \textit{co-evaluation} $b_X\in\mathrm{Hom}
(1,X\otimes X^\vee)$ for which 
\begin{equation}(\mathrm{id}_X\otimes e_X)\circ(b_X\otimes \mathrm{id}_X)=\mathrm{id}_X\,,\ \
(e_X\otimes\mathrm{id}_{X^\vee})\circ(\mathrm{id}_{X^\vee}\otimes b_X)=\mathrm{id}_{X^\vee}\,
%\textcolor{red}{(e_X\otimes\mathrm{id}_X)\circ(\mathrm{id}_X\otimes b_X)=\mathrm{id}_X\,,\ \
%(\mathrm{id}_{X^\vee}\otimes e_X)\circ(b_X\otimes \mathrm{id}_{X^\vee})=\mathrm{id}_{X^\vee}\,.}
\end{equation}
{(where we assume the category is strict, for convenience). Left-dual $^\vee X$ is defined similarly.
In particular in $\mathcal{END}(\cL)$, an object $\beta\in\mathrm{End}(\cL)$ has a right-dual $\beta^\vee\in
\mathrm{End}(\cL)$ if there are elements $e_\beta\in\mathrm{Hom}(\beta\circ\beta^\vee,\mathrm{id})$ and $b_\beta
\in\mathrm{Hom}(\mathrm{id},\beta^\vee\circ\beta)$ in $\cL$ such that
\begin{equation}\beta(e_\beta)b_\beta=1=e_\beta\beta^\vee(b_\beta)\,.\label{coeval}\end{equation}}
A {tensor} category  is called \textit{rigid} if every object $X$ has a right- and left-dual. 

{In a (strict) rigid category, we can define the right-dual $f^\vee\in\mathrm{Hom}(Y^\vee,X^\vee)$ of a 
morphism $f\in\mathrm{Hom}(X,Y)$ by 
\begin{equation} f^\vee=(e_Y\otimes \mathrm{id}_{X^\vee})\circ(\mathrm{id}_{Y^\vee}\otimes f\otimes \mathrm{id}_{X^\vee})\circ(\mathrm{id}_{Y^\vee}\otimes b_X)\,.\end{equation}
In particular in $\mathcal{END}(\cL)$, the right-dual of $r\in\mathrm{Hom}(\alpha,\beta)$ is defined by
\begin{equation}\label{dualmorph}r^\vee=e_\beta\,\beta^\vee(r\,b_\alpha)\,.\end{equation}
 Then $(f\circ g)^\vee=g^\vee\circ f^\vee$ when the composition is defined. Left-dual $^\vee f$ is defined similarly. 
A rigid tensor category is \textit{pivotal} if it is equipped with a natural isomorphism from the identity functor to the double-dual functor $X\mapsto X^{\vee\vee}$.
In a pivotal category we can take ${}^\vee X=X^\vee$.} In a rigid category the (left-)dimension of object $X$ is
$e_{X^\vee} b_{X}$; {a semisimple pivotal category is called \textit{spherical}
if $X$ and $X^\vee$ have the same dimension for all objects $X$ (it suffices to {check this for} simple $X$). 
See e.g. \cite{m3} for the  remaining terminology not explained here.}

Let $\cE$ be a collection of algebra endomorphisms of $\cL_n$ closed under composition. We
require the identity to be in $\cE$. Let $\cC(\cE)$ denote the subcategory of ${\cE\cN\cD}(\cL_n)$ restricted to $\cE$. Then like ${\cE\cN\cD}(\cL_n)$, $\cC(\cE)$ is a $\bbC$-linear tensor category, and the endomorphism algebra of the unit object 1 is $\bbC$. 
By its \textit{idempotent completion} we mean the category $\overline{\cC(\cE)}$ whose objects consist of pairs
$(p,\beta)$ where $\beta\in \cE$ and $p\in \mathrm{End}(\beta)$ is an idempotent, i.e. $p^2=p$, and
whose morphism spaces are Hom$((p,\beta),(q,\gamma))=q\mathrm{Hom}(\beta,\gamma)p$ {with 
composition again given by multiplication}.
$\overline{\cC(\cE)}$ is a tensor category using $(p,\beta)\otimes(q,\gamma):=(p\otimes q,\beta\otimes\gamma)=(p\beta(q),\beta\circ\gamma)$, {and the tensor product of $qrp\in\mathrm{Hom}((p,\beta),(q,\gamma))$ with $q'r'p'\in\mathrm{Hom}((p',\beta'),(q',\gamma'))$ is $(qrp)\otimes(q'r'p')=qrp\beta(q'r'p')$}.
We can {introduce} direct sums into $\overline{\cC(\cE)}$ as follows. Objects {in this new category} consist of ordered
$n$-tuples $((p_1,\beta_1),\ldots,(p_n,\beta_n)) =(p_1,\beta_1)\oplus\cdots
\oplus(p_n,\beta_n)$, and {the morphism spaces are}
\begin{align}\mathrm{Hom}(((p_1,\beta_1),\ldots,(p_n,\beta_n)),&((q_1,\gamma_1),\ldots,(q_m,\gamma_m)))\nonumber\\ =&\,\left(\begin{matrix}q_1\mathrm{Hom}(\beta_1,\gamma_1)p_1&\cdots&
q_1\mathrm{Hom}(\beta_n,\gamma_1)p_n\\ \vdots&\ddots&\vdots\\
q_m\mathrm{Hom}(\beta_1,\gamma_m)p_1&\cdots&
q_m\mathrm{Hom}(\beta_n,\gamma_m)p_n\end{matrix}\right)\,.\label{hom}\end{align}
{Composition is matrix multiplication. Then $((p_1,\beta_1),\ldots,(p_n,\beta_n))\otimes((q_1,\gamma_1),\ldots,(q_m,\gamma_m))$ is the direct sum of $(p_i,\beta_i)\otimes(q_j,\gamma_j)$, while
$$\left(\begin{matrix}q_1r_{11}p_1&\cdots&q_1r_{n1}p_n\\ \vdots&\ddots&\vdots\\
q_mr_{1m}p_1&\cdots&q_mr_{nm}p_n\end{matrix}\right)\otimes\left(\begin{matrix}q'_1r'_{11}p'_1&\cdots&q'_1r'_{n'1}p'_{n'}\\ \vdots&\ddots&\vdots\\
q'_{m'}r'_{1m'}p'_1&\cdots&q'_{m'}r'_{n'm'}p'_{n'}\end{matrix}\right)$$
is the Kronecker product with $(ij,i'j')$-entry $q_ir_{ji}p_j\otimes q'_{i'}r'_{j'i'}p'_{j'}$.}
 We will write $\overline{\cC(\cE)}{}^{ds}$ for the {idempotent completion $\overline{\cC(\cE)}$
extended by direct sums in this way.}

%Define an equivalence relation on ring endomorphisms by $A\sim uAu^{-1}$ where $u$ is a unit in $\cL$. By a  \textit{sector} $[A]$ we mean an equivalence class of endomorphisms.

\medskip\noindent\textbf{Lemma 2.} {Let $\cE$ be a collection of $\cL_n$-endomorphisms as
above, and recall \eqref{tilde}.}  Suppose Hom$(\beta,\gamma)=\mathrm{Hom}(\tilde{\beta},\tilde{\gamma})$ {in $\cL_n$} for all
$\beta,\gamma\in\cE$, and {that} these are all finite-dimensional. Then  {$\overline{\cC(\cE)}{}^{ds}$, the idempotent completion extended by direct sums,} 
is a semisimple {strict} $\bbC$-linear tensor category with finite-dimensional hom-spaces. If {$\cC(\cE)$} is rigid, then so is $\overline{\cC(\cE)}{}^{ds}$.\medskip

\noindent\textit{Proof.} 
The  category $\overline{\cC(\cE)}{}^{ds}$ is manifestly $\bbC$-linear {and strict}. Since all Hom$(\beta,\gamma)$
are finite-dimensional, so are all Hom-spaces {\eqref{hom}} in $\overline{\cC(\cE)}{}^{ds}$.  
Since the anti-linear involution $x\mapsto x'$ sends Hom$(\tilde{\beta},\tilde{\gamma})$ to
Hom$(\gamma,\beta)$, and Hom$(\tilde{\beta},\tilde{\gamma})=
\mathrm{Hom}(\beta,\gamma)$ by hypothesis, then $x\mapsto x'$  bijectively maps 
Hom$({\beta},{\gamma})$ to
Hom$(\gamma,\beta)$. This implies that the (finite-dimensional) algebra
End$(((p_1,\beta_1),\ldots,(p_n,\beta_n)))$ is a $*$-algebra, and hence is  semisimple. Then Corollary 2.3 of \cite{Ya} tells us $\overline{\cC(\cE)}{}^{ds}$ is a semisimple category.

Moreover, suppose {$\cC(\cE)$} is rigid. Then Lemma 3.1 of \cite{Ya} says that  its idempotent
completion $\overline{\cC(\cE)}$ is also rigid: {e.g. $(p,\beta)^\vee=(p^\vee,\beta^\vee)$ where the dual morphism $p^\vee$  is defined in \eqref{dualmorph}, and (co-)evaluation 
is $e_{p,\beta}=p^\vee\beta^\vee(p)\,e_\beta$ and $b_{(p,\beta)}=p\beta(p^\vee)\,b_\beta$}. Hence {$\overline{\cC(\cE)}{}^{ds}$ is also rigid}:  take $((p_1,\beta_1),\ldots,(p_n,\beta_n))^\vee=((p_1,\beta_1)^\vee,\ldots,(p_n,\beta_n)^\vee)$ 
with diagonal (co-)evaluations $e_{(\ldots,(p_i,\beta_i),\ldots)}=\mathrm{diag}(e_{(p_i,\beta_i)})$ etc. QED\medskip

This condition Hom$(\beta,\gamma)=\mathrm{Hom}(\widetilde{\beta},\widetilde{\gamma})$
is crucial for extending the (unitary) Cuntz algebra methods to the (not necessarily unitary) Leavitt setting.
We show {near} the end of section 5 that this condition  holds for the Haagerup--Izumi systems considered here,
and {the same argument should work} for the near-group systems {constructed in} \cite{EG4}.
{Nevertheless, Lemma 2 emphasises} 
 that semisimplicity in the Leavitt picture is \textit{not} automatic, {and this} is very good:
it means {our} context should be flexible enough to include non-semisimple examples
such as those corresponding to the logarithmic CFTs discussed in the Introduction.

\section{Non-unitary Haagerup--Izumi: deconstruction}

 Let $G$ be any abelian group of odd order  $\nu=2n+1$, and define
 $\delta_\pm=(\nu\pm\sqrt{\nu^2+4})/2$,
 the two roots of $x^2=1+\nu x$. {Recall} the {Haagerup--Izumi fusions} \eqref{HIfus}.
% were studied in \cite{iz3} as endomorphisms of the Cuntz algebra. This method reduces the construction to solving a system
%of polynomial equations. \cite{EG2} 
%showed that there are many (conjecturally infinitely many) solutions to these equations, hence
%subfactors,  for larger $G$. Of course these results all required unitarity, i.e. that the endomorphisms
%are $*$-maps.
{A} main result (Theorem 1) of this paper {associates to any system of
algebra endomorphisms} realising these fusions, {a set of numerical invariants.}
 The converse, which {associates a system of endomorphisms and a fusion category
 to these same numerical invariants,} is given next section.

{Suppose $\alpha_g,\rho$ are algebra endomorphisms of an algebra $\cA$ which}
realise the Haagerup--Izumi fusions.  More precisely, {this means}  
{\begin{align}\alpha_g\circ\alpha_h=&\alpha_{g+h}\,,\ 
\alpha_g\circ\rho=\rho\circ\alpha_{-g}\,, \label{alpharho}\\
\rho(\rho(x))=&sxs'+\sum{}_g\,t_g\alpha_g(\rho(x))t_g'\,,\label{rho2}\end{align}}
{where $s,s',t_g,t_g'\in\cA$ satisfy  $s's=1$, $s't_g=t_g's=0$, $t_g't_h=\delta_{g,h}$,
and $1=ss'+\sum_gt_gt'_g$.  We do not assume $\cA$ is a $*$-algebra.} 
{Equation \eqref{alpharho} implies that each $\alpha_g$ is invertible. Note that we have the freedom to rescale the $\nu+1$ elements $s,t_g$ arbitrarily and
independently, provided we then rescale $s',t_g'$ inversely.}
We also require $\alpha_g$ and $\alpha_g\rho=\alpha_g\circ\rho$ to be simple, {i.e.}
{that their intertwiners in $\cA$ are} Hom$(\alpha_g\rho,\alpha_h\rho)=
\mathrm{Hom}(\alpha_g,\alpha_h)=\bbC\delta_{g,h}$ and Hom$(\alpha_g,\alpha_h\rho)=
\mathrm{Hom}(\alpha_g\rho,\alpha_h)=0$. {This implies for instance that $\cA$ has trivial centre, and that
{the representation} $g\mapsto\alpha_g$ of $G$ is faithful.}

{Unless $G$ is cyclic
(in which case $H^2_G(\mathrm{pt};\bbT)=1$), \eqref{alpharho} can be generalised by twisting by 2-cocycles
$\xi\in Z^2_G(\mathrm{pt};\bbT)$ and \eqref{HIfus} will still hold, as explained e.g. in the proof of Theorem 1 in \cite{EG4}}. We will
ignore this generalisation, as it is {conceptually straightforward and merely makes the arithmetic} a little messier, and our primary purpose
with this paper is to explain how to capture non-unitary fusion categories by endomorphisms. \cite{iz3} also
ignored these cocycles, but the unpublished notes \cite{iz4} introduces them (though of course in the unitary setting).

\medskip\noindent\textbf{Theorem 1.} Let $G,\alpha_g,\rho$ {and $s,s',t_g,t_g'\in\cA$} be as above. 
Then
\begin{align}\rho(s)&={\delta_\pm^{-1}}s+b\sum{}_{g}\,t_gt_g\,,\ \ \rho(s')={\delta_\pm^{-1}}s'+{\omega}b\sum{}_{g}\,t'_gt'_g\,,\label{rhos}\\
\rho(t_g)&= bst'_{-g}+{\omega}t_{-g}ss'+\sum{}_{h,k}\,A_{h+g,k+g}t_ht_{h+k+g}t_k'\,,\label{rhotg}\\
\rho(t_g')&= {\omega}bt_{-g}s'+{\overline{\omega}}ss't'_{-g}+\sum{}_{h,k}\,A_{k+g,h+g}t_kt'_{g+h+k}t_h'\,,\label{rhotgp}\\
\alpha_g(s)&\,=s\,,\ \ \alpha_g(s')=s'\,,\ \ 
\alpha_g(t_h)=t_{h+2g}\,,\ \ \alpha_g(t_h')=t'_{h+2g}\,,\label{alphast}
\end{align}
{for some fixed sign $\pm$,} where {$b\in\{1/\sqrt{\omega\delta_\pm}\}$ and $\omega^3=1$}. {In particular, $\alpha_g$ and $\rho$ restrict
to algebra endomorphisms of the Leavitt algebra  $\cL=\cL_{\nu+1}$  with
  generators $s,s',t_g,t_g'$. Moreover, 
 $A_{g,h}\in \bbC$ satisfy} 
\begin{align}
A_{g,h}&\,={\omega}A_{-h,g-h}={\overline{\omega}}A_{h-g,-g}\,,\label{order3sym}\\
\sum{}_h\,A_{h,0}&\,=-{\overline{\omega}\delta_\pm^{-1}}\,,\label{lin1}\\
 \sum{}_g\,A_{h+g,k}&A_{k,g}=\delta_{h,0}-{\delta_\pm^{-1}}\delta_{k,0}\,,\label{quad1}\\
{\sum{}_{l,m}\,A_{l,m}A_{l+g,h}}&{A_{h+m,l+i}A_{i,k+m}}\nonumber\\ &={A_{h-g,i-g}\delta_{k,g}-\overline{\omega}\delta_\pm^{-1}\delta_{h,0}A_{i,k}-{\omega}\delta_\pm^{-1}A_{g,h}\delta_{i,0}}\,.\label{quart}\end{align}

{We will show in Proposition 2 below that in fact
\begin{equation}\label{cubic} \overline{\omega}\sum{}_mA_{m,g+h}A_{g,m+k}A_{h,m+l}=A_{g+l,k}A_{h+k,l}-\delta^{-1}_\pm\delta_{g,0}\delta_{h,0}\,,\end{equation}
for all $g,h,k,l\in G$. We expect that this can be used to derive the more complicated \eqref{quart}, but we haven't
established this yet.}

{According to Izumi \cite{iz3}, a (unitary) Q-system corresponds to the special case of
Theorem 1 with $\omega=1$, $\delta_\pm=\delta_+$, $\overline{A_{g,h}}=A_{h,g}$, $A_{g,0}=\delta_{g,0}-1/(\delta_+-1)$. In this case the quartic identity \eqref{quart} can be replaced with the cubic identity \eqref{cubic}. This special case corresponds to fusion categories coming from one of the even subsystems
of a finite depth finite index subfactor.}

Incidentally, it doesn't matter which square-root is chosen for $b$ in Theorem 1: replacing
$s\mapsto -s$, $s'\mapsto-s'$ shows $b$ is equivalent to $-b$. {This means that} we can require without loss
of generality that $b$ lies on the positive halves of the real or imaginary axes. Which {triples $(\pm,\omega,A)$} yield
isomorphic fusion categories is answered below in Theorem 2, {as is the question of unitarity}.

%\textcolor{red}{We have no examples with $\omega\ne 1$, i.e. $\omega$ a primitive 3rd root of 1. We know there
%are none for $G=\bbZ_1$ and $\bbZ_3$. }

\medskip\noindent\textbf{Lemma 3.} Let $\rho$ be any algebra endomorphism on $\cA$ satisfying \eqref{alpharho}
and \eqref{rho2}{, and assume $\alpha_g$ and $\alpha_g\rho$ are all simple.} Then {Hom$(\alpha_g,\rho^2)=\bbC s\delta_{g,0}$}, Hom$(\alpha_g\rho,\rho^2)=\bbC t_g$, {Hom$(\rho^2,\alpha_g)=\bbC s'\delta_{g,0}$,} and 
Hom$(\rho^2,\alpha_g\rho)=\bbC t_g'$. Moreover, Hom$(\rho^2,\alpha_g\rho^2)=\bbC ss'\delta_{g,0}
+\mathrm{span}_h\{t_{h+g}t'_{h}\}$.\medskip

\noindent\textit{Proof.} Directly from \eqref{rho2} we find $\rho^2(x)s=sx$, $\rho^2(x)t_g=t_g\alpha_g\rho(x)$, $s'\rho^2(x)=xs'$,
and $t_g'\rho^2(x)=\alpha_g\rho(x)t_g'$. In other words, $s\in\mathrm{Hom}(\mathrm{id},\rho^2)$,
$t_g\in\mathrm{Hom}(\alpha_g\rho,\rho^2)$, $s'\in\mathrm{Hom}(\rho^2,\mathrm{id})$, and 
$t_g'\in\mathrm{Hom}(\rho^2,\alpha_g\rho)$.  

{Now suppose $r\in\mathrm{Hom}(\alpha_g,\rho^2)$. Then $s'\in\mathrm{Hom}(\rho^2,\mathrm{id})$
and $t'_h\in\mathrm{Hom}(\rho^2,\alpha_h\rho)$ immediately imply $s'r\in\mathrm{Hom}(\alpha_g,\mathrm{id})
=\bbC\delta_{g,0}$ and $t'_hr\in\mathrm{Hom}(\alpha_g,\alpha_h\rho)=0$ by simplicity. 
 Therefore $r=ss'r+\sum_ht_ht'_hr\in\bbC s\delta_{g,0}$, hence Hom$(\alpha_g,\rho^2)=\bbC s\delta_{g,0}$.}

{
Next, suppose $r\in \mathrm{Hom}(\rho^2,\alpha_g)$. Then $rs\in\mathrm{Hom}(\mathrm{id},\alpha_g)=\bbC\delta_{g,0}$
and $rt_h\in\mathrm{Hom}(\alpha_h\rho,\alpha_g)=0$, which forces $r\in \bbC s'\delta_{g,0}$ as before, and thus
Hom$(\rho^2,\alpha_g)=\bbC s'\delta_{g,0}$.}

{
Now consider $r\in\mathrm{Hom}(\alpha_g\rho,\rho^2)$. Then $s'r\in\mathrm{Hom}(\alpha_g\rho,\mathrm{id})=0$ and  $t'_hr\in\mathrm{Hom}
(\alpha_g\rho,\alpha_h\rho)=\bbC\delta_{g,h}$, and thus   Hom$(\alpha_g\rho,\rho^2)=\bbC t_g$.}
 
 {
 Similarly, let $r\in\mathrm{Hom}(\rho^2,\alpha_g\rho)$. Then $rs\in\mathrm{Hom}(\mathrm{id},\alpha_g\rho)=0$ and 
 $rt_h\in\mathrm{Hom}(\alpha_h\rho,\alpha_g\rho)=\delta_{h,g}\bbC$, which gives us  Hom$(\rho^2,\alpha_g\rho)=\bbC t_g'$. }
 
 {
 Finally, suppose $r\in\mathrm{Hom}(\rho^2,\alpha_g\rho^2)$. Then, using the invertibility of $\alpha$ and the calculation $\alpha_g\rho^2=\rho\alpha_{-g}\rho=\rho^2\alpha_g$, we get $rs\in\mathrm{Hom}(\mathrm{id},\alpha_g\rho^2)=\mathrm{Hom}(\alpha_g,\rho^2\alpha_g)=\bbC s\delta_{g,0}$. Similarly, $rt_h\in\mathrm{Hom}
 (\alpha_h\rho,\alpha_g\rho^2)=\mathrm{Hom}(\alpha_{h+g}\rho\alpha_g,\rho^2\alpha_g)=\bbC t_{h+g}$.
 This {suffices to identify} Hom$(\rho^2,\alpha_g\rho^2)$ in the usual way.} \qquad\textit{QED}\medskip

Note that because $\alpha_g$ is an algebra endomorphism and $s\in\mathrm{Hom}(\mathrm{id},\rho^2)$,
$\alpha_g(s)\in\mathrm{Hom}(\alpha_g,\alpha_g\rho^2)$. {But Hom$(\alpha_g,\alpha_g\rho^2)=
\mathrm{Hom}(\alpha_g,\rho^2\alpha_g)=\mathrm{Hom}(\mathrm{id},\rho^2)$ since $\alpha_g$ is  invertible.}
 By Lemma 3 this means $\alpha_g(s)=\psi(2g)s$ for some {$\psi(2g)\in \bbC$}
 (the 2 is introduced for later convenience; because the order
of $G$ is odd, 2 is invertible). {Because $\alpha_g\alpha_h=\alpha_{g+h}$, we see $\psi\in\widehat{G}$.}
From the Leavitt--Cuntz relation ${s's=1}$, we obtain $\alpha_g(s')=\psi(-2g)s'$. Likewise,
{$\alpha_h(t_g)\in\mathrm{Hom}(\alpha_{h+g}\rho,\alpha_h\rho^2)=\mathrm{Hom}(\alpha_{g+2h}\rho,\rho^2)
=\bbC t_{g+2h}$, and hence}
 $$\alpha_h(t_g)=\epsilon_h(g)t_{g+2h}$$
for some $\epsilon_h(g)\in{\bbC}$. {Again, $\alpha_g\alpha_h=\alpha_{g+h}$ implies  
these numbers $\epsilon_h(g)$ are non-zero and satisfy} 
\begin{equation}\epsilon_{h+k}(g)=\epsilon_h(g)\epsilon_k(g+2h)\,.\label{epsilon}\end{equation}
We can rescale $t_1,\ldots,t_{\nu-1}$ so that $\epsilon_h(0)=1$
for all $h$. But {from \eqref{epsilon} with $g=0$ this implies} $\epsilon_k(2h)=1$ for all $h,k\in G$, and invertibility of 2 then implies all $\epsilon_k(h)=1$. From  $t_g't_g=1$ we likewise get 
$ \alpha_h(t_g')=t'_{g+2h}$. {Thus we know all $\alpha_g$ restrict to endomorphisms of the Leavitt algebra
$\cL_{\nu+1}$ generated by the $s,s',t_g,t'_g$.}

{Since $s\in\mathrm{Hom}(\mathrm{id},\rho^2)$ and $\rho$ is an endomorphism, $\rho(s)\in\mathrm{Hom}
(\rho,\rho^3)$. Hence $s'\rho(s)\in\mathrm{Hom}(\rho,\rho)=\bbC$ and $t'_0\rho(s)\in\bbC t_0$. Write $s'\rho(s)
=a$ and $t_0'\rho(s)=bt_0$ for some $a,b\in\bbC$. Hitting the latter equation with $\alpha_h$, we get
$t'_{2h}\alpha_h(\rho(s))=bt_{2h}$, i.e. $t_g'\rho(s)=\psi(2g)bt_g$. Likewise, $\rho(s')s=a'$ and $\rho(s')t_g=b't'_g$ for some $a',b'\in\bbC$. We thus obtain from $\rho(s)=ss'\rho(s)+\sum_gt_gt_g'\rho(s)$ that}
\begin{equation}
\rho(s)=as+b\sum{}_{g}\,\psi(g)t_gt_g\,,\ \ \rho(s')=a's'+b'\sum{}_{g}\psi(-g)\,t'_gt'_g\,.\label{rhos1st}\end{equation}

{The computation of  $\rho(t_g)$ is similar. First note that $\rho(t_0)\in\mathrm{Hom}(\rho^2,\rho^3)$, 
so    $s'\rho(t_0)\in\mathrm{Hom}(\rho^2,\rho)=\bbC t'_0$ and  $t'_h\rho(t_0)\in\mathrm{Hom}(\rho^2,\alpha_{h}\rho^2)=\mathrm{span}\{\delta_{h,0}ss',t_kt'_{k-h}\}$, using Lemma 3. Write $s'\rho(t_0)=ct_0'$ and
$t_h'\rho(t_0)= \delta_{h,0}dss'+\sum_kA_{h,k}t_{h+k}t'_{k}$, for complex numbers $c,d,A_{h,k}$. Then $\rho(t_0)=cst'_0+dt_0ss'+\sum_{h,k}A_{h,k}t_ht_{h+k}t'_k$. The calculation for $\rho(t_0')t_h$ is identical, and involves
complex numbers $c',d',A'_{h,k}$. Hitting these with $\alpha_{-g/2}$ yields}
\begin{align}\rho(t_g)&\,=\psi(-g)cst'_{-g}+dt_{-g}ss'+\sum{}_{h,k}\,A_{h+g,k+g}t_ht_{g+h+k}t'_k\,,
\label{rhot}\\
\rho(t_g')&\,=\psi(g)c't_{-g}s'+d'ss't'_{-g}+\sum{}_{h,k}\,A'_{h+g,k+g}t_kt'_{g+h+k}t'_h\,.
\label{rhotp}\end{align}
 {Thus we also know  $\rho$ restricts to an endomorphism of the Leavitt algebra
$\cL_{\nu+1}$ generated by the $s,s',t_g,t'_g$.}

Thus the $\cA$-endomorphism $\rho$ is determined from the $2\nu^2+8$ parameters $a,a',b,b',c,c',d,d',A_{h,k},
A'_{h,k}$, as well as the character $\psi\in\widehat{G}$. However there are several consistency 
conditions, {coming from \eqref{rho2} and also the fact that $\rho$ being an endomorphism must preserve the Leavitt--Cuntz relations. To compute various expressions in $\cL_{\nu+1}$, it is convenient to collect our equations
\begin{align}&\ s'\rho(s)=a\,,\ t'_g\rho(s)=b\psi(g)t_g\,,\ \label{strhos}\\
s'\rho(t_g)=&\psi(-g)ct'_{-g}\,,\ t'_g\rho(t_h)=d\delta_{g,-h}ss'+\sum{}_kA_{g+h,k+h}
t_{g+h+k}t_k'\,,\label{strhot}\\
s'\rho(t_g')=&d's't'_{-g}\,,\ t'_h\rho(t_g')=\psi(g)c'\delta_{g,-h}s'+\sum{}_kA'_{k+g,h+g}
t'_{g+h+k}t'_{k}\,.\label{rhotpst}\end{align}
Implicit in the following is Lemma 1, which permits us to compare corresponding coefficients of an expression
in $\cL_{\nu+1}$ in  reduced form (i.e. replace any occurrence of $ss'$ with $1-\sum_gt_gt'_g)$.}

Because $\rho$ satisfies  \eqref{rho2}, we must have $s'\rho(\rho(x))=xs'$.
{But if instead we compute $s'\rho(\rho(s))=ss'$ directly from \eqref{rhos1st} and \eqref{rhot}, using  \eqref{strhos} and \eqref{strhot}, we obtain}
$$s'\rho(\rho(s))=a^2+bc\left(\nu dss'+\sum{}_g\,A_{0,g}\sum{}_{k}\,t_kt'_k\right)\,.$$
{Comparing these expressions for $s'\rho^2(s)$, and performing the analogous calculation for $\rho(\rho(s'))s=ss'$,} we obtain 
\begin{equation}bc\sum{}_g\,A_{0,g}=-a^2=\nu bcd-1\,,\ \ b'c'\sum{}_g\,A'_{0,g}=-a'{}^2=\nu b'c'd'-1\,.\label{rho2ss}\end{equation}
{Likewise,  the} $st_0'$ coefficient of  $t'_0\rho(\rho(s))=\rho(s)t'_0$ becomes
$a=bcd$ ({and similarly we get} $a'=b'c'd'$). {Substituting this into  \eqref{rho2ss}}, we obtain
$-a^2=\nu a-1$ and so $a\in\{1/\delta_\pm\}$ (similarly for $a'$). 

{In particular, $a,a'\ne 0$, so also $b,b',c,c',d,d'\ne 0$.  Hitting {$a=s'\rho(s)$} with $\alpha_g$, we obtain
$$a=\alpha_g(s')\rho(\alpha_{-g}(s))=\psi(-2g)s'\rho(\psi(-2g)s)=\psi(-4g)a$$
for all $g\in G$. Thus, since the order $\nu$ of $G$ is odd, we have that $\psi$ is identically 1. We thus recover \eqref{alphast}.}

Other {coefficients of} $t'_0\rho(\rho(s))=\rho(s)t'_0$ we need now give
{\begin{align}
d\sum{}_h\,A_{h,0}=&\,-a\,,\ \ d'\sum{}_h\,A'_{h,0}=-a'\,,\label{rho2ts1}\\
\sum{}_h\,A_{h,k}A_{k,k'+h}&-d\delta_{k,0}\sum{}_h\,A_{h,0}=\delta_{k',0}\,.\label{AA}
%\\ \sum{}_h\,A'_{h,k}A'_{k,k'+h}&-d'\delta_{k,0}\sum{}_h\,A'_{h,0}=\delta_{k',0}\,.\label{rho2ts2}
\end{align}}
From $\rho(s')\rho(s)=1$ we obtain 
$1=(a's'+b'\sum_gt'_gt'_g)(as+b\sum_gt_gt_g)$, i.e.
\begin{equation}aa'+\nu bb'=1\,.\label{sps}\end{equation}
%Next, expanding $0=\rho(s')\rho(t_g)$ we see that only the $t'_{-g}$ term survives, and so
%$ca'=-b'\sum_hA_{h+g,0}\psi(g-h)$. Together with the condition for $\psi\equiv 1$, this is equivalent to  
%\begin{equation}b'\sum_hA_{h,0}\psi(-h)=-ca'\,.\label{sptg}\end{equation}
%Likewise, $\rho(t'_g)\rho(s)=0$ gives 
 %\begin{equation}b\sum_hA'_{h,0}\psi(h)=-ac'\,.\label{tgps}\end{equation}
%The conditions $\delta_{l,m}=\rho(t'_l)\rho(t_m)$ are equivalent to
%\begin{equation}
%dd'=1\,,\ \ \sum_hA'_{h+m,k}A_{h,k}=\delta_{0,m}-cc'\psi(-m)\delta_{k,0}\,.\label{tlptm}\end{equation}
The $st'_gt'_g$, $t_gt_gs'$, constants, $t_gt'_g$, and $t_ht_{h+k}t'_{l+k}t'_l$ terms {of the 
 Leavitt--Cuntz relation  
$1=\rho(s)\rho(s')+\sum_g\rho(t_g)\rho(t_g')$} give respectively 
\begin{align} c\sum{}_h\,A'_{h,0}&\,=-ab'\,, \ \ c'\sum{}_h\,A_{h,0}=-ba'\,,\label{cuntz3a}\\
aa'+cc'\nu=&\,1\,,\ \ dd'=1\,,\label{cuntz3}\\ \sum{}_g\,A_{h+g,k}A'_{g,k}=&\,\delta_{h,0}-bb'\delta_{k,0}\,.\label{AAprime}\end{align}

Note that we still have the freedom 
to rescale $s\mapsto \lambda s$ and $s'\mapsto s/\lambda$;  choose $\lambda$
so that {$c=b$. Then $bb'=cc'$ (obtained by comparing \eqref{sps} with \eqref{cuntz3}) implies
$b'=c'$. Now, $aa'=(bcd)(b'c'd')=(bb')^2$, so \eqref{sps} implies $bb'\in\{1/\delta_\pm\}$. However, if
$a\ne a'$, then $aa'=1/(\delta_+\delta_-)=-1$, contradicting our value for $bb'$. Thus $a=a'=bb'$. Moreover,
comparing \eqref{cuntz3a} and \eqref{rho2ts1} gives $b'=bd$.}

{Multiplying \eqref{AAprime} by $A_{k,h+m}$ and summing over $h$ using \eqref{AA}} gives
{\begin{equation} A'_{g,h}+\delta_{k,0}d\left(\sum{}_gA'_{g,0}\right)\left(\sum{}_hA_{h,0}\right)
=A_{h,g}-\delta_{k,0}bb'\sum{}_hA_{0,h}\,.\label{ApA}\end{equation}
But the terms proportional to $\delta_{k,0}$ are $d\,(-a'/d')(-a/d)=da^2$ and $-bb'\,(-a^2/bc)=a^2b'/b$, which
we now know are equal. Thus $A'_{g,h}=A_{h,g}$ for all $g,h\in G$. }

{The $st_{h-g}'t_h'$ coefficient of $t_0'\rho^2(t_g)=\rho(t_g)t'_0$ is
\begin{equation}\label{AAsk}c\delta_{h,0}=d^2ab'\delta_{0,g}+cd\sum{}_kA_{g,k+g}A'_{h+k,-g}\,.\end{equation}
Multiplying \eqref{AAsk} by $A_{h+l,-g}$ and summing over $h$ using \eqref{AAprime} collapses to
$A_{l,-g}=dA_{g,l+g}$, which recovers  \eqref{order3sym}; because the permutation $(l,-g)\mapsto(g,l+g)$ is order 3, $d$
must be a 3rd root $\omega$ of 1.}

We  obtain {\eqref{lin1} and \eqref{AA} from \eqref{rho2ts1} and \eqref{quad1}.}
Finally, \eqref{quart} arises from the {$t_ht_{i-g+h}t'_{k+i}t'_k$
coefficient of $t_0'\rho(\rho(t_{-g}))=\rho(t_{-g})t'_{0}$.}
This completes our derivation of Theorem 1.

\section{Non-unitary Haagerup--Izumi: reconstruction}

{This section is devoted to a proof of the following theorem, another main result of our paper.
Recall $\delta_\pm=(\nu\pm\sqrt{\nu^2+4})/2$.}

\medskip\noindent\textbf{Theorem 2.} {Choose any finite abelian group $G$ of odd order $\nu$.}

% and let $\cL$ be the Leavitt algebra with generators $s,s'$ and $t_g,t'_g$ for all $g\in G$.
\smallskip\noindent{\textbf{(a)}} Let {$b\in\{1/\sqrt{\omega\delta_\pm}\}$ and $\omega^3=1$},
and choose any solution $A_{g,h}$ to \eqref{order3sym}--\eqref{quart}.
Define the values of $\rho$ and $\alpha_g$ on the generators $s,s',t_g,t_g'$ by {\eqref{rhos}--\eqref{alphast}}. Then these  extend to algebra endomorphisms $\rho,\alpha_g$  on the 
Leavitt algebra $\cL$ generated by $s,s',t_g,t_g'$. {Then $\overline{\cC(\{\alpha_g\rho^n\})}^{ds}$, 
the idempotent completion extended by direct sums as described in section 3, is} a strict spherical fusion category we'll denote by $\cC(G;\pm,\omega,A)$. {The simple objects of this category are} $\alpha_g=(1,\alpha_g)$ and $\alpha_g\rho=(1,\alpha_g\rho)$ 
 up to equivalence, and they satisfy
the Haagerup--Izumi fusions \eqref{HIfus}. The categorical dimensions of {$\alpha_g$} are 1 and of {$\alpha_g\rho$} are $\delta_\pm$. 

\smallskip\noindent{\textbf{(b)} Two such fusion categories $\cC(G^{(i)};\pm^{(i)},\omega^{(i)},A^{(i)})$
 are equivalent as {tensor} categories iff {$\pm^{(1)}=\pm^{(2)}$, $\omega^{(1)}=\omega^{(2)}$} and there is a group isomorphism
$\pi:G^{(1)}\rightarrow G^{(2)}$  such that $A^{(1)}_{g,h}
=A^{(2)}_{\pi g,\pi h}$ for all $g,h\in G^{(1)}$.}

\smallskip\noindent{\textbf{(c)} $\cC(G;\pm,\omega,A)$ is unitary iff {$\pm=+$} and
$A$ is a hermitian matrix: $A_{g,h}=\overline{A_{h,g}}$ for all $g,h\in G$. $\cC(G;\pm,\omega,A)$ is hermitian iff 
$A_{g,h}$ is hermitian.}\medskip

{We will learn below that the simple objects are all of the form $(uu',\alpha_g\rho^n)$
or $(vv',\alpha_g\rho^n)$ for certain monomials $u=u^{g,n}_{h,i},v=v^{g,n}_{h,j}$ recursively constructed below. The
modular data $S,T$ associated to the double of $\cC(G;\pm,\omega,A)$ is computed next section.} 

%{Let $\cL$ be the Leavitt algebra} {defined in the theorem.} 
%Choose endomorphisms $s,s',t_g,t_g'$ such that Hom$(\mathrm{id},\rho^2)=\bbC s$,
%Hom$(\rho^2,\mathrm{id})=\bbC s'$, and Hom$(\alpha_g\rho,\rho^2)=\bbC t_g$,
%Hom$(\rho^2,\alpha_g\rho)=\bbC t'_g$ for each $g\in G$. 
{By Corollary 1,} it is trivial that the $\alpha_g$ defined by \eqref{alphast} {are algebra endomorphisms of $\cL$.}
 {Similarly, to show that $\rho$ satisfying \eqref{rhos}}--\eqref{rhotgp} extends to an {algebra} endomorphism of $\cL$, {it suffices to verify that} the values of $\rho(s)$ etc  preserve the   Leavitt--Cuntz relations. {It is readily verified that these all reduce to  the identities
$b^4+\nu \omega b^2=1$, \eqref{lin1}, \eqref{quad1}, and
\begin{equation}\label{lin2}\sum{}_g\,A_{0,g}=-{\omega}\delta_{\pm}^{-1}\end{equation}
(the latter follows from \eqref{lin1} and \eqref{order3sym}).}
Thus $\rho$ is an
algebra endomorphism.

To verify  that $\alpha_g\rho=\rho\alpha_{-g}$, we need to show that $\alpha_g(\rho(x))=\rho(\alpha_{-g}(x))$ for $x=s,s',t_h,t_h'$. This is trivial to verify: e.g.
\begin{equation} 
\alpha_g(\rho(t'_l))= bt_{-l+2g}s'+ss't'_{-l+2g}+\sum_{h,k}A_{k+l,h+l}t_{k+2g}t'_{l+h+k+2g}t_{h+2g}'
=\rho(t'_{l-2g})\,.\end{equation}

To see that $\rho$  satisfies \eqref{rho2}, it suffices to verify that
$s'\rho(\rho(x))=xs'$, $t_g'\rho(\rho(x))=\alpha_g\rho(x)t_g'$, $\rho(\rho(y))s= sy$ and
$\rho(\rho(y))t_g=t_g\alpha_g\rho(y)$ for all $g\in G$, $x\in\{s,t_h\}$ and 
$y\in\{s',t_h'\}$. {This is because those equations  imply using $\rho^2(x)=
(ss'+\sum_gt_gt'_g)\rho^2(x)=\rho^2(x)(ss'+\sum_gt_gt'_g)$ that \eqref{rho2} holds when $x$ is any generator,
and this suffices to prove \eqref{rho2} for all $x$ because both sides of \eqref{rho2} are manifestly endomorphisms.} 
In fact, by $\alpha_g$-equivariance, it suffices to establish these for $g=0$. {All of these
equations reduce to $b^4+\nu \omega b^2=1$, \eqref{lin1}, \eqref{quad1}, and \eqref{lin2}, except for the following.}

The equation $s'\rho(\rho(t_g))=t_{g}s'$ {yields} the equations
\begin{align}1=&{\,2\overline{\omega}b^4+\omega b^2\sum{}_{h,k}\,A_{h,k}A_{k,h}}\,,\\
{-\omega b^2A_{h,k}-\overline{\omega}b^2\delta_{h,0}=}&\,{\sum{}_{\ell,m}\,A_{\ell,m}A_{m,\ell+h}A_{h,k+m}\,.}\end{align}
{The former follows from $\sum_{h,k}A_{h,k}A_{k,h}=\nu-\nu b^2$, which in turn follows from \eqref{quad1}. The latter follows directly from \eqref{quad1}.}
The equation $t'_0\rho(\rho(t_g))=\rho(t_{g})t'_0$ {gives} \eqref{quart}
as well as
\begin{align}{\sum{}_{l,m}\,A_{l,m}A_{l+g,k}A_{k+m,l}=}&\,{-b^2\delta_{k,0}
-\omega b^2A_{g,k}\,,}\\
\sum{}_k\,A_{k+g,h}A_{k,-h}=&\,{\omega \delta_{h,g}-\overline{\omega}b^2\delta_{h,0}\,,}\label{AAskew}\\
\sum{}_m\,A_{g,m+g}A_{-g,m+k}=&\,{\overline{\omega}\delta_{k,0}-}b^2\delta_{0,g}\,,\end{align}
which follow from \eqref{quad1} {and \eqref{order3sym}}. 

The simplicity of $\rho$ etc is established by the
following proposition.

\medskip\noindent\textbf{Proposition 1.} Let $\rho$ be as above.  Then for each $g,h\in G$, Hom$(\alpha_g\rho,\alpha_h\rho)=
\mathrm{Hom}(\alpha_g,\alpha_h)=\bbC\delta_{g,h}$ and Hom$(\alpha_g,\alpha_h\rho)=
\mathrm{Hom}(\alpha_g\rho,\alpha_h)=0$.

\medskip\noindent\textit{Proof.} {Write  $a=\delta^{-1}_\pm$. Choose any $x\in\cL$ commuting with 
$\rho(s)=as+b\sum_gt_gt_g$. We will begin by proving that such an $x$}
  must be a polynomial in $\rho(s)$.
  Write $x$ in reduced form (recall Lemma 1). We can assume
without loss of generality that no term in $x$ is a scalar times a power of $s$, i.e. $cs^l$, since 
otherwise we could replace $x$ with  $x-c(\rho(s)/a)^l$
 (the result will still lie in $\cL$ and commute with $\rho(s)$, and will be in $\bbC[\rho(s)]$ iff $x$ is). {Suppose for contradiction that
 $x\ne 0$.}
 
 Assume {first} that not all terms in $x$ begin with $s'$. Amongst those terms, let $w=s^lw'\ne 0$
be the sum of all terms with the maximal leading string of $s$'s ($l$ may be 0). Then $\rho(s)x$ contains
the terms $asw=as^{l+1}w'$, and these are reduced and have longer leading strings of $s$'s than any other terms
in $\rho(s)x-x\rho(s)$ (since no term can be a pure power of $s$). Being reduced, these terms $asw$ cannot cancel
{anything, contradicting} $\rho(s)x=x\rho(s)$.

It remains to consider $x=s'x'$. Then every term in $x'\ne 0$ involves  {only} 
$s'$'s and $t_g'$'s (since $x$ is reduced). Then
$\rho(s)x-x\rho(s)$ when reduced contains terms $-{a}\sum_ht_ht'_h{s'}x'$ with leading {factors $t_ht_h'$}. 
Again, these terms cannot cancel, which contradicts $\rho(s)x=x\rho(s)$.

{These contradictions mean $x=0$.}
Thus any $x\in\cL$ commuting with $\rho(s)$ must be a polynomial in $\rho(s)$, {and hence can contain no
$s',t_k'$.} Likewise,
any $x\in\cL$ commuting with $\rho(s')$ must be a polynomial in $\rho(s')$, and thus contains no
$s,t_k$. Together, they tell us that any $x$ commuting with both $\rho(s)$ and $\rho(s')$  must be a scalar.

Now suppose $x\alpha_g\rho(y)=\alpha_h\rho(y)x$ for all $y$. Then taking $y=s$ tells us
$x\rho(s)=\rho(s)x$, since $\alpha_g\rho(s)=\rho(\alpha_{-g}s)=\rho(s)$, while taking $y=s'$ tells us   
$x\rho(s')=\rho(s')x$. {Therefore $x\in\mathrm{Hom}(\alpha_g\rho,\alpha_h\rho)$ must again be a scalar $\lambda\in\bbC$. Now, for $\lambda\ne 0$, $\lambda \alpha_g\rho(t_0)=\alpha_h\rho(t_0)\lambda$ iff $\lambda\rho(t_{-2g})=
\lambda\rho(t_{-2h})$, iff $g=h$ (since $b\ne0$). Thus Hom$(\alpha_g\rho,\alpha_h\rho)=\delta_{g,h}\bbC$.}

Now turn to ${x}\in \mathrm{Hom}(\alpha_g,\alpha_h)$, i.e. $x\alpha_g({y})=\alpha_h(y)x$ for
all $y\in\cL$. In particular, $xs=sx$ and $xs'=s'x$. {By the identical argument as above, the former requires $x\in\bbC[s]$ 
while the latter requires $x\in\bbC[s']$,} 
 and thus $x$ is a scalar $\lambda\in\bbC$. Of course, $\lambda\ne 0$ intertwines $\alpha_g$
and $\alpha_h$ iff $g=h$, {by evaluating at $y=t_0$. Hence Hom$(\alpha_g,\alpha_h)$}.

Finally, suppose {$x\in\mathrm{Hom}(\alpha_g\rho,\alpha_h)$ and $x\ne 0$ is reduced. Then e.g. $x\rho(s)=sx$.} Assume first that at least one term in $x$
does not begin with $s'$. Amongst those terms, let $y$ be one with a maximal string of leading $s$'s
(this string may be empty, if no term in $x$ begins with $s$). Then $sy$ will be a reduced
term in $sx$, and the only reduced terms in  $x\rho(s)$ with a {leading string of
$s$'s of similar length} are those which are pure monomials in $s$. So $y=rs^n$ for some $n\ge 0$ and some
non-zero scalar $r$. But even those $y$ won't work: the reduced terms in $sx-x\rho(s)$
corresponding to $y$ are $rs^{n+1}-ars^{n+1}$, which can never vanish because $a\ne 1$.
If instead all terms in $x$ begin with an $s'$, then none of them end with an $s$, so
repeat this argument with {$x\rho(s')=s'x$. The proof that    Hom$( \alpha_g,\alpha_h\rho)=0$
is identical.}
\qquad\textit{QED to Prop.1}\medskip

{Recall the category $\mathcal{END}(\cL)$ defined in section 3. Let}
 $\cE$ consist of all monomials of the form $\alpha_g\rho^n$. {Since $(\alpha_g\rho^m)(\alpha_h\rho^n)=
 \alpha_{g\pm h}\rho^{m+n}$, the set $\cE$ is closed under composition. Let $\cC(\cE)$ be the subcategory of $\mathcal{
 END}(\cL)$ with objects $\alpha_g\rho^n$. We want to show $\cC(\cE)$ is rigid.}
Define $(\alpha_g\rho^{2k+1})^\vee=\alpha_g\rho^{2k+1}$ and 
$(\alpha_g\rho^{2k})^\vee=\alpha_{-g}\rho^{2k}$. {Then $(\alpha_g\rho^n)^\vee(\alpha_g\rho^n)=(\alpha_g\rho^n)(\alpha_g\rho^n)^\vee=
\rho^{2n}$ for all $g\in G,n\ge 0$. Define}
$e_{\alpha_g\rho^{n}}={\omega^n}b^{-n}s'\rho(s')\cdots \rho^{n-1}(s')$ and  $b_{\alpha_g\rho^{n}}={\omega^n}b^{-n}\rho^{n-1}(s)\cdots \rho(s)s$. {Since $s'\in\mathrm{Hom}(\rho^{k+2},\rho^k)$ (this is a special case of $s'\rho^2(x)=xs'$),
$\rho^m(s')\in\mathrm{Hom}(\rho^{n+2},\rho^n)$ for any $m\le n$ follows because $\rho$ is an endomorphism. Therefore,
$e_{\alpha_g\rho^n}\in\mathrm{Hom}(\rho^{2n},\mathrm{id})$ as required. Likewise, $b_{\alpha_g\rho^n}\in\mathrm{Hom}(\mathrm{id},\rho^{2n})$.
To see that $e_{\alpha_g\rho^n},b_{\alpha_g\rho^n}$} satisfy \eqref{coeval}, first note that for any $k\ge l$,
\begin{equation}\rho^k(s')\rho^l(s)=\rho^l(\rho^{k-l}(s')s)=\left\{\begin{matrix} 1&k=l\\
{\omega}b^2&k=l+1\\ \rho^l(s)\rho^{k-2}(s')&k\ge l+2\end{matrix}\right.\,.\end{equation}
Using this, it is easy to see that for any $n\ge 2$, we have
 \begin{equation}\rho^n(s')\rho^{n+1}(s')\cdots\rho^{2n-1}(s')\rho^{n-1}(s)\cdots\rho(s)s={\omega}
 b^2\rho^{n-1}(s')\rho^{n}(s')\cdots\rho^{2n-3}(s')\rho^{n-2}(s)\cdots\rho(s)s\,,\nonumber\end{equation}
which by an easy induction on $n$ gives the first  equation of \eqref{coeval}. The second equation
in \eqref{coeval} is handled analogously. {Thus $\cC(\cE)$ is rigid, with (co)evaluations $e,b$.}

{We want to apply Lemma 2. That means we must verify first that Hom$(\alpha_g\rho^m,\alpha_h
\rho^n)=\mathrm{Hom}(\widetilde{\alpha_g\rho^m},\widetilde{\alpha_h\rho^n})$ in $\mathcal{END}(\cL)$, where 
$\tilde{\beta}(x)=\beta(x')'$ is defined by \eqref{tilde}}. Note that $\tilde{\alpha_g}=\alpha_g$ (i.e.
$\alpha_g$ is a $*$-map), but $\tilde{\rho}$ is defined {by \eqref{rhos}--\eqref{rhotgp}} using the adjoint $\overline{A}_{h,g}$
in place of $A_{g,h}$, {$\overline{b}$ in place 
of $b$, and $\overline{\omega}$ in place of $\omega$. It is manifest that $\tilde{\rho}$ is an endomorphism
of $\cL$ satisfying \eqref{rho2}.}

We have $\widetilde{\alpha_g\rho^n}=\alpha_g\tilde{\rho}^n$. An easy induction from \eqref{rho2} (replacing $x$ there with $\rho^{n-2}(x)$ and hitting with $\alpha_g$) verifies
\begin{align} \label{arndecomp}
\alpha_g\rho^n(x)=&\,\sum{}_{h,i}\, u_{h,i}^{g,n}\alpha_h(x)u_{h,i}^{g,n\,\prime}
+\sum{}_{k,j}\,v_{k,j}^{g,n}\alpha_k\rho(x)v_{k,j}^{g,n\,\prime}\,,\\
\alpha_g\tilde{\rho}^n(x)=&\,\sum{}_{h,i}\, u_{h,i}^{g,n}\alpha_h(x)u_{h,i}^{g,n\,\prime}
+\sum{}_{k,j}\,v_{k,j}^{g,n}\alpha_k\tilde{\rho}(x)v_{k,j}^{g,n\,\prime}\,,\end{align}
where $u_{h,i}^{g,n}{\in\mathrm{Hom}(\alpha_h,\alpha_g\rho^n)\cap\mathrm{Hom}(\alpha_h,\alpha_g\tilde{\rho}^n)}$ and $v_{k,j}^{g,n}{\in\mathrm{Hom}(\alpha_k\rho,\alpha_g\rho^n)\cap\mathrm{Hom}(\alpha_k\tilde{\rho},\alpha_g\tilde{\rho}^n)}$ are (finitely many) monomials in the Leavitt generators $s,t_l$ and (for each 
fixed pair
$g,n$) the collection $\{u_{h,i}^{g,n}{,u_{h,i}^{g,n\,\prime}},v_{k,j}^{g,n}{,v_{k,j}^{g,n\,\prime}}\}$ together satisfy the Leavitt--Cuntz relations {$u_{h,i}^{g,n\,\prime}u_{k,j}^{g,n}=\delta_{i,j}\delta_{h,k}$ etc. More precisely, $\{u_{h,\star}^{g,n+1}\}=\{v^{g,n}_{h,i}\}_i$ and $\{v^{g,n+1}_{k,\star}\}=\{u^{g,n}_{k,i}\}_i
\cup\{v^{g,n}_{h,j}t_{k-h}\}_{j,h}$.}
%An easy induction shows that both dim$\,\mathrm{Hom}(\alpha_h,\alpha_g\rho^n)$ and  dim$\,\mathrm{Hom}(\alpha_h\rho,\alpha_g\rho^n)$ (which equal respectively  the numbers of $u_{h,*}^{g,n}$ and $v_{h,*}^{g,n}$)
%depend only on $\delta_{g,h}$ and $n$. This implies   dim$\,\mathrm{Hom}(\alpha_h,\alpha_g\rho^n)=\mathrm{dim\,Hom}(\alpha_{-h},\alpha_{-g}\rho^n)$ and  dim$\,\mathrm{Hom}(\alpha_h\rho,\alpha_g\rho^n)=\mathrm{dim\,Hom}(\alpha_{-h}\rho,\alpha_{-g}\rho^n)$, a fact we will need shortly.

Certainly Hom$(\alpha_g\rho^n,\alpha_{g'}\rho^{n'})$ contains all ${u^{g',n'}_{h,i'}u^{g,n\,\prime}_{h,i}}$
and {$v^{g',n'}_{k,j'}v^{g,n\,\prime}_{k,j}$. In fact,  we will show now using  simplicity (Proposition 1) that together  they span
that Hom-space. To see this,  choose any $x\in \mathrm{Hom}(\alpha_g\rho^n,\alpha_{g'}\rho^{n'})$. Then $u^{g,n\,\prime}_{h,i}x
u^{g',n'}_{h',i'}\in \mathrm{Hom}(\alpha_h,\alpha_{h'})=\bbC\delta_{h,h'}$; when $h=h'$, call this number
$q_{h;i,i'}$. Likewise, $v^{g,n\,\prime}_{h,i}x
v^{g',n'}_{h',i'}\in \mathrm{Hom}(\alpha_h\rho,\alpha_{h'}\rho)=\bbC\delta_{h,h'}$; when $h=h'$, call this number
$r_{h;i,i'}$. Moreover, $u^{g,n\,\prime}_{h,i}x
v^{g',n'}_{h',i'}=v^{g,n\,\prime}_{h,i}x
u^{g',n'}_{h',i'}=0$ since $\mathrm{Hom}(\alpha_h,\alpha_{h'}\rho)=\mathrm{Hom}(\alpha_{h}\rho,\alpha_{h'})=0$.
Thus $x=(\sum_{h,i}u^{g',n'}_{h,i}u^{g,n\,\prime}_{h,i}+v^{g',n'}_{h,i}v^{g,n\,\prime}_{h,i})x(\sum_{h',i'}u^{g',n'}_{h',i'}u^{g,n\,\prime}_{h',i'}+v^{g',n'}_{h',i'}v^{g,n\,\prime}_{h',i'})= \sum_{h,i,i'}(q_{h;i',i}u^{g',n'}_{h,i}u^{g,n\,\prime}_{h,i}+r_{h;i',i}v^{g',n'}_{h,i}v^{g,n\,\prime}_{h,i})$. Thus $\mathrm{Hom}(\alpha_g\rho^n,\alpha_{g'}\rho^{n'})=\mathrm{span}_{h,i,i'}\{u^{g',n'}_{h,i'}u^{g,n\,\prime}_{h,i},v^{g',n'}_{h,i'}v^{g,n\,\prime}_{h,i}\}$. The identical argument shows 
$\mathrm{Hom}(\alpha_g\tilde{\rho}^n,\alpha_{g'}\tilde{\rho}^{n'})$ is also spanned by the same elements,
and so those Hom-spaces are identical (and finite-dimensional).} Thus Lemma 2 applies. 

{Recall $\overline{\cC(\cE)}{}^{ds}$, the idempotent completion of $\cC(\cE)$ extended  by direct sums. 
Note that  all $u^{g,n}_{h,i}u^{g,n\,\prime}_{h,i},v^{g,n}_{k,j}v^{g,n\,\prime}_{k,j}$ are idempotents in $\mathrm{End}(\alpha_g\rho^n)$, thanks to the Leavitt--Cuntz relations. Enumerate these $p_1,\ldots,p_N$. Then $p_i
\mathrm{End}(\alpha_g\rho^n)p_j=\delta_{ij}\bbC p_i$ and $\sum_ip_i=\mathrm{id}$, using the
above spanning set (in fact basis) for End$(\alpha_g\rho^n)$, so the $p_i$ form a
complete set of minimal idempotents in  $\mathrm{End}(\alpha_g\rho^n)$. All $(p_i,\alpha_g\rho^n)$ 
 are objects in $\overline{\cC(\cE)}{}^{ds}$. Since
End$(p_i,\alpha_g\rho^n):=p_i\mathrm{End}(\alpha_g\rho^n)p_i=\bbC p_i$ is 1-dimensional, the
$(p_i,\alpha_g\rho^n)$ are simple  in $\overline{\cC(\cE)}{}^{ds}$. These $(p_i,\alpha_g\rho^n)$ (as $i,g,n$ vary)
exhaust  all  simple objects in $\overline{\cC(\cE)}{}^{ds}$, as any other idempotent in $\mathrm{End}(\alpha_g\rho^n)$ is a disjoint sum of the $p_i$. Moreover,} 
 $(u^{g,n}_{h,i}u^{g,n\,\prime}_{h,i},\alpha_g\rho^n)$ and $(1,\alpha_h)$ are 
isomorphic, with isomorphism $u^{g,n}_{h,i}$
{and} inverse $u^{g,n\,\prime}_{h,i}$, since $u^{g,n}_{h,i}u^{g,n\,\prime}_{h,i}$ is the identity in End($(u^{g,n}_{h,i}{u^{g,n\,\prime}_{h,i}},\alpha_g\rho^n)$).
{Likewise,
$(v^{g,n}_{k,j}v^{g,n\,\prime}_{k,j},\alpha_g\rho^n)$ and $(1,\alpha_k\rho)$ are 
isomorphic. }
We thus get a fusion category, because there are only finitely many isomorphism classes
of simple objects, namely the  $[(1,\alpha_g)],[(1,\alpha_g\rho)]$.

{To show that $\overline{\cC(\cE)}{}^{ds}$ is pivotal, note first that $(\alpha_g\rho^n)^{\vee\vee}=
\alpha_g\rho^n$. We want to show also that the double-dual on all intertwiner spaces $\mathrm{Hom}(\alpha_g\rho^n,
\alpha_{g'}\rho^{n'})$ is also the identity map. We must be careful here (and elsewhere) to  
 keep track of the Hom-space we are working in by
writing $(\xi|x|\eta)$ for $x\in \mathrm{Hom}(\xi,\eta)$.  For convenience abbreviate 
$1_\xi=(\xi|1|\xi)$, $s=(\mathrm{id}
|s|\rho^{2})$, $s^{\prime}=(\rho^{2}|s^{\prime}|\mathrm{id})$, $t_g=(\alpha_g\rho|t_g|\rho^{2})$ and $t_g^{\prime}=(\rho^{2}|t^{\prime}_g|
\alpha_g\rho)$. We can compute directly from \eqref{dualmorph} that $(1_\xi)^\vee=1_{\xi^\vee}$, $s^\vee=s'$,
$t_g^\vee=t'_g$, $s^{\prime\vee}=s$, and $t_g^{\prime \vee}=t_g$, and so the double-dual leaves unchanged 
all of these. But the double-dual is a monoidal functor, so it will also leave unchanged the morphisms
 $(\alpha_k\rho^l|s|\alpha_k\rho^{l+2})=s\otimes 1_{\alpha_k\rho^l}$,  $(\alpha_{k+h}\rho^{l+1}|t_{k+2h}|\alpha_{h}\rho^{l+2})=t_{k+2h}\otimes 1_{\alpha_h\rho^l}$, $(\alpha_k\rho^{l+2}|s'|\alpha_k\rho^{l+2})=s'\otimes 1_{\alpha_k\rho^l}$,  and $(\alpha_{h}\rho^{l+2}|t'_{k+2h}|\alpha_{k+h}\rho^{l+1})=t_{k+2h}'\otimes 1_{\alpha_h\rho^l}$.
 By writing $u_{h,i}^{g,n}\in\mathrm{Hom}(\alpha_h,\alpha_g\rho^n)$ and $v_{h,i}^{g,n}\in\mathrm{Hom}(\alpha_h\rho,\alpha_g\rho^n)$ as monomials in $s,t_k$, they can be written as a sequence of compositions of these 
morphisms $(\alpha_k\rho^l|s|\alpha_k\rho^{l+2})$ and $(\alpha_{k+h}\rho^{l+1}|t_{k+2h}|\alpha_{h}\rho^{l+2})$
(this is manifest in the recursions given earlier). Hence the double-dual also leaves unchanged  
 $u_{h,i}^{g,n}$ and $v_{h,i}^{g,n}$. Identical conclusions applies to  $u_{h,i}^{g,n\prime}$ and $v_{h,i}^{g,n\prime}$, and hence to the compositions  $u_{h,i'}^{g',n'} u_{h,i}^{g,n\prime}$ and  $v_{h,i'}^{g',n'}v_{h,i}^{g,n\prime}$.
But those compositions span Hom$(\alpha_g\rho^n,\alpha_{g'}\rho^{n'})$. Thus the double-dual fixes every
morphism $r\in\mathrm{Hom}(\alpha_g\rho^n,\alpha_{g'}\rho^{n'})$. From this we get that the double-dual
functor $X\mapsto X^{\vee\vee}$ is the identity functor on  $\overline{\cC(\cE)}{}^{ds}$, and so  $\overline{\cC(\cE)}{}^{ds}$
is pivotal.
%Recall that the number of $u_{h,*}^{g,n}$ always equals the number of $u_{-h,*}^{-g,n}$, and
%similarly for  $v_{h,*}^{g,n}$ and $v_{-h,*}^{-g,n}$. Define  $p_{h,i}^{g,n}=u_{h,i}^{g,n}u_{h,i}^{g,n\,\prime}$ and $q_{h,i}^{g,n}=v_{h,i}^{g,n}v_{h,i}^{g,n\,\prime}$. Then the choice
%\begin{equation}
%(p_{h,i}^{g,n},\alpha_g\rho^n)^\vee= \left\{\begin{matrix}(p_{h,i}^{g,n},\alpha_g\rho^n)&n\ \mathrm{odd}\\
%(p_{-h,i}^{-g,n},\alpha_{-g}\rho^n)&n\ \mathrm{even}\end{matrix}\right.\ , \ \ 
%(q_{h,i}^{g,n},\alpha_g\rho^n)^\vee= \left\{\begin{matrix}(q_{h,i}^{g,n},\alpha_g\rho^n)&n\ \mathrm{odd}\\
%(q_{-h,i}^{-g,n},\alpha_{-g}\rho^n)&n\ \mathrm{even}\end{matrix}\right.\,,\end{equation}
%where the (co-)evaluation...
%With this choice, $(p_{h,i}^{g,n},\alpha_g\rho^n)^{\vee\vee}=(p_{h,i}^{g,n},\alpha_g\rho^n)$ and $(q_{h,i}^{g,n},\alpha_g\rho^n)^{\vee\vee}=(q_{h,i}^{g,n},\alpha_g\rho^n)$ and our category is pivotal.
The dimension calculation is now trivial: $e_{\alpha_g}b_{\alpha_g^\vee}=1$ and $e_{(\alpha_g\rho^{n})^\vee}b_{\alpha_g\rho^{n}}=\overline{\omega}^nb^{-2n}=\delta_{\pm}^{n}$, {from which we read off that
$X$ and $X^\vee$ have the same dimension for any simple $X$. This means} 
 that $\cC$ is spherical.}

{Now turn to the proof of part (b) {of Theorem 2}. Suppose there is a tensor category equivalence
between $\cC(G^{(i)};\pm^{(i)},\omega^{(i)},A^{(i)})$. Because $\alpha^{(1)}=(1,\alpha^{(1)})$ is
simple, the equivalence must send
$\alpha_g^{(1)}$ to $(p^{(2)},\alpha_h^{(2)}\rho^{(2)\,m})$ for some {(minimal) idempotent $p^{(2)}$ 
and some $\alpha_h^{(2)}\rho^{(2)\,m}$.} Then $\mathrm{id}^{(1)}=\alpha_g^{(1)\,
\nu}{\mapsto}(x,\alpha^{(2)}_k\rho^{(2)\,m\nu})$  for some {$x\in\cL^{(2)},k\in G^{(2)}$.  
But if $m>0$, this}  can never equal id$^{(2)}=(1,\mathrm{id}^{(2)})$.
Similarly, if $\rho^{(1)}\mapsto (p,\alpha^{(2)}_h\rho^{(2)\,m})$ for some $m>1$, then
no object in $\cC(G^{(1)};\pm^{(1)},\omega^{(1)},A^{(1)})$ can get sent to $\rho^{(2)}=(1,\rho^{(2)})$.}

So {our tensor equivalence}  defines a bijection $\pi:G^{(1)}\rightarrow G^{(2)}$ and an element $r\in G^{(2)}$ by
$\alpha_g^{(1)}\mapsto\alpha_{\pi(g)}^{(2)}$ and $\rho^{(1)}\mapsto \alpha_r^{(2)}\rho^{(2)}$.
Thanks to the fusion rules, $\pi$ must be a group isomorphism, and the tensor equivalence
must send $\alpha^{(1)}_g\rho^{(1)}\mapsto \alpha_{\pi(g)+r}^{(2)}\rho^{(2)}$. 

{Although the tensor equivalence will map Hom-spaces to Hom-spaces,
we don't know \textit{a priori} whether it  lifts to a well-defined algebra homomorphism between
the Leavitt algebras, so as above we will be careful to keep track of the Hom-space we are working in by
using the $(\xi|x|\eta)$ notation.  For convenience abbreviate 
$1_\xi^{(i)}=(\xi|1^{(i)}|\xi)$, $s^{(i)}=(\mathrm{id}
|s^{(i)}|\rho^{(i)\,2})$, {and $s^{\prime(i)}$, %=(\rho^{(i)2}|s^{\prime(i)}|\mathrm{id})$, 
$t^{(i)}_g$, %=(\alpha^{(i)}_g\rho^{(i)}|t^{(i)}_g|\rho^{(i)\,2})$ 
$t_g^{\prime(i)}$ %=(\rho^{(i)\,2}|t^{\prime(i)}_g|\alpha^{(i)}_g\rho^{(i)})$. 
similarly.}  Note if the tensor equivalence sends object $\xi$ to object $\xi'$, then it must take the identity $1^{(1)}_\xi$ in  End$(\xi)$  to the identity $1^{(2)}_\xi$ in End$(\xi')$.} 

{By simplicity (Proposition 1), we know}
 $s^{(1)}$ {(which spans Hom$(\mathrm{id}^{(1)},\rho^{(1)\,2})$)} is sent to $\lambda s^{(2)}$ 
  {(which spans Hom$(\mathrm{id}^{(2)},(\alpha^{(2)}_r{\rho^{(2)}})^2)=\mathrm{Hom}(\mathrm{id}^{(2)},\rho^{(2)\,2})$)} and likewise
$t_g^{(1)}\in \mathrm{Hom}(\alpha_g^{(1)}\rho^{(1)},\rho^{(1)\,2})$  to $\mu_g t_{r+\pi g}^{(2)}
\in \mathrm{Hom}(\alpha_{\pi g+r}^{(2)}\rho^{(2)},\rho^{(2)\,2})$ for some non-zero $\lambda,\mu_g\in\bbC$. { Since $1_{\xi}^{(1)}$ is sent to
$1_{\xi'}^{(2)}$, the relations $s'\circ s=1_{\mathrm{id}}$ and $t_g'\circ t_h=1_{\alpha_g\rho}\delta_{g,h}$ give  $s^{\prime\,(1)}
\mapsto \lambda^{-1}s^{\prime\,(2)}$  and 
$t^{\prime\,(1)}_g\mapsto \mu_g^{-1} t^{\prime\,(2)}_{r+\pi{g}}$. From $1_{\alpha_k}\otimes t_h=(\alpha_{k+h}\rho|
\alpha_k(t_h)|\alpha_k\rho^2)$ and $t'_{2k+h}\otimes 1_{\alpha_k}=(\alpha_{k}\rho^2|t'_{2k+h}|\alpha_{k+h}\rho)$ we
obtain 
\begin{equation}\label{tpt12}(t'_{2k+h}\otimes 1_{\alpha_k})\circ(1_{\alpha_k}\otimes t_h)=1_{\alpha_{h+k}\rho}\,;
\end{equation}
hence $1^{(1)}_{\alpha^{(1)}_{h+k}\rho^{(1)}}$ gets sent to both $1^{(2)}_{\alpha^{(2)}_{\pi(h+k)+r}\rho^{(2)}}$ and 
$\mu^{-1}_{2k+h}\mu_h1^{(2)}_{\alpha^{(2)}_{\pi(h+k)+r}\rho^{(2)}}$. Thus $\mu_g=\mu$ is independent of $g$.}
Comparing dimensions {of $\rho^{(1)}$ and $\alpha_r^{(2)}\rho^{(2)}$, we get $\omega^{(1)}b^{(1)\,2}
=\omega^{(2)}b^{(2)\,2}$, i.e.}
we must have $b^{(1)}=\pm b^{(2)}$ {(hence $b^{(1)}=b^{(2)}$ and the signs $\pm^{(1)}$ and
$\pm^{(2)}$ are equal)} {and $\omega^{(1)}=\omega^{(2)}$. The calculation
\begin{equation}\label{ttrs12}t_{0}'\circ (t'_0\otimes \rho)\circ(1_\rho\otimes s)=t_0'\circ(\rho^3|t'_0|\rho^2)\circ (\rho|\rho(s)|\rho^3)=b\,1_\rho\end{equation}
means, computing the image of the tensor equivalence in two ways, $b^{(1)}=\mu^{-2}\lambda b^{(2)}$, which fixes 
the value of $\lambda$. Similarly, the calculation
\begin{equation}\label{ttrtt12}(t_{h+k}'\otimes 1_{\alpha_h}) (t'_h\otimes 1_\rho)(1_\rho\otimes t_0) t_k=(\alpha_h\rho^2|t'_{h+k}|
\alpha_k)(\rho^3|t'_h|\alpha_h\rho^2)(\rho^2|\rho(t_0)|\rho^3)t_k=A_{h,k}1_{\alpha_k\rho}\end{equation}
gives $A_{h,k}=A_{\pi h,\pi k}$.} 

{Note that $s,s',t_g,t'_g$ obey the Leavitt--Cuntz relations \eqref{leavitt}, iff $\pm \mu^2 s,\pm \mu^{-2}s',\mu t_{g+r},\mu t_{g+r}'$ do, for any sign $\pm$, $\mu\in\bbC^\times$ and $r\in G$. These choices leave
unchanged the algebra $\cL_{\nu+1}$ and its endomorphisms $\rho,\alpha_g$. Part (b) follows.}

{Finally, let us turn to part (c) {of Theorem 2}. Suppose $A$ is hermitian.
Define a conjugate-linear map on $\cL$ by ${s}^*=\pm s'$, ${s}^{\prime\,*}=\pm s$, ${t}^*_g=t_g'$
and ${t}^{\prime\,*}_g=t_g$, extended so that $({cxy})^*=\bar{c}{y}^*{x}^*$ for all $c\in\bbC$ and
$x,y\in\cL$, where the sign in these expressions  is as in ${\omega^2}b^{-2}=\delta_\pm$. Then 
$\overline{b}=\pm\omega b$ so} $({\rho(x)})^*
=\rho({x}^*)$ for the $2+2\nu$ generators $x$, and hence that relation holds for all $x\in\cL$.
It is easy to see that this determines a $*$-operation on  $\cC(G;\pm,\omega,A)$, in the sense defined in the 
introduction. If in addition {$\pm=+$},
then this conjugate linear map is the usual $*$-operation on $\cL$, and so taking completions
we get a system of endomorphisms on the Cuntz algebra {which extend to the infinite factor $N$} and thus {we possess} a unitary category.

{Conversely, suppose $\cC(G;\pm,\omega,A)$ possesses a $*$-operation. {Again,} 
we don't know \textit{a priori} whether the $*$-operation (which by definition is defined
only on individual Hom-spaces) lifts to a well-defined $*$-operation on
the Leavitt algebra, so {again write} $(\xi|x|\eta)$ for $x\in \mathrm{Hom}(\xi,\eta)$ {as before}.
 Note that the $*$-operation must take the identity $(\xi|1|\xi)$ in  End$(\xi)$ to
itself.  {From simplicity (Proposition 1),} we may write $t_g^*=\beta_g t'_g$ and $t_g^{\prime\,*}=\beta'_gt_g$ for some non-zero $\beta_g,\beta'_g\in\bbC$. Then taking $*$ of \eqref{tpt12} gives
$\mu_h\mu'_{2k+h}=1$, i.e. that $\mu_g=\mu^{\prime-1}_h=\mu$ is independent of $g,h\in G$.
Now taking $*$ of \eqref{ttrtt12}, we get $\overline{A_{h,k}}=A_{k,h}$, 
and we see that for the category to be hermitian, the matrix $A$ must be hermitian.}

Finally, in a unitary category the categorical dimensions must all be positive.
But $d_{\rho}={\delta_\pm}$, and $\delta_-<0$. {This concludes the proof of Theorem 2.}

\section{The tube algebra and {modular data}}

{\subsection{The tube algebra and its centre}}

{We will now determine the quantum double or centre of our categories $\cC(G;\pm,\omega,A)$ using the tube algebra approach of \cite{iz1}. That approach assumes unitarity, but
\cite{m1} categorises the method, generalising it beyond the context we need, and all of our equations
come from there.}

Let $\Delta=\{\alpha_g,\alpha_h\rho\}_{g,h\in G}$ be a 
finite system of endomorphisms associated to a solution of our equations \eqref{order3sym}-\eqref{quart}. 
{Write $\Sigma \Delta$ for the objects in $\cC(G;\pm,\omega,A)$, and write $[\sigma]$ for
the sector or equivalence class of an object {(where the conjugation now need not be by a unitary)}.}
 The categorical dimension {$d_\sigma=d_{[\sigma]}$} of any object $\sigma\in\Sigma\Delta$ was computed last section.
We found there the dimensions $d_{[\alpha_g]}=1$ and  $d_{[\alpha_g\rho]}=\delta_\pm$ for the simple objects
{(note that $\delta_+>0>\delta_-$, so these dimensions can be negative)}.  The global dimension is then $\lambda_{\pm}=\nu(1+\delta_\pm^2)=
2\nu+\nu^2{\delta_\pm}$, which is strictly  positive as it must be.

The tube algebra 
Tube $\Delta$ is a finite-dimensional algebra over $\bbC$, defined as a vector space by
\begin{equation}\mathrm{Tube}\,\Delta=\oplus_{\xi,\eta,\zeta\in\Delta}\mathrm{Hom}{(\xi\zeta,
\zeta\eta)}\,.\end{equation}
{It will be semisimple even if the fusion category is non-unitary \cite{m1}. As in section 5,} given an element $X$ of Tube$\,\Delta$, it is convenient to write $(\xi\zeta|X|\zeta\eta)$ for the restriction to
Hom$(\xi \zeta,\zeta \eta)$, 
since the same operator may belong to distinct intertwiner spaces. 
For readability we will often write $g$ and $g\rho$ for  $\alpha_g$ and $\alpha_g\rho$, respectively.
In our case the intertwiner spaces are {computed by Lemma 3.}
%\begin{align}\mathrm{Hom}(\alpha_{g+h},\alpha_g\cdot\alpha_h)=&\,
%\mathrm{Hom}({}_{g+h}\rho,g\cdot {}_h\rho)=\mathrm{Hom}({}_{g-h}\rho,{}_g\rho\cdot h)=\bbC 1\,,\\  \mathrm{Hom}(\alpha_{g-h},
%{}_g\rho\cdot{}_h\rho)=&\,\bbC s\,,\ \ \mathrm{Hom}({}_g\rho\cdot{}_h\rho,\alpha_{g-h})=\bbC s'\,,\  \mathrm{Hom}({}_k\rho,{}_g\rho\cdot{}_h\rho)=\bbC t_{k+g-h}\,.\nonumber\end{align} 
Then a basis for  Tube $\Delta$ consists of 
$\cA_{gh}=(g,h|1|h,g)$, $\cB_{gh}=(g,h\rho|1|h\rho, -g)$, $\cC_{gh}=(g\rho,(g-h)/2|1|(g-h)/2,h\rho)$, $\cD_{gkh}=(g,k\rho|t_{2k+g-h}|k\rho,h\rho)$,
$\cE_{gkh}=(g\rho,k\rho|t_{g-h}'|k\rho, h)$, $\cF_{gh}=(g\rho, {(g+h)/2}\rho|ss'|{(g+h)/2}\rho,h\rho)$, and $\cG_{gh}^{kl}=(g\rho,k\rho|t_{l-h+k}t'_{l+g-k}|
k\rho,h\rho)$ (note that the vector space structure of Tube$\,\Delta$ given at the bottom
of p.655 of \cite{iz3} is incomplete). {Thus Tube$\,\Delta$ is $\nu^4+2\nu^3+4\nu^2$-dimensional.}
%\begin{equation}
%\mathrm{Tube}\,\Delta=\bigoplus_{g}\cA_{g,g}\oplus\bigoplus_{g\ne 0}\cA_{g,-g}\oplus\bigoplus_{g,h}
%\cA_{{}_g\rho,h}\oplus\bigoplus_{g,h}\cA_{g,{}_h\rho}\oplus\cB_{\rho,\rho}\,,\end{equation}
%where $\cA_{0,0}=\mathrm{span}$\cA_{g,g}=\bbC B_{g,g}\oplus\mathrm{span}_{k\in G}A_{g,k}$,
%$\cA_{g,h}=\bbC B_{g,h}$ (for $g\ne h$), $\cA_{g,\rho}=\mathrm{span}_{z\in\cF}
%C_{g,z}$, $\cA_{\rho,g}=\mathrm{span}_{z}D_{g,z}$,
%$\cA_{\rho,\rho}=\mathrm{span}_{k}E_k\oplus\mathrm{span}_{k}E'_k\oplus\mathrm{span}_{w,z}E''_{w,z}$. 

The multiplicative structure of Tube $\Delta$ is  given by
\begin{equation}
(\xi\zeta|X|\zeta\eta)(\bar{\xi}\bar{\zeta}|Y|\bar{\zeta}\bar{\eta})=\delta_{\eta,\bar{\xi}}\sum_{\nu\prec\zeta\bar{\zeta}}
(\xi\nu|T(\nu)'\zeta(Y)X\xi(T(\nu))|\nu\bar{\eta})\,,\label{mult}\end{equation}
where we continue to write ${\alpha_g\rho}=g\rho$ and $g$ for $\alpha_g$, and $T(\nu)$ denotes whichever $1,s,t_l$
lies in Hom$(\nu,\zeta\,\bar{\zeta})$. 
In particular, we obtain: $\cA_{gh}\cA_{gl}=\cA_{g,h+l}$, $\cA_{gh}\cB_{gl}=\cB_{gl}\cA_{-g,-h}=\cB_{g,h+l}$, $\cB_{gh}\cB_{-gl}=\cA_{g,h-l}+\delta_{g0}\sum_m\cB_{0m}$,
$\cC_{gh}\cC_{hk}=\cC_{gk}$, %$AC=CA=BC=CB=0$,
$\cA_{gh}\cD_{gkh'}=\cD_{g,h+k,h'}$, 

\noindent $\cB_{g0}\cD_{-gk0}=\sum_lA_{l+k-g,l+k+g}\cD_{g,l,0}$, 
$\cD_{gkh}\cC_{hk'}=\cD_{g,k+(k'-h)/2,k'}$, %$DA=BD=CD=DD=0$, 
$\cE_{gkh}\cA_{hl}=\cE_{g,k-l,h}$, $\cE_{0kh}\cB_{h0}=
\sum_mA_{m+k-h,m+k+h}\cE_{0m,-h}$,  $\cC_{gh}\cE_{hkg'}=\cE_{g,(g-h+2k)/2,g'}$, 

\noindent $\cD_{g0h}\cE_{h0l}=
\delta_{gl}{\bar{\omega}}\cA_{g,0}+\delta_{g,-l}\sum_mA_{m+g+h,2g}\cB_{gm}$,
$$\cE_{0kh}\cD_{hl0}={\bar{\omega}\delta^{-1}_\pm}\delta_{l,k}\cC_{00}+ {\bar{\omega}}\delta_{k,l+h}  \cF_{00}+\sum{}_{g,m}\,A_{m-k+l+h,g-k+l}
A_{m-h+k-l,g+k-l}\cG_{00}^{mg}\,,$$
%$AE=BE=EC=EE=0$,
 $\cC_{gh}\cF_{hk}=\cF_{gk}$, $\cF_{gh}\cC_{hk}=\cF_{gk}$,
$\cD_{k0h}\cF_{h0}={\delta^{-1}_\pm}\cD_{k,-k-h/2,0}$, $\cF_{0h}\cE_{h0l}={\delta^{-1}_\pm}\cE_{0,l-h/2,l}$, 

\noindent$\cF_{0h}\cF_{h0}={\delta^{-3}_\pm}\cC_{00}+\omega\delta^{-2}_\pm\sum_l\cG^{l0}_{00}$, %$AF=FA=BF=FB=FD=EF=0$, 
$\cC_{gh}\cG^{kl}_{hh'}=\cG_{gh'}^{(g-h+2k)/2,(2l+g-h)/2}$,

\noindent $\cG_{hh'}^{kl}\cC_{h'g'}=\cG_{hg'}^{(2k+g'-h')/2,(2l+g'-h')/2}$, 
$\cD_{k0h}\cG_{h0}^{k'l}=\sum_mA_{m+l,m+k+k'}A_{m+k+k',m+l}\cD_{km0}$, 

\noindent $\cG_{0g}^{k'l}\cE_{g0k}=\sum_mA_{m+g-k-k',l-k-k'}A_{l+m,m-k+k'}\cE_{0mk}$, 
$$\cG_{0g}^{kl}\cF_{g0}={\delta^{-2}_\pm}\delta_{2k,g}\cC_{00}+
{\delta^{-1}_\pm}\sum{}_m\,A_{l+m-g/2,2k-g}\cG_{00}^{m,k-g/2}\,,$$
$\cF_{0k}\cG_{k0}^{k'l}={\delta^{-2}_\pm}\delta_{2k',k}\cC_{00}+{\delta^{-1}_\pm}\sum_mA_{m+l-k/2,2k'-k}\cG_{00}^{m,k'-k/2}$, and
\begin{align}\cG_{0h}^{kl}\cG_{h0}^{k'l'}=&\,{\bar{\omega}\delta^{-1}_\pm\delta_{k,k'}A_{l+l'-h,2k-h}}\cC_{00} +\delta_{kl'}{\omega}\delta_{k'l}\cF_{00}\nonumber\\
&+\sum{}_{m,m'}\,A_{m+l'-k,m'+k'-k}A_{m'-k'+h-k,l'-k'-k+l}A_{m+l-k',m'+k-k'}\cG_{00}^{mm'}\,,\nonumber
\end{align} 
%, and $GA=GA=BG=GB=GD=EG=0$.
where we've only written the non-zero products. Note that we have $G$ actions, by multiplying by $\cA$ 
or $\cC$, so for simplicity we restrict to subscripts equal 
to  0 when this $G$-action can yield the other values. 

{Unless $A$ is hermitian, we can't expect Tube$\,\Delta$ to have a natural structure as
a $*$-algebra.} 
%\textcolor{red}{Because Hom$(\tilde{\beta},\tilde{\gamma})=\mathrm{Hom}(\beta,\gamma)$ holds for all
%$\beta,\gamma\in\Sigma \Delta$ (proved in section 5),
%Tube$\,\Delta$ also possesses a $*$-operation (even when the fusion category $\cC(G;\pm,\omega,
%A)$ is not Hermitian). Define first} $R_\zeta\in\mathrm{Hom}(1,{\zeta^\vee}\zeta)$ and $\textcolor{red}{\overline{R}_\zeta\in\mathrm{Hom}(1,\zeta{\zeta^\vee})$ through the equations} $\textcolor{red}{\overline{R}}'_\zeta\zeta(R_\zeta)=d_\zeta^{-1}=R'_\zeta{\zeta^\vee}(\textcolor{red}{\overline{R}}_\zeta)$.
%\textcolor{red}{In particular}, $R_{\alpha_g}=\textcolor{red}{\overline{R}}_{\alpha_g}=1$
%and $R_{\alpha_g\rho}=\textcolor{red}{\overline{R}}_{\alpha_g\rho}=s$.
%The $*$-operation is given by
%\begin{equation}(\xi\zeta|X|\zeta\eta)^*=d_\zeta(\eta{\zeta^\vee}|{\zeta^\vee}(\xi(R^{\vee\,\prime}_\zeta)X')R_\zeta|\textcolor{red}{\zeta^\vee}\xi)\,,\end{equation}
 %so $\cA_{g,h}^*=\cA_{g,-h}$, $\cB^*_{g,h}=
%\cB_{-g,h}$, $\cC_{g,h}^*=\cC_{h,g}$, $\cD^*_{g,h,k}=\cE_{h,k,g}$, $\cF^*_{g,h}=a\cF_{h,g}+a\sum_l
%\cG_{hg}^{(g+h)/2,l}$ and $\cG_{gh}^{kl\,*}=\delta_{g+h,2k}\cF_{h,g}+\sum_{k'}A_{k'+l+g-k,h+g-2k}\cG_{h,g}^{k,k'-g+k}$.

Let {$\sigma\in\Sigma\Delta$}. A half-braiding for $\sigma$ is a choice of
invertible  $\cE_\sigma(\xi)\in\mathrm{Hom}(\sigma\xi,\xi\sigma)$ for each $\xi\in\Delta$, such that for every $\eta,\zeta\in\Delta$ and any $X\in\mathrm{Hom}(\zeta,\xi\eta)$,
\begin{equation} X\cE_\sigma(\zeta)=\xi(\cE_\sigma(\eta))\,\cE_\sigma(\xi)\,\sigma(X)\,.
\label{halfbraiding}\end{equation}
%For our systems this reduces to
%\begin{align}
%&&\cE_\sigma(g+h)=\alpha_g(\cE_\sigma(h))\,\cE_\sigma(g)\,,\label{HBa}\\
%&&\cE_\sigma(\rho)=\alpha_g(\cE_\sigma(\rho))\,\cE_\sigma(g)\,,\label{HBb}\\
%&&U_g\cE_\sigma(\rho)=\rho(\cE_\sigma(g))\,\cE_\sigma(\rho)\sigma(U_g)\,,\label{HBc}\\
%&&S_g\cE_\sigma(g)=\rho(\cE_\sigma(\rho))\,\cE_\sigma(\rho)\,\sigma(S_g)\,,\label{HBd}\\
%&&T_z\cE_\sigma(\rho)=\rho(\cE_\sigma(\rho))\,\cE_\sigma(\rho)\sigma(T_z)\,,\label{HBe}\end{align}
%for all $g,h\in G$, $z\in \cF$. 
{In general, $\sigma$ will be a formal direct sum of $\cL$-endomorphisms, so the values 
$\cE_\sigma(\xi)$ will be matrices with entries in $\cL$. In this case, by $\sigma(X)$ in \eqref{halfbraiding} we mean
the diagonal matrix with entries $\eta(X)$ as $\eta$ runs over all simples in $\sigma$, with multiplicities, and by
$\xi(\cE_\sigma(\zeta))$ we mean to evaluate each entry of the matrix $\cE_\sigma(\zeta)$ by $\xi$.
This equation makes sense as the morphisms ({matrices over} $\cL$) on the left side are intertwiners
for $\sigma\zeta\rightarrow \zeta\sigma\rightarrow\xi\eta\sigma$ while the right side intertwines $\sigma\zeta
\rightarrow \sigma\xi\eta\rightarrow\xi\sigma\eta\rightarrow\xi\eta\sigma$ --- composition is just matrix multiplication. Invertibility of $\cE_\sigma(\xi)$
is equivalent to $\cE_\sigma(\mathrm{id})=1$, the identity matrix.}
There may be more than one half-braiding associated to a given
$\sigma$; in that case we denote them by $\cE_\sigma^j$. 

{The quantum double or centre of the fusion category $\cC=\cC(G;\pm,\omega,A)$ is a strict modular
tensor category (MTC) with objects $(\sigma,\cE_\sigma)$ where $\sigma\in\Sigma\Delta$ and $\cE_\sigma$ is a half-braiding.
The morphisms $x\in\mathrm{Hom}((\sigma,\cE_\sigma),(\tau,\cE_\tau))$  are $x\in\mathrm{Hom}_\cC
(\sigma,\tau)$ satisfying $\sigma(x)\cE_\sigma(\zeta)=\cE_\tau(\zeta)x$ $\forall\zeta\in\Sigma\Delta$; composition is 
as in $\cC$.
The tensor product of objects is given by $(\sigma,\cE_\sigma)\otimes(\tau,\cE_\tau)=(\sigma\tau,\cE_{\sigma\tau})$
where $\cE_{\sigma\tau}(\zeta)=\cE_\sigma(\zeta)\cE_\tau(\zeta)$; the tensor product of morphisms is 
multiplication as in $\cC$. The unit is $(\mathrm{id},1)$. The braiding is $c_{(\sigma,\cE_\sigma),(\tau,\cE_\tau)}
=\cE_\sigma(\tau)$. Duals are $(\sigma,\cE_\sigma)^\vee=(\sigma^\vee,\cE_{\sigma^\vee})$ where
$\cE_{\sigma^\vee}(\zeta)= e_{\sigma^\vee}\otimes\mathrm{id}_{\zeta\sigma^\vee}(\sigma^\vee(\cE_{\sigma}(\zeta)^{-1})\sigma^\vee\zeta(b_\sigma))$; (co-)evaluation is as in $\cC$.
If $\cC$ is hermitian (resp. unitary), one should require the $\cE_\sigma(\xi)$ to be unitary and not merely
invertible, in which case the resulting category will be  a hermitian (resp. unitary) MTC.  See  \cite{m1} for details.}
 
 {Tube$\,\Delta$, being a finite-dimensional semisimple algebra over $\bbC$, decomposes 
 into a direct sum $\oplus_iM_{k_i\times k_i}(\bbC)$ of matrix algebras. The (indecomposable) half-braidings $\cE_\sigma^j$ make this explicit.
 Decompose the sector $[\sigma]$ into a sum $\sum_{i=1}^{k'}[g_i]+\sum_{i=1}^{k''}[h_i\rho]$
of simples, repetitions allowed. In $\cC$, $\sigma$ is the formal direct sum 
$$\sigma=((1,g_1),\ldots,(1,g_{k'}),(1,h_1\rho),
\ldots,(1,h_{k''}\rho))\,,$$
where the 1's denote the identity idempotent (and will be dropped for readability).
Let $k=k'+k''$. Then by \eqref{hom}, for each simple $\xi$ $\cE_\sigma^j(\xi)$ will be a $k\times k$ matrix
with entries $\cE_\sigma^j(\xi)_{\eta,\bar{\eta}}\in\mathrm{Hom}(\eta\xi,\xi\bar{\eta})\subset \cL$,
as $\eta,\bar{\eta}$ run over all simples $\{g_i,h_i\rho\}$ in $\sigma$, repetitions included. The resulting
$k\times k$ matrix algebra $\{\cE_\sigma^j(\xi)\}$ (with entries contained in $\cL$) is isomorphic as a $\bbC$-algebra to an irreducible summand of Tube$\,\Delta$,
and all irreducible summands are of that form.}
  
{We will determine the possible half-braidings $\cE_\sigma^j$, by determining 
 the \textit{matrix units} in Tube $\Delta$ of the corresponding simple summand $M_{k\times k}$.}
Matrix units $e_{i,j}$ of   $M_{k\times k}$
are a basis satisfying $e_{i,j}e_{m,l}=\delta_{j,m}e_{i,l}$. 
 {The relation between the matrix units and the corresponding half-braidings is \cite{iz1}}
\begin{equation}%\label{matrent}
%&&\cE_\sigma^j(\xi)_{{(\eta,i),(\bar{\eta},j)}}=\xi(W_\sigma({\bar{\eta},j})')\,\cE_\sigma^j(\xi)\,W_\sigma({\eta,i})\,,\\ &&
e(\sigma^j)_{{\eta,\bar{\eta}}}=\frac{d_\sigma}{\lambda_{\pm}
\sqrt{d_\eta\,d_{\bar{\eta}}}}\sum{}_\xi\, d_\xi\,(\eta\xi|\cE_\sigma^j(\xi)_{{\eta,\bar{\eta}}}|
\xi\bar{\eta})\,,\end{equation}
{where the sum is over $\xi\in\Delta$, and again $\eta,\bar{\eta}$ run through the simples $\{g_i,h_i\rho\}$ in $\sigma$.
The corresponding central projection (the unit of that simple summand) is then $z(\sigma^j)=\sum_{\eta}e(\sigma^j)_{\eta,
\eta}$.
Our primary interest {this section is in} determining the modular data of the double, and for this
purpose the diagonal matrix units are all that we need.}

{As a $\bbC$-algebra, Tube$\,\Delta$ decomposes as a direct sum 
\begin{equation}\mathrm{Tube}\,\Delta\cong M_{1\times 1}\oplus M_{\nu+1\times\nu+1}\oplus
\frac{\nu-1}{2}M_{\nu+2\times \nu+2}\oplus\frac{\nu^2-\nu}{2}M_{\nu+2\times\nu+2}\oplus\frac{\nu^2+3}{2}M_{\nu\times
\nu}\label{Tubedecomp}\end{equation}
corresponding to half-braidings $\cE_{\mathrm{[id]}}$, $\cE_{[\mathrm{id}]+\Sigma_g[g\rho]}$,
$\cE_{2[\mathrm{id}]+\Sigma_g[g\rho]}^{\psi}$, $\cE_{[h]+[-h]+\Sigma_g[g\rho]}^{\phi}$, and $\cE_{\Sigma_g[g\rho]}^{l}$ respectively, where $\psi,\phi\in\widehat{G}$ but $\psi$ is non-trivial and $\psi,\bar{\psi}$ give the
same half-braiding, $h\in G$ but $h$ is non-trivial and $\pm h$ give same half-braiding, and $1\le l\le \frac{\nu^2+3}{2}$ is some parameter to be interpreted later.}

{The proof of \eqref{Tubedecomp}  for $\cC(G;\pm,\omega,A)$, in particular the determination of the associated matrix units,
 follows the analysis in section 8 of \cite{iz3}, which does this for the (unitary) Q-systems. The main differences
 are the presence of $\omega$ and the absence of a $*$-structure. The central projection of the unique 
 half-braiding of $\sigma=\mathrm{id}$ is again given by
 $z^{i}=\lambda_{\pm}^{-1}\sum\cA_{0g}+\delta\lambda_{\pm}^{-1}\sum\cB_{0g}$, and so
 \begin{equation}\cE^i(g)_{0,0}=\cE^i(g\rho)_{0,0}=1\,.\label{cEi}\end{equation}
 The matrix units corresponding
 to the second summand of \eqref{Tubedecomp} are
\begin{align}e^{ii}_{0,0}&= \frac{\delta_\pm}{\lambda_{\pm}}(\delta_\pm\sum{}_g\,\cA_{0,g}-\sum{}_g\,\cB_{0,g})\,,\nonumber\\
e_{0,g\rho}^{ii}&=\frac{\overline{\omega}\delta_\pm}{\nu\sqrt{\nu\delta_\pm+2}}\sum{}_k\,\cD_{0kg}\,,\ \ 
e_{g\rho,0}^{ii}=\frac{\overline{\omega}\delta_\pm}{\nu\sqrt{\nu\delta_\pm+2}}\sum{}_k\,\cE_{gk0}\,,\ \\
e_{g\rho,h\rho}^{ii}&=\frac{\delta_\pm}{\lambda_{\pm}}\left(\cC_{gh}+\delta_\pm\cF_{gh}+\omega\delta_\pm\sum_{k,l,m}A_{k+m,l+m}A_{k-m,l-m}
\cG_{gh}^{k+g/2+h/2,l+g/2+h/2}\right)\,.\nonumber\end{align}
Compare with Proposition 8.2(2) of \cite{iz3}. This corresponds to
\begin{align}\cE^{ii}(h)_{0,0}&=1\,,\ \cE^{ii}(h)_{g\rho,g\rho}=\delta_{h,0}\,,\  \cE^{ii}(h\rho)_{0,0}=-{\delta_{\pm}^{-2}}\,,\\
\cE^{ii}((k+g)\rho)_{g\rho,g\rho}&=\delta_{k,g}ss'+\omega\sum{}_{l,m}\,A_{k+m,l+m}A_{k-m,l-m}t_{g+k+l}t'_{g-k+l}\,.\nonumber\end{align}
The third class of half-braidings is parametrised by pairs $\{\psi,\overline{\psi}\}$ of non-trivial characters $\psi\in\widehat{G}$, and has diagonal matrix units
\begin{align}e^{iii;\psi}_{0,0}&= \nu^{-1}\sum{}_g\,\psi(g)\cA_{0,g}\,,\ e^{iii;\psi}_{0',0'}= \nu^{-1}\sum{}_g\,\overline{\psi(g)}\cA_{0,g}
%\,,\ e^{iii,\psi}_{\infty,\bar{\infty}}= \nu^{-1}\sum{}_g\psi(g)\cB_{0,g}\,,\ 
%e^{iii,\psi}_{0',0'}= \nu^{-1}\sum{}_g\psi(g)\cB_{0,g}
\,,\\
%e_{\infty,g}^{iii}&=\frac{\overline{\omega}\delta_\pm}{\nu\sqrt{\nu\delta_\pm+2}}\sum{}_k\cD_{0kg}\,,\ \ 
%e_{g,\infty}^{iii}=\frac{\overline{\omega}\delta_\pm}{\nu\sqrt{\nu\delta_\pm+2}}\sum{}_k\cE_{gk0}\,,\ \\
e_{g\rho,g\rho}^{iii;\psi}&=\frac{1}{\nu\delta_\pm}\left(\cC_{gg}+\delta_\pm\cF_{gg}+\omega\delta_\pm\sum_{k,l,m}\psi(m)A_{k+m,l+m}A_{k-m,l-m}
\cG_{gg}^{k+g,l+g}\right)\,.\nonumber\end{align}
This corresponds to
\begin{align}\cE^{iii;\psi}(g)_{0,0}&=\psi(g)\,,\ \cE^{iii;\psi}(g)_{0',0'}=\overline{\psi(g)}\,,\ \cE^{iii;\psi}(g)_{k\rho,k\rho}=\delta_{g,0}\,,\nonumber\\ \cE^{iii;\psi}(g\rho)_{0,0}&=\cE^{iii;\psi}(g\rho)_{0',0'}=0\,,\\ 
\cE^{iii;\psi}((g+k)\rho)_{k\rho,k\rho}&=\delta_{g,0}ss'+\omega\sum{}_{l,m}\,\psi(m)A_{k+m,l+m}A_{k-m,l-m}t_{k+l+g}t'_{k+l-g}\,.\nonumber\end{align}    
The fourth class of half-braidings is parametrised by all characters $\psi\in\widehat{G}$ and non-trivial pairs $\pm h\in G$, and has diagonal matrix units
\begin{align}e^{iv;h,\psi}_{h,h}&= \nu^{-1}\sum{}_g\,\psi(g)\cA_{h,g}\,,\ e^{iv;h,\psi}_{-h,-h}= \nu^{-1}\sum{}_g\,\overline{\psi(g)}\cA_{-h,g}\,,\\
e_{g\rho,g\rho}^{iv;h,\psi}&=\frac{1}{\nu\delta_\pm}\left({\cC_{gg}}+\delta_\pm\overline{\psi(h)}\cF_{gg}+\omega\delta_\pm\sum_{k,l,m}\psi(m){A_{k+h+m,l+m}A_{k-h-m,l-m}}
\cG_{gg}^{k+g,l+g}\right)\,.\nonumber\end{align}
Compare with  Proposition 8.2 of \cite{iz3} (where there is a minor typo there for $e^{iv;h,\psi}_{-h,-h}$). This corresponds to
\begin{align}\cE^{iv;h,\psi}(g)_{h,h}&=\psi(g)\,,\ \cE^{iv;h,\psi}(g)_{-h,-h}=\overline{\psi(g)}\,,\ \cE^{iv;h,\psi}(g)_{k\rho,k\rho}=\delta_{g,0}\,,\nonumber\\ \cE^{iv;h,\psi}(g\rho)_{h,h}&=\cE^{iv;h,\psi}(g\rho)_{-h,-h}=0\,,\\ 
\cE^{iv;h,\psi}((g+k)\rho)_{k\rho,k\rho}&=\delta_{g,0}\overline{\psi(h)}ss'+\omega\sum{}_{l,m}\,\psi(m){A_{k+h+m,l+m}A_{k-h-m,l-m}}t_{k+l+g}t'_{k+l-g}\,.\nonumber\end{align}    
The matrix units for the final summand of \eqref{Tubedecomp} are addressed
next subsection.}

{\subsection{The half-braidings for $\sigma=\sum\alpha_h\rho$}}

{Define $n,\mu,m$ by $\nu=2n+1$ and $\mu=\nu^2+4=2m+1$. The analysis of \cite{iz3} is not complete in determining the matrix units for the final summand of \eqref{Tubedecomp},  
 even if one is only interested
in unitary Q-systems as in \cite{iz3}. Because of  this, it is not possible to determine the modular $S,T$ matrices in general.
(To be fair, \cite{iz3} was mainly interested in the solution for $\nu=3$
corresponding to the Haagerup subfactor, and for this solution his equations do uniquely determine the matrix
units.) In this  subsection we supplement the equations given in \cite{iz3}, Lemma 8.3.}

{Generalising Lemma 8.3 of \cite{iz3} to our context, we learn that the matrix units corresponding to the final summand,
i.e. to the half-braidings with $[\sigma]={\sum_h[\alpha_h\rho]}$, are of the form
\begin{equation}e^{v;j}_{g\rho,h\rho}=\frac{\nu}{\lambda_{\pm}}\left(\cC_{gh}+\overline{w_j}\delta_\pm\cF_{gh}+\delta_\pm\sum{}_{k,l}\,C^j_{k,l}
\cG_{gh}^{k+g/2+h/2,l+g/2+h/2}\right)\,,\label{ejrr}\end{equation}
for $1\le j\le (\nu^3+3)/2$, where the $m(\nu^2+1)$ variables $w_j,C^j_{k,l}\in\bbC$ satisfy
 the $m(\nu^2+1)$ equations:
 \begin{equation}
 \sum{}_{g}\,C_{0,g}^j=w_j-\overline{w_j}\,\delta^{-1}_\pm\,,\
 w_j\,C^j_{g,h}-\sum{}_{k}\, A_{g+k,2h}C^j_{h,k}= \delta_{h,0}\,{\omega} \overline{w_j}\,\delta^{-1}_\pm\,, \label{8.10}\end{equation}
for all $g,h,k\in G$. This half-braiding corresponds to
\begin{equation}
\cE^{v;j}(k)_{g\rho,g\rho}=\delta_{k,0}\ \ \mathrm{and}\ \ 
\cE^{v;j}((k+g)\rho)_{g\rho,g\rho}=\delta_{k,0}\overline{w_j}ss'+\sum{}_l\,C_{k,l}^jt_{l+k+g}t'_{l-k+g}\,,
\end{equation}
for all $g,k\in G$.}

{The $w_j$ {are the corresponding diagonal entries} of the modular matrix $T$, and so must be roots of unity.
 Some  solutions to \eqref{8.10}, occurring for $w_j$ of small order,
  are redundant (i.e. correspond to the previous summands of \eqref{Tubedecomp}) and should be 
dropped. %In particular, 
%precisely $\sigma_1(\nu)-1$ (i.e. the sum of divisors $d>1$ of $\nu$) 
%solutions to \eqref{8.10} with $w_l= 1$ are redundant, and 
% for any other root of unity $w_l$ of order $d$ dividing $\nu$, 
%the number of redundant solutions is $\sum d'\phi(\nu/d')/2$ where $\phi$ is
%Euler's totient and the sum is over all divisors $d'<d$ of $d$. 
}

{Note that, when there is more than one half-braiding $\cE^j_{\Sigma\,g\rho}$ with the same value of
$w_j$, say $w_j=w_{j'}$, then there will be infinitely many different solutions to \eqref{8.10}
with $w_j$, namely $tC^j_{g,h}+(1-t)C^{j'}_{g,h}$ for any $t\in \bbC$. This is because \eqref{8.10}
are linear, for fixed $w_j$. Such $w_j$ can indeed occur --- in \cite{EG2}, 6 fusion categories (in fact Q-systems)  $\cC(\bbZ_\nu;+,{1},A)$
were found with $w_j$ of higher multiplicities. Those examples correspond to $\nu=9,11,19$; for reasons
explained in \cite{EG2}, we expect
there to be higher multiplicities, and hence ambiguities, whenever $\mu$ is composite. Whenever we cannot
uniquely determine the $C^j_{g,h}$, we cannot uniquely determine e.g. the modular $S,T$ matrices.}

{The situation will only get worse as we generalise the context beyond Q-systems to not-necessarily-unitary fusion categories. For this reason, we need to supplement Izumi's \eqref{8.10} with non-linear constraints.
{This is done next.}}

\medskip\noindent\textbf{Proposition 2.} {Let $\cC(G;\pm,\omega,A)$ be any category in Theorem 2. 
Then \eqref{cubic} holds. Moreover, in addition to \eqref{8.10}, $w_j,C^j_{g,h}$ must satisfy:}
\begin{align}
{\omega}w_jC^j_{p,s}C^j_{h,r}\delta_\pm=&\,\delta_{s,h}\delta_{r,p}+\overline{w_j}A_{p+h,2s}\delta_{r,s}
\\&+\delta_\pm\sum{}_{k,l}\,C^j_{k,l}A_{h+l-s,r+k-s}A_{r-k-s,l-k-s+p}{A}_{h+p-k,r+s-k}\,,\label{nonlin1}\\
\frac{\lambda_{\pm}}{\nu}\delta_{j,j'}=&\,1+\overline{w}_{j'}w_j+\delta_\pm\omega w_j
\sum{}_{t,q}\,C^{j'}_{t,q}C^j_{q,t}\,,\label{nonlin2}\\
0=&\,1+\overline{\psi(g)}w_j+\delta_\pm\omega w_j
\sum_{t,q,m}\,\psi(m)\,C^j_{q,t}A_{t+m+g,q+m}A_{t-m-g,q-m}\,,\label{nonlin3}
\end{align}
{for all $\psi\in\hat{G}$, $g\in G$. Conversely, these equations uniquely determine 
$C^j_{g,h}$ and $w_j$.} 

\medskip\noindent\textit{Proof.} {Consider the subalgebra $\mathcal{A}_\rho=\mathrm{span}\{\cC_{00},
\cF_{00},\cG_{00}^{gh}\}$ of Tube$\,\Delta$.  From the products calculated in the previous subsection,
we find that $\cA_\rho$ is commutative with unit $\cC_{00}$.
Now, the diagonal matrix unit $e^j:=e^{v;j}_{\rho,\rho}$ in \eqref{ejrr} is a minimal projection in $\cA_\rho$, and hence 
for any $P\in\cA_\rho$,  $e^jPe^j$ (which equals $Pe^j$ by commutativity) must be a scalar multiple of $e^j$.
Write $e^j\cG_{00}^{gh}=x_{g,h}e^j$ for scalars $x_{g,h}\in\bbC$. We compute 
$$\frac{\lambda_{\pm}}{\nu}x_{s,p}e^j=\cG_{00}^{sp}e^j=\cG_{00}^{sp}+\overline{w_j}\left(\delta_{s,0}\frac{\cC_{00}}{\delta_\pm}
+\sum{}_{k}\,A_{p+k,2s}\cG_{00}^{ks}\right)+\delta_\pm\sum{}_{k,l}\,C_{k,l}^j$$
$$\times\left({\bar{\omega}}\delta_{s,k}A_{p+l,2s}\frac{\cC_{00}}{\delta_\pm}+\omega\delta_{s,l}\delta_{k,p}
\cF_{00}+\sum{}_{h,r}\,A_{h+l-s,r+k-s}A_{r-k-s,l-k-s+p}{A_{h+p-k,r+s-k}}\cG_{00}^{hr}\right)$$
Therefore $x_{s,p}=\omega w_jC^j_{p,s}$ and we recover \eqref{nonlin1}.}

{Similarly, we compute $e^j\cC_{00}=e^j$ and  $e^j\cF_{00}=\delta_\pm^{-1}w_je^j$.
\eqref{nonlin2} and \eqref{nonlin3} now immediately follow 
from  $e^je^{j'}=\delta_{jj'}e^{j'}$ and $e^je^{iii;\psi}_{\rho,\rho}=e^je^{iv;g,\psi}_{\rho,\rho}=0$ respectively.} 
Comparing the $\cG_{00}^{mh}$ coefficients of the associativity of ${\cF_{00}\cF_{00}\cG_{00}^{kl}}$ gives
{\begin{equation}
A_{h+l,2k}A_{m+k,2h}=\frac{\delta_{k,m}\delta_{l,h}}{\delta_\pm}+\sum{}_{r}\,A_{m-k,{h+r-k}}A_{h-r-k,l-r-k}{A}_{m+l-r,h+k-r}\,,\label{nonlin3}\end{equation}
which is equivalent to \eqref{cubic} using \eqref{order3sym}}. 

{Conversely,  the matrix units $e^j$ are uniquely determined by their orthogonality to
those for the other half-braidings, as well as the relations $e^je^{j'}=\delta_{jj'}e^j$. These are equivalent to
\eqref{nonlin2} and \eqref{nonlin3}, once we know  $e^j\cC_{00}=e^j$,  $e^j\cF_{00}=\delta_\pm^{-1}w_je^j$, $e^j\cG_{00}^{gh}=\omega w_jC^j_{h,g}e^j$. These latter equations follow from \eqref{nonlin1} and \eqref{8.10}.}
\textit{QED}\medskip

Curiously, the right-side of \eqref{nonlin2} isn't manifestly symmetric in $j\leftrightarrow j'$, even though the
left-side is. {We know we have a complete list of identities satisfied by $A$, $\omega$ and
$\delta_\pm$, so 
\eqref{nonlin3} (equivalently \eqref{cubic}) is redundant, but it doesn't seem to be trivially redundant. Conversely, we expect that it, in conjunction with  \eqref{order3sym}-\eqref{quad1}, implies the more complicated \eqref{quart}, and
so can replace it in Theorems  1 and 2, but we haven't verified this yet.}

{\subsection{Modular data for the double of $\cC(G;\pm,\omega,A)$}}

\noindent\textbf{Definition 1.} \textit{Modular data} {consists of} a pair $S,T$ of unitary matrices satisfying:

\smallskip\noindent{(i)} $S$ is symmetric (i.e. $S^t=S$) and $T$ is diagonal and of finite order ($T^N=I$);

\smallskip\noindent{(ii)} $S^2$ is a permutation matrix of  order $\le 2$, and $(ST)^3=S^2$;

\smallskip\noindent{(iii)} $S_{1i}\in \bbR\setminus\{0\}$ and some index {$1'$} has $S_{1'i}>0$, for all $i$;

\smallskip\noindent{(iv)} for each $i,j,k$, the numbers $N_{ij}^k$ defined by
\begin{equation}N_{ij}^k:=\sum_{l}\frac{S_{il}S_{jl}\overline{S_{kl}}}{S_{1l}}\label{verl}\end{equation}

are nonnegative integers.\medskip

Any MTC has modular data. The  {index $i$ parametrises} the simple objects (primaries) $X_i$. 
The 
entries {$T_{i,i}$} of the diagonal matrix $T$  (up to normalisation) are eigenvalues of the twist {$\theta_{X_i}=(\mathrm{tr}_{X_i}\otimes\mathrm{id}_{\mathrm{End}\,{X_i}})(c_{X_i,X_i})\in\bbC\,\mathrm{id}_{X_i}$} while those of the  symmetric matrix $S$ are associated
to the Hopf link: up to normalisation, {$S_{i,j}=\mathrm{tr}_{X_i\otimes X_j}(c_{X_i,X_j}\circ c_{X_j,X_i})$}. 
By Proposition 2.12 of \cite{ENO}, the matrices $S$ and $T$ will be unitary in any MTC, even when the category is not unitary (or even hermitian).
`1' corresponds to  the tensor identity $X_1$
and the permutation $S^2$ sends {$i$ to  $i^\vee$, where $[X_i^\vee]=[X_{i^\vee}]$}. \eqref{verl} is called Verlinde's formula, and
the numbers $N_{ij}^k$ are the structure constants {$[X_i\otimes X_j]=\sum_j N_{ij}^k[X_k]$} of the Grothendieck ring of the MTC.

Ignoring the normalisation, those matrices $S$ and $T$ in a MTC define through  $\left({0\atop 1}{-1\atop 0}\right)\mapsto S$, $\left({1\atop 0}{1\atop 1}\right)\mapsto T$
 a projective representation of the modular group SL$_2(\bbZ)=\langle \left({0\atop 1}{-1\atop 0}\right),{\left({1\atop 0}{1\atop 1}\right)}\rangle$, but it is always possible to choose a normalisation
 so that it defines  a linear (i.e. true) representation of SL$_2(\bbZ)$. This choice uniquely determines $S$ up to a sign and then $T$ up to a third root of 1. Property (iv) says the $S$ matrix diagonalises the fusion coefficients
 $N_{ij}^k$, so some column of $S$  (a
 Perron--Frobenius eigenvector) will have constant phase. We require that column (which we call the $1'$th)
 to be strictly positive, as this is necessary  for the existence
of a character vector, as explained in section 7.3. {This will be the case e.g. in a rational CFT.}

From the point of view of modular data, there is little difference between unitary and non-unitary MTCs. In a unitary category, $1'$ must equal $1$.

{Now specialise to the MTCs which are the doubles of the fusion categories $\cC(G;\pm,\omega,A)$
of Theorem 2.}
Write as before $\nu=|G|=2n+1$,  $\mu=\nu^2+4=2m+1$, $\delta_{\pm}=(\nu\pm\sqrt{\mu})/2$ and $\lambda_{\pm}=2\nu+
\nu^2\delta_{\pm}$. The main reason for introducing the tube algebra in section 6.1 is to construct its modular data. The simple objects of the MTC are in one-to-one
correspondence with the simple summands in \eqref{Tubedecomp}, or equivalently with the
irreducible half-braidings.  
 As mentioned 
earlier, in the tube algebra picture, the braidings are given by the half-braidings, and (co-)evaluations
hence traces are as in $\cC$. In particular, we obtain the normalised $S,T$ matrices 
 from the diagonal entries $\cE_\sigma^j(\xi)_{\eta,\eta}$:
 \begin{align}
T_{\sigma^j,\sigma^j}=&\,d_\xi\,\phi_\xi(\cE_\sigma^j(\xi)_{\xi,\xi})\,,\label{Tdef}\\
\overline{S_{ \sigma^i,\bar{\sigma}^{j}}}=&\,\frac{d_\sigma}{\lambda_{\pm}}\sum{}_{\xi}\,d_\xi\,\phi_\xi(
\cE_{\bar{\sigma}}^{j}(\eta)_{\xi,\xi}\,\cE_\sigma^i(\xi)_{\eta,\eta})
\,,\label{Sdef}\end{align}
{for any $j$ in \eqref{Tdef} and \eqref{Sdef}, and any simple $\eta\prec\sigma$ in \eqref{Sdef}.}
 In \eqref{Tdef},
$\xi$ can be any simple in $\sigma$, and in \eqref{Sdef} the sum is over all simple $\xi$  in $\overline{\sigma}$ while $\eta$ is any (fixed) simple in $\sigma$. The
 \textit{standard left inverse} $\phi_\xi$ of the endomorphism $\xi$ is
$\phi_\xi(x)=R_{\xi}'{\xi}^\vee(x)\,R_{\xi}$, where
$R_\zeta\in\mathrm{Hom}(1,{\zeta^\vee}\zeta)$ and {$\overline{R}_\zeta\in\mathrm{Hom}(1,\zeta{\zeta^\vee})$ 
are normalised by ${\overline{R}}'_\zeta\zeta(R_\zeta)=d_\zeta^{-1}=R'_\zeta{\zeta^\vee}({\overline{R}}_\zeta)$. Note that for $x\in\mathrm{End}(\eta\xi)$, $\phi_\xi(x)\in\mathrm{End}(\eta)=\bbC 1$. 
{\eqref{Tdef},\eqref{Sdef} have the desired normalisation built in --- as computed
in section 5.3 of \cite{m1}, the normalisation of $T$ is trivial (i.e. $T_{1,1}=1$)} for the double of any (not necessarily unitary)
fusion category. The derivation of  \eqref{Tdef},\eqref{Sdef} is as in \cite{iz1}, except that the complex conjugate in
\eqref{Sdef} replaces the $*$'s in his Lemma 5.3: his formula assumes $\phi_\xi$ is a $*$-map; equation (5.6)
of \cite{m1} writes this as $S_{ \sigma^{i\,\vee},\bar{\sigma}^{j}}$, which is equivalent to our complex conjugation.}

{In our case, $R_{\alpha_g}=\overline{R}_{\alpha_g}=1$ and
 $R_{\alpha_g\rho}=\overline{R}_{\alpha_g\rho}=s$, so $\phi_{\alpha_g}=\alpha_{-g}$ and $\phi_{\alpha_g\rho}(x)
=s'\alpha_g(\rho(x))s$. We see from \eqref{Tdef},\eqref{Sdef} and the matrix units computed earlier this section
that the modular data is formally identical to that of \cite{iz3} (e.g. $\omega$ doesn't explicitly appear),
except for a trivial dependence on the sign $\pm$. In particular, using \eqref{Tubedecomp},} the primaries fall into four classes:\smallskip

\begin{itemize} 
\item[(i)] two {primaries}, denoted $\mathbf{0}$ and $\mathfrak{b}$;

\item[(ii)]  $n$ primaries, denoted
 $\mathfrak{a}_\psi=\mathfrak{a}_{\overline{\psi}}$ for non-trivial $\psi\in \widehat{G}$;
 
\item[(iii)] $n\nu$ primaries, denoted $\mathfrak{c}_{h,\phi}=\mathfrak{c}_{-h,\phi}$ for $h\in
G,h\ne 0$ and $\phi\in\widehat{G}$;

\item[(iv)] $m$ primaries, denoted $\mathfrak{d}_l$.\end{itemize}

Breaking $S$ and $T$ into 16 blocks, we get
\begin{align}
T=&\,{\rm diag}(1,1;1,\ldots,1;\phi( h);w_1,\ldots,w_m)\,,\cr
S=&\,{1\over \nu} \left(\begin{matrix}B&1_{2\times n}&1_{2\times n\nu}&C\cr 
1_{n\times 2}&2_{n\times n}&D&0_{n\times m}\cr 1_{n\nu\times 2}&D^t&E&
0_{n\nu\times m}\cr C^t&0_{m\times n}&0_{m\times n\nu}&F\end{matrix}\right)\,,\label{genhS2}\end{align}
where  
$k_{a\times b}$ for any $k\in\bbC$ is the $a\times b$ matrix with constant entry $k$, 
$D_{\psi,(h,\phi)}=\psi(h)+\overline{\psi(h)}$, $E_{(h,\phi),(h',\phi')}=\phi'(h)\phi(h')+\overline{\phi'(h)\phi(h')}$, 
{$$B=\frac{1}{2}\left(\begin{matrix}1\mp y&1\pm y\cr 1\pm y&1\mp y\end{matrix}\right)\ {\rm and}\ 
C=\pm y\left(\begin{matrix}1&1&\cdots&1\cr -1&-1&\cdots&-1\end{matrix}\right)$$}
for $y={\nu\over\sqrt{\mu}}$. 
We denote transpose with `$t$'.

{Much more difficult} is to identify the $m\times m$ matrix $F$ and the phases $w_l$. 
Once the solutions $C_{k,l}^j$ and  $w_l$ {to \eqref{8.10},\eqref{nonlin1}-\eqref{nonlin3} have
been found, we conclude}
\begin{equation}
F_{\mathfrak{d}_j,\mathfrak{d}_l}=\frac{\nu}{\lambda_{\pm}}\left(w_jw_l+\delta_\pm
\sum{}_{g,p}\,\overline{C^j_{-g,p}}\,\overline{C^l_{g,p+g}}\right)\,.\end{equation}
{Incidentally,
\eqref{nonlin1} gives an alternate expression for the diagonal entries of $S$:
$$S_{\mathfrak{d}_j,\mathfrak{d}_j}=\frac{1}{\lambda_{\pm}}\left(\omega w_j 
n_3+w_j^2(1-\delta_\pm)+\delta_\pm \omega w_j\sum_{g,h,k,l}\overline{C_{k,l}^j}
\overline{A_{l-p-2g,k-g}}\,\overline{A_{-k-g,l-k-g}}\,\overline{A_{-k,2p+g-k}}\right)$$
where $n_3$ is the number of $g\in G$ with order dividing 3.}

{We have thus identified the $S$ and $T$ matrices for any fusion category $\cC(G;\pm,\omega,A)$,
although the numbers $w_j$ and the submatrix $F$ seem at this point completely opaque.
However, in the following section we list all known fusion categories (unitary or otherwise) of type $\cC(G;\pm,
\omega,A)$, 
and identify their modular data. We will find that the mysterious matrix $F$ and phases $w_l$  always seem to take a  remarkably simple form. For this reason we conjecture:}

\medskip\noindent{\textbf{Conjecture 1.} Choose any finite abelian group $G$ of odd order $\nu$, and choose any
fusion category $\cC=\cC(G;\pm,\omega,A)$. Then there is an abelian group $H$ of order $\mu=\nu^2+4$
and a nondegenerate  bilinear form $\beta$ on $H$, which determines {the
submatrix $F$ and the phases $w_j$ for the double of $\cC$ explicitly.
In particular}, the $m=(\nu^2+3)/2$ primaries $\mathfrak{d}_{ l}$ of class (iv) are parametrised by pairs $\pm l$ of elements in $H$, $l\ne 0$, and 
\begin{align}w_l&=\exp[2\pi\i m\, \beta(l,l)]\,,\\ F_{l,l'}&=\mp {\frac{{2}}{\sqrt{\mu}}}\cos(2\pi\beta(l,l'))\,.\end{align} }

{By a nondegenerate bilinear form $\beta$ on $H$, we mean
$\beta:H\times H\rightarrow
\bbQ/\bbZ$ obeys  $\beta(g+g',h+h')\equiv
 \beta(g,h)+\beta(g,h')+ \beta(g',h)+\beta(g',h')$ (mod 1) for all 
$g,g',h,h'\in H$, and for any non-zero $g\in H$ there is an $h\in H$ such that
$\beta(g,h)\not\equiv 0$ (mod 1).}

{It is possible that not all $G$ and $H$ arise in Conjecture 1. For example, we know of no fusion categories of type $\cC(\bbZ_3\times\bbZ_3;\pm,\omega,A)$ {(\cite{EG2} showed there are no Q-systems for $\bbZ_3\times\bbZ_3$)}, and we know of no fusion categories $\cC(G;\pm,\omega,A)$ whose corresponding modular data has $H=\bbZ_5\times\bbZ_5\times\bbZ_5$ (it would necessarily
have $G=\bbZ_{11}$). But in both cases, we haven't come close to an exhaustive search.}

{This conjecture fits into the grafting framework of section 3.3 of \cite{EG2}. This modular data can be twisted
by $H^3({G\sdprod\bbZ_2};\bbT)$, as explained in section 3 of \cite{EG2}, and indeed as explained
there in section 3.3, non-unitarity is the natural context for some of these twists. We have nothing more to
add to this discussion. As mentioned earlier, the method of this paper can be generalised to even-order $G$ 
{\cite{EGm,EGn,iz4}}, and a
very small number of solutions are known at present. Although the corresponding elements of $S$ and $T$ 
also appear {to be surprisingly} simple, they do not fit into Conjecture 1, and we are not yet prepared to extend the conjecture to
cover them.}

\section{Explicit solutions}

{\subsection{The fusion category classification for small $G$}}

{This subsection obtains all fusion categories $\cC(G;\pm,\omega,A)$ for $|G|\le 5$.
Recall $\delta_\pm=(\nu\pm\sqrt{\nu^2+4})/2$, where
$|G|=\nu$, and Conjecture 1 from section 6.3.}

\medskip\noindent{\textbf{Theorem 3.} The complete list of fusion categories $\cC(G;\pm,\omega,A)$ appearing Theorem 2 for $G=\bbZ_1,\bbZ_3,\bbZ_5$ are (up to equivalence):}

\begin{itemize}
\item[(i)] \textit{for $G=\bbZ_1$:}  exactly one for either  sign; $A=(-1/\delta_\pm)$; both have $\omega=1$; their 
modular data has $H=\bbZ_5$ and $\beta(k,l)=kl/5$ (for `+'), $\beta(k,l)=2kl/5$ (for `--');

\item[(ii)] {\textit{for $G=\bbZ_3$:}   two inequivalent unitary ones with `+', and two inequivalent hermitian but non-unitary ones with `--'; all four have $\omega=1$ and 
\begin{equation}\label{A3x3}A=\left(\begin{matrix} c&d&e\\ d&e&f\\ e&g&d\end{matrix}\right)\,,\end{equation}
where the parameters for these four solutions are
\begin{align} &+:(c,d,e,f,g)=(c_1,d_1,d_2,f_5,f_5)\,,\nonumber\\ 
&+:(c,d,e,f,g)=(c_2,  d_5,d_5,f_1,f_2)\,,\nonumber\\
&-:(c,d,e,f,g)=(c_3,d_6,d_6,f_3, f_4)\,,\nonumber\\ 
&-:(c,d,e,f,g)=(c_4,d_3,d_4,f_6,f_6)\,,\nonumber\end{align}
for $c_i,d_j,f_k$ explicitly defined below; the 
modular data for all four has $H=\bbZ_{13}$, and $\beta(k,l)=kl/13$ resp.  $\beta(k,l)=2kl/13$ for `+' resp. `--'; } 

\item[(iii)] {\textit{for $G=\bbZ_5$:}  two inequivalent unitary ones with `+', and two inequivalent hermitian but non-unitary ones with `--'; all four have $\omega=1$ and 
\begin{equation}\label{A5x5}A=\left(\begin{matrix}c&d&e&f&g\\ d&g&h&i&h\\ e&j&f&i&i\\ f&k&k&e&h\\ g&j&k&j&d
\end{matrix}\right)\,\end{equation}
where the parameters for these four solutions are
\begin{align}
+: (c,d,e,f,g,h,i,j,k)=&\, (c_2,d_1,d_1,d_1,d_1,h_7,h_{11},h_8,h_{10})\,,\nonumber\\
%(aa1,cc2,dd1,dd1,dd1,dd1,hh8,hh10,hh7,hh10), (aa1,cc2,dd1,dd1,dd1,dd1,hh10,hh7,hh11,hh8), (aa1,cc2,dd1,dd1,dd1,dd1,hh11,hh8,hh10,hh7),  (all 4 equiv to )
+:(c,d,e,f,g,h,i,j,k)=&\,(c_4,d_4,d_3,d_6,d_5,h_4,h_2,h_4,h_2)\,,\nonumber\\ %(aa1,cc4,dd3,dd5,dd4,dd6,hh2,hh4,hh2,hh4), (aa1,cc4,dd5,dd6,dd3,dd4,hh4,hh2,hh4,hh2), (aa1,cc4,dd6,dd4,dd5,dd3,hh2,hh4,hh2,hh4) (all 4 are equiv)
-:(c,d,e,f,g,h,i,j,k)= &\,(
c_1,d_2,d_2,d_2,d_2,h_5,h_{12},h_9,h_6)\,,\nonumber\\ %(and the 3 outer automorphism associates). %(aa2,cc1,dd2,dd2,dd2,dd2,hh9,hh6,hh5,hh12),(aa2,cc1,dd2,dd2,dd2,dd2,hh12,hh9,hh6,hh5), (aa2,cc1,dd2,dd2,dd2,dd2,hh6,hh5,hh12,hh9), (all 4 equiv to galois assoc of Hg_5)
-:(c,d,e,f,g,h,i,j,k)=&\,(c_3,d_7,d_{10},d_9,d_8,h_3,h_1,h_3,h_1)\,,\nonumber\end{align}
for $c_i,d_j,h_k$ explicitly defined below; the 
modular data for all four has $H=\bbZ_{29}$, and $\beta(k,l)=kl/29$ resp.  $\beta(k,l)=2kl/29$ for `+' resp. `--'.} 
\medskip
\end{itemize}

 {The two fusion categories for $\nu=1$ are realised by affine $G_2$ at level 1 (`+'), and Yang--Lee (`--').
The first two fusion categories for $\nu=3$ are realised by an even 
subsystem of the Grossman--Snyder system $H_3$ \cite{GrSn} and an even subsystem
 of the Haagerup subfactor. The other two are their Galois associates.
The first two fusion categories for $\nu=5$ are realised by an
 even subsystem of the  Haagerup--Izumi subfactor for $G=\bbZ_5$ found
in \cite{iz3}, and to one of the even subsystems of the Grossman--Snyder system described in 
section 6.6 of \cite{GrSn}. The other two are their Galois associates.}

Our proof of Theorem 3 uses Gr\"obner basis techniques
as implemented in Maple 17.02. First, we find a basis for the ideal generated by the identities of
Theorem 1. Using it, the eigenvalues are found corresponding
to multiplication by each of the variables in the quotient of the polynomial ring by our ideal.
The eigenvalues are the possible values of the variables. All of these steps are completed in a fraction of a second 
for $\nu=3,5$. We then have to determine {(by trial and error) 
which eigenvalues go together to form solutions}. 

$\nu=1$ was worked out in section 2, so
turn to $G=\bbZ_3$.  Consider first $\omega=1$. The order-3 symmetry \eqref{order3sym}
gives us \eqref{A3x3}. These variables $(c,d,e,f,g)$ satisfy \eqref{lin1}-\eqref{quart}. The Gr\"obner basis method
tells us there are precisely 8 solutions. However by Theorem 2, two solutions
$A^{(1)}$, $A^{(2)}$ {yield equivalent fusion categories} if they can be obtained from each other by the
action of Aut$(\bbZ_3)\cong\{\pm 1\}$, i.e. if $A^{(1)}_{i,j}=A^{(2)}_{-i,-j}$ {for all $i,j\in\bbZ_3$.
In other words,}  the 5-tuples $(c,d,e,f,g)$ and $(c,e,d,g,f)$ are equivalent. Up to this equivalence, we then get 4 
solutions, as given in Theorem 3.
{There, $c_1=(2-\sqrt{13})/3$, $c_2=(7-\sqrt{13})/6$, $c_3=(7+\sqrt{13})/6$, $c_4=(2+\sqrt{13})/3$.
$d_1,\ldots,d_4 \approx -.321,.554,{.717}-.329i,.717+.329i$ respectively are the roots of $9x^4-15x^3+7x^2+x-1$, while
$d_5=(1-\sqrt{13})/6$ and $d_6=(1+\sqrt{13})/6$. Finally,
$f_1,\ldots,f_4\approx  .217+.758i, .217-.758i,-.954,.186$ respectively are the roots of $9x^4+3x^3+x^2+5x-1$, 
and $f_5=(1+\sqrt{13})/6$, $f_6=(1-\sqrt{13})/6$.}

{Now consider $G=\bbZ_3$ with $\omega\ne 1$, a nontrivial third root of 1. Then \eqref{order3sym} gives
$$A=\left(\begin{matrix} 0&d&e\\ \overline{\omega}d&\omega e&0\\ \overline{\omega}e&0&\omega d\end{matrix}\right)\,,$$
where the zeros arise for any entry of $A$ fixed by the order-3 symmetry. The quadratic identities \eqref{quad1} give
e.g. $\overline{\omega}(d^2+e^2)=1-1/\delta_{\pm}$ and $d^2+e^2=1$, which are incompatible. Thus there are no solutions
for $G=\bbZ_3$ with $\omega\ne 1$.}

Now turn to $G=\bbZ_5$,  with $\omega=1$.  \eqref{order3sym} gives \eqref{A5x5}. The Gr\"obner basis
method tells us \eqref{lin1}-\eqref{quart} have exactly 16 solutions (as always, half with `+' and half with `--').
As before, we must identify solutions related by the action of Aut$\,G\cong\bbZ_4$, which sends
$(c,d,e,f,g,h,i,j,k)\mapsto(c,e,g,d,f,i,j,k,h)$. This yields the  4 inequivalent fusion categories given in Theorem 3. Explicitly, $c_1=(13+ \sqrt{29})/10$, $c_2=(13- \sqrt{29})/10$, $c_3=(7+\sqrt{29})/5$, and $c_4=(7-\sqrt{29})/5$. 
Also, $d_1=(3-\sqrt{29})/10$, $d_2=(3+\sqrt{29})/10$, 
$d_3\approx-.537$, $d_4\approx -.426$, $d_5\approx -0.032$, $d_6\approx .480$, $d_7\approx .400-.282i$, $d_8\approx .400+.282i$, $d_9\approx .957-.983i$, $d_{10}\approx .957+.983i$, where the final 8 of these $d_i$ are the  roots of the irreducible polynomial 
$625x^8-1375x^7+1275x^6+245x^5-654x^4+152x^3+75x^2-29x-1$. Finally, $h_1\approx-.675$, 
$h_2\approx .218$, $h_3\approx .437$, $h_4\approx .620$, $h_5\approx-1.270$, $h_6\approx -.095$, 
$h_7\approx 0.084-.536i$, $h_8\approx .084+.536i$, $h_9\approx .106$, $h_{10}\approx .534-.099i$, $h_{11}\approx .534+0.099i$, $h_{12}\approx 1.420$, where $h_1$ to $h_4$ are solutions to the
irreducible polynomial $25x^4-15x^3-9x^2+7x-1$, while $h_5$ to $h_{12}$ are solutions to
the irreducible $625x^8-875x^7-525x^6+1110x^5-789x^4+402x^3-95x^2-3x+1$.
%$c=(13\pm \sqrt{29})/10$ (note that $(\delta-2)/(\delta-1)=(13-\sqrt{29})/10$) or $(7\pm\sqrt{29})/5$ (multiplicity 4 for each of these 4 possibilities).

%$d,e,f,g$: $-1/(\delta-1) and $-1/(\delta'-1)$ (mult 4 each)
%-.537550540161056, -.426673514623090, -0.0320262666601846, .480700879303979, .400665722187905-.282754393434430*I, .400665722187905+.282754393434430*I,  .957108998882271-.983969285031469*I, .957108998882271+.983969285031469*I

%$h,i,j,k$: (25*LL^4-15*LL^3-9*LL^2+7*LL-1)^2(625*LL^8-875*LL^7-525*LL^6+1110*LL^5-789*LL^4+402*LL^3-95*LL^2-3*LL+1)
%-.675701517051091, .218333180143626, .437185036337640, .620183300569825 (mult 2)
%-1.27011310232277, -0.0953351438734511, 0.0849880736684989-.536827364819919*I, 0.0849880736684989+.536827364819919*I, .106609869654718, .534270166688226-0.0998097253393428*I, .534270166688226+0.0998097253393428*I, 1.42032189582805
%
%(aa2,cc3,dd10,dd8,dd7,dd9,hh1,hh3,hh1,hh3), 
%(aa2,cc3,dd8,dd9,dd10,dd7,hh3,hh1,hh3,hh1), (aa2,cc3,dd9,dd7,dd8,dd10,hh1,hh3,hh1,hh3) (all 4 are equiv)

{Finally, turn to $G=\bbZ_5$ and $\omega\ne 1$ a nontrivial third root of 1. Write
$$A=\left(\begin{matrix}0&d&e&f&g\\ \bar{\omega}d&\omega g&h&i&\bar{\omega}h\\ \bar{\omega}e&j&\omega f&\omega i&\bar{\omega}i\\ \bar{\omega}f&k&\bar{\omega}k&\omega e&\omega h\\ \bar{\omega}g&\omega j&\omega k&\bar{\omega}j&\omega d
\end{matrix}\right)\,\,.$$
Using the Gr\"obner basis method, it can be shown that \eqref{lin1} and \eqref{quad1} with $h=0,1$ are inconsistent. This concludes the proof of Theorem 3.}

{As is curious aside,  the linear and quadratic identities \eqref{lin1},\eqref{quad1} suffice to fix $A$ for $G=\bbZ_5$,
but for $G=\bbZ_3$ there are 8 spurious solutions which run afoul of the quartic \eqref{quart} (or cubic \eqref{cubic}) identities.} 

{We know of no examples of fusion categories with $\omega\ne 1$.}

Of course, the set of all fusion categories $\cC(G;\pm,\omega,A)$ for fixed $G$ is closed under Galois actions.
Theorem 3 is disappointing, in that all fusion categories for $G=\bbZ_1,\bbZ_3,\bbZ_5$ are Galois associates
of known unitary fusion categories. But we see no reason at all to expect this to continue for larger $G$, {and expect it is an accident of small $G$.} 

{\subsection{Some Q-systems and their doubles}}

{Q-systems are unitary fusion categories coming
from an even part of a subfactor. After Theorem 1 we explained they correspond here to $\omega=1$, `+', and 
$A$ with specified values for $A_{0,g},A_{g,0},A_{g,g}$. \cite{EG2} found several new Q-systems of  {type $\cC(G;\pm,\omega,A)$}, although was unable to identify the modular data of some of them. In this subsection we {use Proposition 3 to} explain how
they all fit into Conjecture 1.}

{A convenient way to express the matrix $A$ of a Q-system, for $G=\bbZ_\nu$, is in terms of numbers $j_2,j_3,\ldots,j_{n+1}\in\bbR$ (recall $\nu=2n+1$): for $0<g<h<\nu$ we have
$$A_{g,h}=\overline{A_{h,g}}=\frac{\sqrt{\delta}}{\delta-1}\exp[\i(j_h-j_g-j_{h-g})]\,,$$
where $j_1=0$ and $j_{n+1+i}=j_{n+1}+j_n-j_{n-i}$ for $1\le i<n$
(see Lemma 7.3 of \cite{iz3}). The Q-systems found in \cite{EG2} correspond to
\begin{align}
%j_2^{(3)}\approx&\, 1.292076\,;\nonumber\\
%(j_2^{(5)},j_3^{(5)})\approx&\,(0.1846862,1.5984702)\,;\nonumber\\
(j_2^{(7)},j_3^{(7)},j_4^{(7)})\approx &\,(2.471228,0.51685555,0.2137724)\,;\nonumber\\
(j_2^{(9)},\ldots,j_5^{(9)})\approx&\,
(2.396976693,2.079251103,-0.2079168419,-2.508673987)\,;\nonumber\\
(j_2^{(9)\prime},\ldots,j_5^{(9)\prime})\approx&\,
(-2.364737070,1.031057162,1.569692175,0.3383837765)\,;\nonumber\\
j^{(11)}\approx &\,(0.9996507,2.7258434,-0.5714203,-1.7797340,
1.2675985)\,,\nonumber\\
j^{(11)\prime}\approx&\,(-2.6444397,-1.7629598,-2.6444440,
2.7572657,0.1128260)\,;\nonumber\\
j^{(13)}\approx &\,( -3.1050384,0.5993399,-0.111708,
-0.969766,1.336848,1.00483129)\,;\nonumber\\
j^{(15)}\approx&\,(-1.0777623,-.7748018,-2.171863,-1.6068402,-.257508,2.092502,.72289565)\,;\nonumber\\
j^{(17)}\approx&\,(-1.466074,.291489,3.130735,-2.693185,1.398153,-.611938,-1.667078,-1.754821)\,;\nonumber\\
j^{(19)}\approx&\,(-2.677465,1.088972,-.899442,.015448,-1.240928,-.493394,1.839879,\nonumber\\
&\qquad -1.525884,-2.084374)\,;\nonumber\\
j^{(19)\prime}\approx&\,( .896858,-.882585,-2.369855,-1.873294,-1.711620,-.119360,2.972018,\nonumber\\
&\qquad-2.460652,.041334)\,, \nonumber\end{align}
where the superscript $(7)$ etc refers to the value of $\nu$. These approximate values suffice to determine
the exact algebraic values of the $A_{g,h}$, and to verify that these do indeed satisfy all the equations
\eqref{lin1}-\eqref{quart}, using the method described in section 3.5 of  \cite{EG4}. (The solutions for $11\le \nu\le 19$
were announced as conjectural in \cite{EG2}, but using  \cite{EG4} have now been
shown to yield exact solutions.)
This list constitutes the complete classification of Q-systems for $\bbZ_7,\bbZ_9$, up to equivalence. There is no Q-system solution for $G=\bbZ_3\times\bbZ_3$ {(more precisely, any such solution would require
nontrivial 2-cocycle twists of \eqref{alpharho})}.}

\medskip\noindent\textbf{Proposition 3.} {The modular data for the 10 Q-systems listed above,
is given by Conjecture 1 with abelian group $H$ and bilinear form $\beta$ given by:}

$j^{(7)}$: $H=\bbZ_{53}$, $\beta(l,l')=ll'/53$;

$j^{(9)}$: $H=\bbZ_{85}$, $\beta(l,l')=ll'/85$;

$j^{(9)}{}'$: $H=\bbZ_{85}$, $\beta(l,l')=12ll'/85$;

$j^{(11)}$: $H=\bbZ_{125}$, $\beta(l,l')=ll'/125$;

$j^{(11)}{}'$: $H=\bbZ_{25}\times \bbZ_5$, $\beta((l_1,l_2),(l'_1,l'_2))={2l_1l'_1/25+2l_2l'_2/5}$;

$j^{(13)}$: $H=\bbZ_{173}$, $\beta(l,l')=ll'/173$;

$j^{(15)}$: $H=\bbZ_{229}$, $\beta(l,l')=ll'/229$;

$j^{(17)}$: $H=\bbZ_{293}$, $\beta(l,l')=ll'/293$;

$j^{(19)}$: $H=\bbZ_{365}$, $\beta(l,l')={ll'/365}$;

$j^{(19)}{}'$: $H=\bbZ_{365}$, $\beta(l,l')=22ll'/365$.\medskip

Given a nondegenerate bilinear form $\beta$ on some abelian group
of order $\nu^2+4$, let $S^\beta,T^\beta$ denote the modular data described in Conjecture 1.
{Section 4.1 of \cite{EG2} proved this proposition for these Q-systems at $\nu=7,13,15,17$, and conjectured
the correct $H$ and $\beta$ for 5 of the 6 remaining. It
 was unable to determine the modular data for the $2+2+2$ Q-systems at $G=\bbZ_9,\bbZ_{11},\bbZ_{19}$, because of the ambiguity described in section 6.2 above. It had no guess for the modular data for $j^{(11)}{}'$
 because it did not think of trying noncyclic $H$.}

{Our proof of Proposition 3 followed very closely what we used in \cite{EG2}, section 4.1. In particular,
a floating point proof is possible and effective, since the integrality of the fusion coefficients $N_{ij}^k$ in
\eqref{verl} serves as error-correction. More precisely,  equation (1.3) of \cite{EG2} shows $S$ in modular
data is uniquely determined from the fusion coefficients, $T$ and the entries $S_{1,i}$. Our strategy here
is to guess at a phase $w_j$ consistent with Conjecture 1, use the linear equations \eqref{8.10} and \eqref{nonlin3} to determine the corresponding $C^j_{g,h}$ up to a small number of parameters (for almost all choices of $w_j$,
this linear system will be inconsistent and we can throw away that choice). For typical examples, the choice
$w_j=e^{2\pi i 182/365}$ for solution $j^{(19)}$ identifies $C^j_{g,h}$ up to 1 parameter, while the choice
$w_j=e^{2\pi i 2/5}$ for solution $j^{(11)\prime}$ needs 4 parameters. Then we chose at random some
nonlinear equations from \eqref{nonlin1} to fix those parameters.}

\subsection{Character vectors}

A natural question is to realise the doubles of these fusion categories by completely rational nets of factors
and/or by rational vertex operator algebras (VOAs). As a first step, one should consider the corresponding character
vectors. This is quite accessible, and {provides considerable information}.

\medskip\noindent\textbf{Definition 2.} Let $\rho$ be a $d$-dimensional representation of SL$_2(\bbZ)$, such that
$T:=\rho{\left({1\atop 0}{1\atop 1}\right)}$ is diagonal. By a \textit{character vector}
$\bbX(\tau)=(\chi_1(\tau),\ldots,\chi_d(\tau))^t$ for $\rho$, we mean a holomorphic function $\bbX$ from the
upper half-plane $\bbH=\{\tau\in\bbC\,|\,\mathrm{Im}\,\tau>0\}$ to $\bbC^d$, which obeys
\begin{equation}\bbX\left(\frac{a\tau+b}{c\tau+d}\right)=\rho\left(\begin{matrix}a&b\\ c&d\end{matrix}\right)\,\bbX(\tau)
\label{translaw}\end{equation}
for all $\tau\in\bbH$ and $\left({a\atop c}{b\atop d}\right)\in\mathrm{SL}_2(\bbZ)$, and for which there exist
exponents $\lambda_k\in\bbR$ and coefficients $\chi_{k;n}\in\bbZ_{\ge 0}$,
such that
\begin{equation}\label{qexp}e^{-2\pi i \lambda_k\tau}\chi_k(\tau)=\sum_{n=0}^\infty \chi_{k;n}q^n\end{equation}
converges absolutely for $|q|<1$, for $k=1,\ldots,d$, where we write $q=e^{2\pi i\tau}$. We also require $\chi_{1;0}=1$.\medskip

Choosing any $\lambda_k'$ so that $T_{kk}=e^{2\pi i\lambda'_k}$, it is clear from holomorphicity and the transformation law \eqref{translaw} at $\left({1\atop 0}{1\atop 1}\right)$,
that $e^{-2\pi i \lambda_k'\tau}\chi_k(\tau)$
is holomorphic in the punctured disc $0<|q|<1$ with an isolated singularity at $q=0$, so \eqref{qexp} should be regarded as a meromorphicity condition
at the so-called cusp $\tau=i\infty$. Any holomorphic $\bbX:\bbH\rightarrow \bbC^d$ obeying \eqref{translaw} and \eqref{qexp}
is called a weakly-holomorphic vector-valued modular {function (vvmf)} for SL$_2(\bbZ)$ with multiplier $\rho$. The
characters of the irreducible modules {$M_j$} for any completely rational conformal net of factors on $S^1$, or 
for any strongly rational VOA,  form a character vector, where $\rho$ is the modular data
coming from the corresponding MTC. The label 1 is the vacuum module $M_1=\cV_1$ (the VOA or net itself), and $T_{11}=e^{-\pi i c/12}$ for a parameter $c$  called the central charge. We can assume without loss of generality that all $\chi_{k;0}\ne 0$, in which case
$h_k=\lambda_k+c/24$ is called the conformal weight of the module $M_j$.
 Because $T$ is only determined by the MTC up to a third root of 1, the  category determines the
central charge only mod 8. For the doubles of fusion categories, as mentioned previously, 
the central charge $c$ is known to be in $8\bbZ$.

{The existence of a character vector is not at all automatic. For one thing, it requires that all $\lambda_j\in\bbQ$,
but that holds in any MTC. Moreover, given any character vector $\bbX(\tau)$, the vector $\mathbf{v}:=\bbX(i)$
exists and is strictly positive (since at $\tau=i$ we have $q=e^{-2\pi}>0$); then \eqref{translaw} says $\mathbf{v}
=S\mathbf{v}$ and hence $S$ must have a strictly positive eigenvector with eigenvalue 1. But we know that
in any modular data, some column (equivalently row, since $S=S^t$) of $S$, namely the common Perron--Frobenius eigenvector of the fusion matrices $N_i=(N_{ij}^k)$, must have constant phase. This is why we 
demanded that a column of $S$ be strictly positive, in section 6.3.}  

When the MTC is unitary, we must have $c\ge 0$ and $h_k>h_1=0$ for $k\ne 1$. {The only unitary VOA or net at
 $c=0$ is the  trivial theory}. In the unitary case, the positive row
of $S$ must be the first (=vacuum) row. When unitarity is dropped,  then
{$h_k\ge h_{1'}$ for all $k$}. The quantity $c_{\mathrm{eff}}=-24h_{1'}$ is called the effective central charge, and must be nonnegative. Again, $c_{\mathrm{eff}}=0$ can only
occur for the trivial VOA and conformal net. To our knowledge, all known examples have $h_j>h_{1'}$ for $j\ne 1'$,
{but this is not yet a theorem.}

The Hauptmodul $j(\tau)=q^{-1}+744+196884q+\cdots$ of SL$_2(\bbZ)$ is a weakly-holomorphic modular function for the trivial
multiplier.   
For any $\rho$ in Definition 2, the space $\cM^!(\rho)$ of weakly-holomorphic vvmfs is
trivially a module over the polynomial ring $\bbC[j(\tau)]$. It turns out that this module is always free of
rank $d$ {(Theorem 3.3(a) of \cite{vvmf})}. Put another way, there is a $d\times d$ matrix 
$$\Xi(\tau)=q^\Lambda\sum_{n=0}\Xi_nq^n\,,$$ with coefficients $\Xi_n\in M_{n\times n}(\bbC)$, with the property that $\bbX(\tau)\in\cM^!(\rho)$ iff there is a vector-valued polynomial
$p(x)\in\bbC^d[x]$ such that $\bbX(\tau)=\Xi(\tau)p(j(\tau))$. So knowing all weakly-holomorphic
vvmfs for $\rho$ is equivalent to knowing $\Xi(\tau)$. We can and will require $\Xi_0=I_{d\times d}$.
The matrix $\Lambda$ will be diagonal, with entries satisfying $T_{kk}=e^{2\pi i\Lambda_{kk}}$.
There is a recursion uniquely determining each $\Xi_n$ from the complex matrices $\Lambda$ and $\Xi_1$
{(equation (36) of \cite{vvmf})}.
In short,  knowing all weakly-holomorphic
vvmfs for $\rho$ is equivalent to knowing the exponents $\Lambda$ and the first
nontrivial coefficient matrix $\Xi_1$.

Once $\Xi(\tau)$ (or equivalently $\Lambda,\Xi_1$) are known, it is then just combinatorics to find all
character vectors for a given effective central charge (since $c_{\mathrm{eff}}$ directly gives bounds for the degrees of all component polynomials $p_k$ in $p(x)\in\bbC^d[x]$).
In \cite{EG2}, this procedure was done for several doubles, including the double of the Haagerup fusion
categories, for central charges $8,16,24$.

To illustrate this for a non-unitary example, in this subsection we give $\Xi(\tau)$ for the non-unitary
cousin of the Haagerup ($G=\bbZ_3$). Its fusion category and MTC is a Galois associate of that
of the Haagerup. By contrast, $\Xi(\tau)$ and hence the corresponding VOA or conformal net, are not at all related in
an obvious way to those of the Haagerup, as we'll see.

The double of either of the unitary fusion categories for $G=\bbZ_3$, at any (effective) central charge $c=c_{\mathrm{eff}}\equiv 8$
(mod 24) (one of the three possibilities), was
found in \cite{EG2} to have $\Lambda$ resp. $\Xi_1$ be
\begin{align}&\,{\rm
diag}({\scriptsize -1/3,-1/3,-1/3,-1/3,-1,-2/3,-34/39,-19/39,-5/39,-37/39,-31/39,-28/39})\,,\nonumber\\
&\,{\scriptsize \left(\begin{array}{cccccccccccc}6& 80& 81& 81& 8748& 1215
& 3549& 273& 13& 5538& 2275& 1378\cr 
80& 6& 81& 81& 8748& 1215& -3549& -273& -13& -5538& -2275& -1378\cr 
81& 81& 167& -81& -8748& -1215& 0& 0& 0& 0& 0& 0\cr 
81& 81& -81& 167& -8748& -1215& 0& 0& 0& 0& 0& 0\cr 
3& 3& -3& -3& -12& 18& 0& 0& 0& 0& 0& 0\cr 
27& 27& -27& -27& 1458& -152& 0& 0& 0& 0& 0& 0\cr
7& -7&0& 0& 0& 0& -88& -14& -1& 50& 63& 64\cr
42& -42&0&0&0& 0&-1484&92&16&2940& -192&-1041\cr
119&-119&0& 0&0&0&-2142&987&11&-24990&-6035&4641\cr
5& -5& 0& 0& 0& 0&17&13&-3&-2& 35&-14\cr
13&-13&0& 0&0&0&174&-1&-5&294&-147& 51\cr
14&-14&0&0& 0& 0&448&-77& 7&-343& 125& -24\end{array}\right)}\end{align}
%7& -7& 0& 0& 0& 0& -1& -14& 50& 64& -88& 63\cr 
%42& -42& 0& 0& 0& 0& -1484& 92& 16& 2940& -192& -1041\cr 
%119& -119& 0& 0& 0& 0& -2142& 987& 11& -24990& -6035& 4641\cr 
%5& -5& 0& 0& 0& 0& 17& 13& -3& -2& 35& -14\cr 
%13& -13& 0& 0& 0& 0& 174& -1& -5& 294& -147& 51\cr 
%14& -14& 0& 0& 0& 0& 448& -77& 7& -343& 125& -24
(we are following the conventions of \cite{vvmf}, which has $\Lambda$ shifted by the identity from the $\Lambda$ 
used in \cite{EG2}). Here, the positive row of $S$ is $1'=1$, the vacuum 0.
At (effective) central charge $c=8$, the polynomial $p(x)$ will be $(\alpha,0,0,0,0,0,0,0,\beta,0,0,0)^t$ for constants
$\alpha,\beta\in\bbC$ (otherwise $\lambda_1$ would not be the unique minimum). But $\alpha=1$, since $\chi_{1;0}=1$.
Thus the only possible character vectors at central charge $c=8$ are
$${\scriptsize\left(\begin{matrix} \chi_0(\tau)\cr \chi_{\mathfrak{b}}(\tau)\cr \chi_{\mathfrak{a}}(\tau)=\chi_{\mathfrak{c}_0}(\tau)\cr \chi_{\mathfrak{c}_1}(\tau)\cr
\chi_{\mathfrak{c}_2}(\tau)\cr \chi_{\mathfrak{d}_1}(\tau)\cr \chi_{\mathfrak{d}_2}(\tau)\cr \chi_{\mathfrak{d}_3}(\tau)\cr \chi_{\mathfrak{d}_4}(\tau)\cr \chi_{\mathfrak{d}_5}(\tau)\cr \chi_{\mathfrak{d}_6}(\tau)
\end{matrix}\right)}=\Xi(\tau){\scriptsize\left(\begin{matrix}1\\ 0\\ 0\\ 0\\ 0\\ 0\\ 0\\ 0\\ \beta\\ 0\\ 0\\ 0\end{matrix}\right)}=
{\scriptsize\left(\begin{matrix}q^{-1/3}\left(1+(6+13\beta)q
+(120+78\beta)q^2+(956+351\beta)q^3+(6010+1235\beta)q^4+\cdots\right)\cr 
q^{2/3}\left((80-13\beta)+(1250-78\beta)q+(10630-351\beta)q^2+ 
(65042-1235\beta)q^3+\cdots\right)\cr 
q^{2/3}\left(81+1377q+11583q^2+71037q^3+\cdots\right)\cr 
 3+243q+2916q^2+21870q^3+\cdots\cr 
q^{1/3}\left(27+594q+5967q^2+39852q^3+\cdots\right)\cr q^{5/39}\left((7-\beta)
+(292-6\beta)q+(3204-43\beta)q^2+(23010-146\beta)q^3+\cdots\right)\cr
q^{20/39}\left((42+16\beta)+(777+121\beta)q+(7147+547\beta)q^2+
(45367+2000\beta)q^3+\cdots\right)\cr 
q^{-7/39}\left(\beta +(11\beta+119)q+(73\beta+1623)q^2+(300\beta+12996)q^3+
(76429+1063\beta)q^4+\cdots\right)\cr q^{2/39}\left((5-3\beta)+(229-50\beta)q
+(2738-252\beta)q^2+(19942-1032\beta)q^3+\cdots\right)\cr
 q^{8/39}\left((13-5\beta)+(347-37\beta)q+(3804-212\beta)q^2+
(26390-794\beta)q^3+\cdots\right)\cr
 q^{11/39}\left((14+7\beta)+(441+61\beta)q+(4445+303\beta)q^2+
(30329+1167\beta)q^3+\cdots\right)
\end{matrix}\right)}$$
The first coefficient of $\chi_{\mathfrak{d}_3}(\tau)$ (i.e. $\chi_{9;0}$) tells us $\beta\in\bbZ_{\ge 0}$,
while the first coefficient of $\chi_{\mathfrak{d}_4}(\tau)$ (i.e. $\chi_{10;0}$) then implies $\beta=0,1$.
Thus there are only two possible character vectors for the Haagerup modular
data at $c=8$,  as stated in \cite{EG2}.

The double of either of the non-unitary fusion categories for $G=\bbZ_3$, at effective central charge $c_{\mathrm{eff}}\equiv 8$ (mod 24) (one of three possibilities), has $\Lambda$  resp. $\Xi_1$ equal to
\begin{align}&\mathrm{diag}({\scriptsize-1/3, -1/3, -1/3, -1/3, -1, -2/3, -16/39, -25/39, -40/39, -22/39, -49/39, -4/39})\,,\nonumber\\
&\,{\scriptsize \left(\begin{array}{cccccccccccc}110& -24& 81& 81& -4374& 1215& -390& -1820& -16770& -910& -53872& 52\cr -24& 110& 81& 81& -4374& 1215& 390& 1820& 16770& 910& 53872& -52 \cr 81& 81& 167& -81& 4374& -1215& 0& 0& 0& 0& 0& 0\cr 81& 81& -81& 167& 4374& -1215& 0& 0& 0& 0& 0& 0\cr -6& -6& 6& 6& -12& -36& 0& 0& 0& 0& 0& 0\cr 27& 27& -27& -27& -729& -152& 0& 0& 0& 0& 0& 0\cr -28& 28& 0& 0& 0& 0& 143& -405& -9580& -518& 3654& -1\cr -1/2& 1/2& 0& 0& 0& 0& -81& -262& 1457& 56& -3832& 26\cr -1/2& 1/2& 0& 0& 0& 0& -7& 7& -12& 6& -7& 1\cr -28& 28& 0& 0& 0& 0& -35& 120& 1820& -314& 7224& -27\cr -1/2& 1/2& 0& 0& 0& 0& 2& 2& 1& -1& 0& -1\cr -57/2& 57/2& 0& 0& 0& 0& 399& 2660& 8436& -854& -204212& 79\end{array}\right)}
\nonumber\end{align}
Here, the positive row of $S$ is $l=2$, the primary $\mathfrak{b}$.
At effective central charge $c_{\mathrm{eff}}=8$ for this $\rho$, the polynomial $p(x)$ will be $(\alpha,\beta,\gamma,\delta,0,0,0,0,0,0,0,\epsilon)^t$ for constants
$\alpha,\beta,\gamma,\delta,\epsilon\in\bbC$ (otherwise $\lambda_2$ would not be the unique minimum). 
Thus the only possible character vectors at effective central charge $c_{\mathrm{eff}}=8$ are
$$%{\scriptsize\left(\begin{matrix} \chi_0(\tau)\cr \chi_{\mathfrak{b}}(\tau)\cr \chi_{\mathfrak{a}}(\tau)=\chi_{\mathfrak{c}_0}(\tau)\cr \chi_{\mathfrak{c}_1}(\tau)\cr
%\chi_{\mathfrak{c}_2}(\tau)\cr \chi_{\mathfrak{d}_1}(\tau)\cr \chi_{\mathfrak{d}_2}(\tau)\cr \chi_{\mathfrak{d}_3}(\tau)\cr \chi_{\mathfrak{d}_4}(\tau)\cr \chi_{\mathfrak{d}_5}(\tau)\cr \chi_{\mathfrak{d}_6}(\tau)
%\end{matrix}\right)}=
{\scriptsize\left(\begin{matrix}q^{-1/3}\left(\alpha+(110\alpha+52\epsilon-24\beta+81\gamma+81\delta)q+(1589\alpha-219\beta+1377\gamma+1377\delta+650\epsilon)q^2+(12721\alpha-1135\beta+11583\gamma+11583\delta+4108\epsilon)q^3+\cdots%(-3988\beta+71037\gamma+71037\delta+75040\alpha+19812\epsilon)q^4
\right)\\ 
q^{-1/3}\left(\beta+(110\beta-24\alpha+81\gamma+81\delta-52\epsilon)q+(1589\beta+1377\gamma+1377\delta-650\epsilon-219\alpha)q^2+(12721\beta+11583\gamma+11583\delta-4108\epsilon-1135\alpha)q^3+\cdots%(75040\beta+71037\gamma+71037\delta-3988\alpha-19812\epsilon)q^4
\right)\\ 
q^{-1/3}\left(\gamma+(167\gamma+81\alpha-81\delta+81\beta)q+(2747\gamma-1377\delta+1377\alpha+1377\beta)q^2+(23169\gamma-11583\delta+11583\alpha+11583\beta)q^3+\cdots%(71037\alpha+142089\gamma+71037\beta-71037\delta)q^4
\right)\\ 
q^{-1/3}\left(\delta+(167\delta+81\alpha+81\beta-81\gamma)q+(2747\delta+1377\alpha+1377\beta-1377\gamma)q^2+(23169\delta+11583\alpha+11583\beta-11583\gamma)q^3+\cdots%(71037\alpha+71037\beta-71037\gamma+142089\delta)q^4
\right)\\ 
-6\alpha-6\beta+6\gamma+6\delta+(-486\alpha-486\beta+486\gamma+486\delta)q+(-5832\alpha-5832\beta+5832\gamma+5832\delta)q^2+\cdots%(-43740\alpha-43740\beta+43740\gamma+43740\delta)q^3
\\ 
q^{1/3}\left(27\alpha+27\beta-27\gamma-27\delta+(594\alpha+594\beta-594\gamma-594\delta)q+(5967\alpha+5967\beta-5967\gamma-5967\delta)q^2+\cdots%(39852\alpha+39852\beta-39852\gamma-39852\delta)q^3
\right)\\ 
q^{23/39}\left(-28\alpha+28\beta-\epsilon+(-1025\alpha/2+1025\beta/2-52\epsilon)q+(-4359\alpha+4359\beta-378\epsilon)q^2+\cdots%(-27016\alpha+27016\beta-2029\epsilon)q^3
\right)\\ 
q^{14/39}\left(-\alpha/2+\beta/2+26\epsilon+(-95\alpha+95\beta+352\epsilon)q+(-1416\alpha+1416\beta+2431\epsilon)q^2+\cdots%(-11220\alpha+11220\beta+12337\epsilon)q^3
\right)\\ 
q^{-1/39}\left(-\alpha/2+\beta/2+\epsilon+(-67\alpha+67\beta+53\epsilon)q+(-932\alpha+932\beta+431\epsilon)q^2+\cdots%(-7440\alpha+7440\beta+2485\epsilon)q^3
\right)\\ 
q^{17/39}\left(-28\alpha+28\beta-27\epsilon+(-512\alpha+512\beta-378\epsilon)q+(-8585\alpha/2+8585\beta/2-2510\epsilon)q^2+\cdots%(-(52359/2)\alpha+(52359/2)\beta-12638\epsilon)q^3
\right)\\ 
q^{-10/39}\left(-\alpha/2+\beta/2-\epsilon+(-67\alpha+67\beta-53\epsilon)q+(-904\alpha+904\beta-457\epsilon)q^2+\cdots%(-6956\alpha+6956\beta-2785\epsilon)q^3
\right)\\ 
q^{-4/39}\left(\epsilon+(79\epsilon-57\alpha/2+57\beta/2)q+(756\epsilon-579\alpha+579\beta)q^2+(4513\epsilon-5196\alpha+5196\beta)q^3+\cdots%(-33040\alpha+33040\beta+20622\epsilon)q^4
\right)
\end{matrix}\right)}$$
We see that $\alpha,\beta,\gamma,\delta,\epsilon\in\bbZ_{\ge 0}$ and $\alpha\equiv\beta$ (mod 2); in fact $\beta>0$
since $c_{\mathrm{eff}}=8$.
Comparing the leading terms of $\chi_5(\tau)$ and $\chi_6(\tau)$, we must have $\alpha+\beta=\gamma+\delta$
and hence also $\gamma\equiv\delta$ (mod 2). This means that the $q,q^2,q^3$ coefficients of $q^{1/3}\chi_1(\tau)$ are all even and thus cannot equal 1. Hence either $c\le -88$ or $\alpha=1$. Assume $\alpha=1$. Then 
all coefficients of e.g. $\chi_5(\tau)$ up to at least $q^4$ vanish. We don't have a proof yet that there is no
character vector with $c_{\mathrm{eff}}=8$ for this $\rho$, but it seems highly likely.

This calculation is meant to give further evidence that, even though these unitary and non-unitary  fusion categories and hence MTCs are related simply by a Galois automorphism, the relation if any between corresponding VOAs or
conformal nets will be far from straightforward.

\newcommand\biba[7]   {\bibitem{#1} {#2:} {\sl #3.} {\rm #4} {\bf #5,}
                    {#6 } {#7}}
                    \newcommand\bibx[4]   {\bibitem{#1} {#2:} {\sl #3} {\rm #4}}

\def\ASENS            {Ann. Sci. \'Ec. Norm. Sup.}
\def\AM   {Acta Math.}
   \def\AnM              {Ann. Math.}
   \def\CMP              {Commun.\ Math.\ Phys.}
   \def\IJM              {Internat.\ J. Math.}
   \def\JAMS             {J. Amer. Math. Soc.}
\def\JFA              {J.\ Funct.\ Anal.}
\def\JMP              {J.\ Math.\ Phys.}
\def\JRA              {J. Reine Angew. Math.}
\def\JSP              {J.\ Stat.\ Physics}
\def\LMP              {Lett.\ Math.\ Phys.}
\def\RMP              {Rev.\ Math.\ Phys.}
\def\RNM              {Res.\ Notes\ Math.}
\def\RIMS             {Publ.\ RIMS.\ Kyoto Univ.}
\def\Inv              {Invent.\ Math.}
\def\npbp             {Nucl.\ Phys.\ {\bf B} (Proc.\ Suppl.)}
\def\nupb             {Nucl.\ Phys.\ {\bf B}}
\def\nup              {Nucl.\ Phys. }
\def\nupp             {Nucl.\ Phys.\ (Proc.\ Suppl.) }
\def\adma             {Adv.\ Math.}
\def\coma             {Con\-temp.\ Math.}
\def\PAMS             {Proc. Amer. Math. Soc.}
\def\PJM              {Pacific J. Math.}
\def\ijmp             {Int.\ J.\ Mod.\ Phys.\ {\bf A}}
\def\jpa              {J.\ Phys.\ {\bf A}}
\def\PLB              {Phys.\ Lett.\ {\bf B}}
\def\RIMS             {Publ.\ RIMS, Kyoto Univ.}
\def\Top               {Topology}
\def\TAMS             {Trans.\ Amer.\ Math.\ Soc.}

\def\Duke              {Duke Math.\ J.}
\def\K                 {K-theory}
\def\JOP               {J.\ Oper.\ Theory}

\vspace{0.2cm}\addtolength{\baselineskip}{-2pt}
\begin{footnotesize}
\noindent{\it Acknowledgement.}

{The authors thank  Cardiff School of Mathematics, University of Alberta and University of Warwick  for generous hospitality while researching this
paper. They also benefitted greatly from Research-in-Pairs held at Oberwolfach and Research in Teams at BIRS. Their research was supported in part by  EPSRC, PIMS
and NSERC. The first named author would like to thank Chris Phillips for information and  references to the $K$-theoretic aspects of Leavitt algebras.}

\end{footnotesize}


\begin{thebibliography}{99}
    \begin{scriptsize}

    \addcontentsline{toc}{section}{References}

    \setlength{\parskip}{-1ex}
    
 \bibx{ABC}{Ara, P.,  Brustenga, M, Corti\~{n}as, G}
{K-theory of Leavitt path algebras.}
M\"{u}nster J. of Math. \textbf{2}  (2009), 5--34.    
    
\biba{ahaag}{Asaeda, M., Haagerup, U.}{Exotic subfactors of finite depth
with Jones indices $(5+\sqrt{13})/2$ and
$(5+\sqrt{17})/2$} \CMP{202} {1--63}{(1999)}.


%\bibx{BG1}
%{Bantay, P., Gannon, T.}
%{Vector-valued modular functions for the modular group and the hypergeometric equation.} {Commun. Number Theory Phys.} \textbf{1}, {637--666} 
%{(2008)}.


%\bibx{BMPS} {Bigelow, S., Morrison, S., Peters, E., Snyder, N.} {Constructing the extended  
%Haagerup planar algebra.} arXiv: math.OA0909.4099.

%\bibx{Bisch} {Bisch, D.} {On the structure of finite depth subfactors.} {In: Algebraic methods in 
%operator theory.}  Birkh\"auser, Boston 1994, pp. 175--194.

%\bibx{BiHa} {Bisch, D., Haagerup, U.} {Composition of subfactors: new examples of infinite depth subfactors.} {Ann. Sci. \'Ecole Norm. Sup.}  {\textbf{29}}, {329--383} (1996). 

%\biba{BE1}
%{B\"ockenhauer, J., Evans, D.E.} {Modular invariants, graphs and
%$\alpha$-induction for nets of subfactors, I} \CMP{197} {361--386}
%{(1998)}.

%\biba{BE2}
%{B\"ockenhauer, J., Evans, D.E.} {Modular invariants, graphs and
%$\alpha$-induction for nets of subfactors, II}\CMP{200}
%{57--103}{(1999)}.

%\biba{BE3}
%{B\"ockenhauer, J., Evans, D.E.} {Modular invariants, graphs and
%$\alpha$-induction for nets of subfactors,
%III}\CMP{205}{183--228}{(1999)}.

\biba{BE4}
{B\"ockenhauer, J., Evans, D.E.} {Modular invariants from
subfactors: Type I coupling matrices and intermediate
subfactors}\CMP{213} {267--289} {(2000)}.

%\bibx{BE5}
%{B\"ockenhauer, J., Evans, D.E.} {\sl Modular invariants from
%subfactors.} {In: {Quantum Symmetries in Theoretical Physics
%and Mathematics} (Bariloche, 2000),  Contemp. Math.
%294, Amer. Math. Soc., Providence 2002, pp.\ 95--131}.

%\bibx{BE6}
%{B\"ockenhauer, J., Evans, D.E.} {Modular invariants and
%subfactors.} In: {Mathematical Physics in Mathematics and Physics
%(Siena, 2000). Fields Inst.\ Commun. 30, Amer.\
%Math.\ Soc., 2001, pp.\ 11--37}.

%\biba{BEK1}
%{B\"ockenhauer, J., Evans, D.E., Kawahigashi, Y.} {On
%$\a$-induction, chiral generators and modular invariants for
%subfactors} \CMP{208} {429--487} {(1999)}.

%\biba{BEK2}
%{B\"ockenhauer, J., Evans, D.E., Kawahigashi, Y.} {Chiral
%structure of modular invariants for subfactors} \CMP{210}
%{733--784} {(2000)}.

%\biba{BEK3}{B\"ockenhauer, J., Evans, D.E., Kawahigashi, Y.}
%{Longo-Rehren subfactors arising from $\alpha$-induction}\RIMS{37}
%{1--35} {(2001)}.

%\bibx{BrVi} {Brugui\`eres, A., Virelizier, A.} {On the center of fusion categories.} {} {} {}

%\bibx{CMS} {Calegari, F., Morrison, S., Snyder, N.} {Cyclotomic integers, fusion categories,
%and subfactors.} \CMP{} \textbf{303}, {845--896} {(2011)}.%; arXiv: 1004.0665.

\bibx{Cardy} {Cardy, J. L.}
 {Conformal invariance and the Yang-Lee edge singularity in two dimensions.} Phys. Rev. Lett. \textbf{54}  1354-1356 (1985)

%\bibx{CS} {Conway, J.H., Sloane, N.J.A.} {Low-dimensional lattices. I: Quadratic forms of small 
%determinant.} {Proc. R. Soc. Lond.} \textbf{A418},  17--41 (1988).

%\bibx{CRT} {R. Coquereaux, R. Rais, E.H. Tahri} {Exceptional quantum subgroups for the rank two Lie algebras B2 and G2.} J. Math. Phys. \textbf{51}:092302, (2010);  arXiv:1001.5416. 


%\bibx{gancoste} {Coste, A., Gannon, T.}
%{Remarks on Galois symmetry in rational conformal field theories.}
%{Phys. Lett.} \textbf{B323}, {316--321} {(1994)}.


%\bibx{CGR} {Coste, A., Gannon, T., Ruelle, P.} {Finite group modular data.}
%Nucl. Phys. \textbf{B581}, 679--717 (2000).

\bibx{CKLW} {Carpi, S., Kawahigashi, Y.,   Longo, R.,   Weiner, M.} {From vertex operator algebras to conformal nets and back.} arXiv:1503.01260.

\bibx{CrRi} {Creutzig, T., Ridout, D.}
{Logarithmic conformal field theory: beyond an introduction.}
J. Phys. A 46 (2013), no. 49, 494006, 72 pp. 

\bibx{DaRu} {Davydov, A., Runkel, I.} {A braided monoidal category for symplectic fermions.}
In: Symmetries and Groups in Contemporary Physics.  Proceedings, Tianjin 2012 (World Scientific,  Hackensack 2013)
399--404.% arXiv:1301.1996.

\bibx{DMS} {Di Francesco, P., Mathieu, P., S\'en\'echal, D.} {Conformal Field Theory.} (Graduate Texts in Contemporary Physics, Springer New York) 1997.

%\bibx{dpr}{Dijkgraaf, R., Pasquier, V., Roche, P.}
%{Quasi-quantum groups related to orbifolds models.}
%{Modern quantum field theory (Bombay, 1990), 375--383,
%World Sci.\ Publishing, River Edge, 1991}.

%\bibx{DMNO} {Davydov, A., M\"uger, M., Nikshych, D., Ostrik, V.} {The Witt group of non-degenerate
%braided fusion categories.} arXiv: 1009.2117v2.

%\bibx{DaRu} {Davydov, A., Runkel, I.} {$\bbZ/2\bbZ$-extensions of Hopf algebra module
%categories by their base categories.} arXiv:1207.3611v1.

%\bibx{DV3} {Dijkgraaf, R., Vafa, C., Verlinde, E., Verlinde, H.} {The operator algebra of
%orbifold models.} Commun. Math. Phys. \textbf{123}, 485--526 (1989).

%\bibx{DW} {Dijkgraaf, R., Witten, E.} 
%{Topological gauge theories and group cohomology.} 
%{Commun. Math. Phys.} \textbf{129}, 393--429 (1990). 

%\bibx{DM} {Dong, C., Mason, G.} {Integrability of $C_2$-cofinite vertex operator algebras.}
 %Int. Math. Res. Not.  2006, Art. ID 80468, 15 pp; {arXiv: math/0601569}.


%\bibx{DoGe} {Dovgard, R., Gepner, D.}
%{Conformal field theories with a low number of primary fields.} 
%J. Phys. A, Math. Theor. \textbf{42}, 304009, 8 pp. (2009).

    
%\bibx{Ad} {Adams, J. F.} {The fundamental representations of $E_8$.} In:
%Piccinini, R., Sjerve, D. (eds.)
%Algebraic topology in honor of Peter Hilton. Proceedings, St. John's 1983,
%Contemp. Math. 37. Providence: American Mathematical Society, 1985, pp. 1--10.




%\bibx{Ati} {Atiyah, M. F.} {Bott periodicity and the index of elliptic operators.} {Quart. J. Math. Oxford Ser.} 19  (1968), 113--140.


%\bibx{Ati} {Atiyah, M.} {K-theory past and present.} Sitzungsberichte der
%Berliner Math. Gesellschaft, Berlin: Berliner Math. Gesellschaft,
%2001, pp. 411--417; math.KT/0012213. 

%\bibx{AH} {Atiyah, M., Hopkins, M.} {A variant of $K$-theory: $K_{\pm}$.} In:
%Tillman, U. (ed.) Topology, geometry and quantum field theory.
%Proceedings, Oxford 2002, London Math. Soc. 
%Lecture Note Ser. 308. Cambridge: Cambridge University Press, 2004, pp. 5--17. 

%\bibx{AS} {Atiyah, M. F., Segal, G.} {Twisted K-theory.} 
%Ukr. Mat. Visn. 1, 287--330 (2004); translation in 
%Ukr. Math. Bull. 1 (2004), 291--334. 

%\bibx{BaBo} {Bais, F. A., Bouwknegt, P. G.} {A classification of subgroup
%truncations of the bosonic string.} Nucl. Phys. B276, 561--570 (1987).

%\bibx{BB} {Blackadar, B.} {K-theory for operator algebras.}
%Math. Sci. Res. Inst. Publ. 5. New York: 
%Springer-Verlag, 1986. 

%\bibx{BHW} {Baez, J. C., Hoffnung, A. E., Walker, C. D.}
%{Higher dimensional algebra VII: Groupoidification.} 
%{Theory Appl. Categ.} 24 (2010), 489--553. 

%\bibx{Bant} {Bantay, P.} {Orbifolds and Hopf algebras.} {Phys. Lett.} B245 (1990), 477--479.

%\bibx{Ban} {Bantay, P.} {Permutation orbifolds and their applications.} In:
%Vertex Operator Algebras in Mathematics and Physics. Proceedings, Toronto 2000.
%Fields Institute Communications 39. Providence: American Mathematical Society,
%2003, pp. 13--23.

%\bibx{BeGa} {Beltaos, E.} {Fixed points and fusion rings.} {JHEP to appear,}
 %{arXiv:1111.3099} 

%\bibx{BFS} {Birke, L., Fuchs, J., Schweigert, C.} {Symmetry breaking boundary conditions
%and WZW orbifolds.} Adv. Theor. Math. Phys. 3, 671--726 (1999).

%\bibx{BauCon} {Baum, P., Connes, A.} {Chern character for discrete groups.} {F\^ete of Topology.} {(Academic Press, 1988)}, 163--232.

%\bibx{Black} {Blackadar, B.} {K-theory for Operator Algebras.} {Math. Sci. Res. Inst. Publ., \textbf{5}} (Springer-Verlag,
%New York 1996).

%\bibx{BE1}
%{B\"ockenhauer, J., Evans, D. E.}
%{Modular invariants, graphs and alpha-induction for
%nets of subfactors. I.} {\CMP \,{197}, {361--386} {(1998)}; II. \CMP \, {200},
%{57--103} {(1999)}; III. \CMP \, {205}, {183--228} {(1999)}}. 

%\bibx{BE4}
%{B\"ockenhauer, J., Evans, D. E.}
%{Modular invariants from subfactors: Type I coupling
%matrices and intermediate subfactors.} \CMP\,{213}, {267--289} {(2000)}.


%\bibx{BE5}
%{B\"ockenhauer, J., Evans, D. E.} {Modular invariants from
%subfactors.} In: Coquereaux, R., Garcia, A., Trinchero, R. (eds.)
%{{Quantum symmetries in theoretical physics
%and mathematics}. Proceedings, Bariloche 2000.  Contemp. Math. 294. 
%Providence: American Mathematical Society, 2002, pp. 95--131}.

%\bibx{BE6} {B\"ockenhauer, J., Evans, D. E.} {Modular invariants and subfactors.} In:
%Longo, R. (ed.) {{Mathematical physics in mathematics and physics. Quantum and operator 
%algebraic aspects}. Proceedings, Siena 2000. Providence: American Mathematical Society. Fields 
%Inst. Commun. 30, 2001, pp. 11--37}.  

%\bibx{BEK1}
%{B\"ockenhauer, J., Evans, D. E., Kawahigashi, Y.}
%{On $\alpha$-induction, chiral generators and modular
%invariants for subfactors.} \CMP\,{208}, {429--487} {(1999)}.

%\bibx{BEK2}
%{B\"ockenhauer, J., Evans, D. E., Kawahigashi, Y.}
%{Chiral structure of modular invariants for subfactors.}
%\CMP\,{210}, {733--784} {(2000)}.

%\bibx{Bot} {Bott, R.} {Introduction to equivariant cohomology.} {In: 
%DeWitt-Morette, C., Zuber, J.-B. (eds.)} 
%Quantum Field Theory: Perspective and Prospective. 
%Dordrecht: Kluwer Academic Publications, 1999, pp. 35--58.

%\bibx{BN} {Bouwknegt, P., Nahm, W.} {Realizations of the exceptional modular
%invariant $A^{(1)}_1$ partition functions.} Phys. Lett. B184, 359--362
%(1987).

%\bibx{BDR} {Bouwknegt, P., Dawson, P., Ridout, D.} {D-branes on group manifolds
%and fusion rings.} J. High Energy Phys. 12, 065 (2002).

%\bibx{Bor} {Borel, A.} {Homology and cohomology of compact connected Lie 
%groups.} Proc. Nat. Acad. Sci. U.S.A.  39,  1142--1146 (1953).

%\bibx{BS} {Braun, V., Sch\"afer-Nameki, S.} {Supersymmetric WZW models and 
%twisted K-theory of SO(3).} Adv. Theor. Math. Phys. 12, 217--242 (2008);
%{hep-th/0403287}.

%\bibx{Bra} {Braun, V.} {Twisted K-theory of Lie groups.} J. High Energy Phys. 03, 029 (2004).

%\bibx{BMRS} {Brodzki, J., Mathai, V., Rosenberg, J., Szabo, R. J.}
%{D-Branes, RR-Fields and duality on noncommutative manifolds.}
%Commun. Math. Phys. 277 (2008), 643--706.

%\bibx{BMRS2} {Brodzki, J., Mathai, V., Rosenberg, J., Szabo, R. J.} {Noncommutative
%correspondences, duality and D-branes in bivariant K-theory.} {Adv. Theor. Math.\ Phys.} 13 (2009), 497--552.%{arXiv: 0708.2648v1} {}

%\bibx{Br} {Brown, E. H.} {The cohomology of BSO$_n$ and BO$_n$ with integer
% coefficients.} {Proc. Amer. Math. Soc.} 85, 283--288 (1982).

%\bibx{Brown} {Brown, K, S.}
%{Cohomology of groups.}
%{Graduate Texts in Mathematics, 87. Springer-Verlag, New York, 1994.}

%\bibx{Bry} {Brylinski, J.-L.} {A correspondence dual to McKay's.} alg-geom/9612003.

%\bibx{BrZ} {Brylinski, J.-L., Zhang, B.} {Equivariant K-theory of compact connected
%Lie groups.} K-Theory 20, 23--36 (2000).

%\bibx{CIZ} {Cappelli, A., Itzykson, C., Zuber, J.-B.} {The $A$-$D$-$E$
%classification of minimal and $A_1^{(1)}$ conformal invariant 
%theories.} Commun. Math. Phys. {113}, 1--26 {(1987)}. 


%\bibx{CW} {Carey, A. L., Wang, B.-L.} {Thom isomorphism and push-forward
%map in twisted K-theory.} math.KT/0507414.

%\bibx{CCH} {Carpi, S., Conti, R., Hillier, R.} {Conformal nets and KK-theory.} {Ann. Funct. Anal.}
%{4} {(2013)}, 11--17.

%\bibx{Connes book} {Connes, A.} {Noncommutative geometry.} San Diego: Academic 
%Press, 1994.

%\bibx{CoSk} {Connes, A., Skandalis, G.} {The longitudinal index theorem for foliations.} {Publ. Res. Inst. Math. Sci.} {20 (6)
%(1984), 1139Ð1183.}

%\bibx{CGR} {Coste, A., Gannon, T., Ruelle, P.} {Finite group modular data.}
%Nucl. Phys. B581 (2000), 679--717.

%\bibx{CM} {Coxeter, H. S. M., Moser, W. O. J.} {Generators and Relations for
%Discrete Groups (4th ed.).} Berlin-Heidelberg: Springer, 1984.

%\bibx{CST} {Cuntz, J., Skandalis, G., Tsygan, B.} {Cyclic Homology in Non-Commutative Geometry.} (Encycl. Math. Sci. 121, Springer 2004).

%\bibx{DW} {Dijkgraaf, R., Witten, E.}  {Topological gauge theories and group cohomology.} {Commun. Math.
%Phys.} 129 (1990), 393--429.

%\bibx{DV}
%{Dijkgraaf, R., Verlinde, E.}
%{Modular invariance and the fusion algebras.}
%\nupp {5B}, 87--97 (1988).

%\bibx{DoMa} {Dong, C., Mason, G.} {Quantum Galois theory for compact Lie groups.}
%{J. Algebra} {214} (1999), 92--102. 


%\bibx{Doug} {Douglas, C. L.} {On the twisted $K$-homology of simple Lie groups.}
%Topology 45, 955--988 (2006).

%\bibx{Doug2} {Douglas, C. L.} {Fusion rings of loop group representations.} arxiv: 0901.0391v1.

%\bibx{DK} {Donovan, P., Karoubi, M.} {Graded Brauer groups and $K$-theory with local 
%coefficients.} Inst. Hautes \'Etudes Sci. Publ. Math. 38, 5--25 (1970).

%\bibx{DJ} {Dunbar, D. C., Joshi, K. G.} {Characters for coset conformal
%field theories and maverick examples.} {Int. J. Mod. Phys.} {A8} 4103--4121 (1993). 

%\bibx{EEK} {Echterhoff, S., Emerson, H., Kim, H. J.} {KK-theoretic duality
%for proper twisted actions.} { Math. Ann. 340 (2008), 839--873.} % math/0610044.


%\bibx{EM} {Emerson, H., Meyer, R.} {Bivariant K-theory via correspondences.} {Adv. Math.} 225
%(2010), 2883--2919.

%\bibx{EGO} {Etingof, P., Gelaki, S., Ostrik, V.} {Classification of fusion categories of dimension $pq$.}
%{Int. Math. Res. Notices 2004,} no. 57, 3041--3056 (2004).

\bibx{ENO} {Etingof, P., Nikshych, D., Ostrik, V.} {On fusion categories.} 
 {Ann. Math.} \textbf{162}, {581--642} {(2005)}.


%\biba{E1}
%{Evans, D.E.} {Fusion rules of modular invariants}\RMP{14}
%{709--732}{(2002)}.

%\bibx{Ev2}
%{Evans, D. E.}
%{Critical phenomena, modular invariants and operator algebras.} In:
%{Operator algebras and mathematical physics (Constan\c ta 2001).
%Cuntz, J., Elliott, G. A., Stratila, S. et al (eds.)
% The Theta Foundation, Bucharest 2003, pp. 89--113}.

%\bibx{E2} {Evans, D.E.} {From Ising to Haagerup.} Markov Processes Relat. Fields \textbf{13},
%267--287 (2007).

%\bibx{Ev} {Evans, D.E.} {Twisted K-theory and modular invariants: I. Quantum
%doubles of finite groups.} In: Bratteli, O., Neshveyev, S., Skau, C. (eds.)
%Operator Algebras: The Abel Symposium 2004. Springer,   Berlin-Heidelberg 2006,
%pp. 117--144.


%\bibx{E1}
%{Evans, D. E.}
%{Fusion rules of modular invariants.} \RMP{} {14}, {709--732} (2002).

%\bibx{E2}
%{Evans, D. E.}
%{Critical phenomena, modular invariants and operator algebras.} In:
%Cuntz, J., Elliott, G. A., Stratila, S. et al (eds.)
%{Operator Algebras and Mathematical Physics. Proceedings, Constan\c ta 2001.
%Bucharest: The Theta Foundation, 2003, pp. 89--113}.

%\bibx{E3} {Evans, D. E.} {Modular invariant partition functions in statistical
%mechanics, conformal field theory and their realisation by subfactors.} In:
%Zambrini, J.-C. (ed.)
%XIVth International Congress on Mathematical Physics. Proceedings, Lisbon 
%2003. Singapore: World Scientific, 2005, pp. 464--475.

%\bibx{Ev} {Evans, D. E.} {Twisted K-theory and modular invariants: I Quantum
%doubles of finite groups.} In: Bratteli, O., Neshveyev, S., Skau, C. (eds.)
%Operator Algebras: The Abel Symposium 2004. Berlin-Heidelberg: Springer, 2006,
%pp. 117--144.

%\bibx{EG1} {Evans, D. E., Gannon, T.}
%{Modular invariants and twisted equivariant K-theory.} {Commun. Number Th. 
%Phys.} {3} (2009), 209--296.

%\bibx{EG3} {Evans, D. E., Gannon, T.}
%{Modular invariants and twisted equivariant K-theory II: Dynkin diagram symmetries.} {J. K-Theory} 12 (2013), 273--330.


\bibx{EG2} {Evans, D. E., Gannon, T.}
{The exoticness and realisability of twisted Haagerup-Izumi modular data.} 
Commun. Math. Phys. {\bf 307} (2011), 463--512. %{arXiv:1006.1326.}


\bibx{EGm} {Evans, D. E., Gannon, T.} {The search for the exotic -- subfactors and conformal field theories.}
In: Progress in Analysis: Proceedings of the 8th Congress of the ISAAC (Moscow, 2011), Vol. 1, Burenkov et al  (eds) (People's Friendship University, Moscow 2012) pp.8--25.

\bibx{EGn} {Evans, D. E., Gannon, T.} {Generalised Haagerup and $E_6$ subfactors and their modular data.}
April 2012.

\bibx{EG4} {Evans, D. E., Gannon, T.}
{Near-group categories and their doubles.} 
Adv. Math.  255 (2014), 586--640.%{arXiv:1006.1326.}

%\bibx{EvGould} 
%{Evans, D. E., Gould, J. D.} 
%{Dimension groups and embeddings of graph algebras.}
%{International J. Math.} {5}, 291--327 (1994).

\bibx{EKaw} {Evans, D. E., Kawahigashi, Y.} {Quantum Symmetries on
Operator Algebras.} {Oxford: Oxford University Press, 1998}.

%\bibx{EP1}
%{Evans, D. E., Pinto, P. T.}
%{Subfactor realisation of modular invariants.} {\CMP} {237} 
%{(2003)}, {309--363}.

%\bibx{EvPui} {Evans, D. E., Pugh, M.} {Ocneanu cells and Boltzmann weights for 
%the $SU(3)$ $\mathcal{ADE}$ graphs.} M\"{u}nster J. Math. {2}, 95--142 (2009);
 %arXiv:0906.4307.

%\bibx{EvPuii} {Evans, D. E., Pugh, M.}  {$SU(3)$-Goodman-de la Harpe-Jones 
%subfactors and the realisation of $SU(3)$ modular invariants.} 
%Rev. Math. Phys. {21}, 877--928 (2009); arXiv:0906.4252.

%\bibx{FFFS} {Felder, G., Fr\"ohlich, J., Fuchs, J., Schweigert, C.} {The  
%geometry of WZW branes.} J. Geom. Phys. 34, 162--190 (2000); hep-th/9909030.

%\bibx{FGK} {Felder, G., Gawedzki, K., Kupiainen, A.} {Spectra of Wess-Zumino-Witten
%models with arbitrary simple groups.} Commun. Math. Phys. 117, (1988) 127--158.

%\bibx{FRS2} {Fredenhagen, K., Rehren, K.-H., Schroer, B.}
%{Superselection sectors with braid group statistics and exchange
%algebras. II.} {\RMP} \textbf{Special issue}, {113--157} {(1992)}.

%\bibx{Fr} {Fredenhagen, S.} {D-brane charges on SO(3).} {J. High Energy Phys.} 
%11, 082 (2004); {hep-th/0404017.}

%\bibx{FGM} {Fredenhagen, S., Gaberdiel, M. R.,  Mettler, T.} {Charges of twisted branes: the 
%exceptional cases.} {J. High Energy Phys.} 05, 058 (2005); hep-th/0504007. 


%\bibx{FHT} {Freed, D. S., Hopkins, M. J., Teleman, C.} 
%{Twisted equivariant K-theory
   % with complex coefficients.} {J. Topol. 1 (2008), 16--44.}%; math.AT/0206257.

%\bibx{FHTi} {Freed, D. S., Hopkins, M. J., Teleman, C.} 
%{Loop groups and twisted
%K-theory I.} {  J. Topol.} 4 (2011),  737--798.%; math.AT/0711.1906. 

%\bibx{FHTii} {Freed, D. S., Hopkins, M. J., Teleman, C.} {Loop groups and
%twisted K-theory II.} {J. Amer. Math. Soc.} 26 (2013), 595--644.%{math.AT/0511232 v.2}.

%\bibx{FHTlg} {Freed, D. S., Hopkins, M. J., Teleman, C.} {Loop groups and
%twisted K-theory III.} {Ann. of Math. (2) 174 (2011), 947--1007.} %{math.AT/0312155 v.3}.

%\bibx{FHLT} {Freed, D. S., Hopkins, M. J., Lurie, J., Teleman, C.} {Topological quantum field
%theories from compact Lie groups.} arXiv:/ 0905.0731.

%\bibx{FFRS}  {Fr\"ohlich, J., Fuchs, J., Runkel, I., Schweigert, C.} {Defect lines, dualities,
%and generalized orbifolds.} {XVIth International Congress on Mathematical Physics, 608Ð613, World Sci. Publ., Hackensack, NJ, 2010;} arXiv: math-ph/0909.5013.


%\bibx{FG1} {Fr\"ohlich, J., Gabbiani, F.} {Braid statistics in
%local quantum theory.} {\RMP} {2}, {251--353}  {(1990)}.

%\bibx{Fu} {Fuchs, J.} {Simple WZW currents.} Commun. Math. Phys. 136, 345--356 (1991).

%\bibx{FSS} {Fuchs, J., Schellekens, A. N., Schweigert, C.} {From Dynkin
%diagram symmetries to fixed point structures.} Commun. Math. Phys. 180, 39--97 (1996).

%\bibx{GPRV} {Gaberdiel, M.R., Persson, D., Ronellenfitsch, H., Volpato, R.} {Generalised Mathieu Moonshine.} {Commun. Number Th. Phys.} {7} (2013), 145--223.%{arXiv:1211.7074}.

%\bibx{GaGa} {Gaberdiel, M. R., Gannon, T.} {Boundary states for WZW models.}
%Nucl. Phys. B639, 471--501 (2002).

%\bibx{GaGa1} {Gaberdiel, M. R., Gannon, T.} {The charges of a twisted brane.}
%J. High Energy Phys. 01, 018 (2004); hep-th/0311242. 


%\bibx{GaGa2} {Gaberdiel, M. R., Gannon, T.} {D-brane charges on non-simply 
%connected groups.}  J. High Energy Phys.  2004,  no. 4, 030, 27 pp.

%\bibx{GaGa3} {Gaberdiel, M. R., Gannon, T.} {Twisted brane charges for 
%non-simply connected groups.}  J. High Energy Phys.  2007,  no. 1, 035, 30 pp.


%\bibx{GGR} {Gaberdiel, M. R., Gannon, T., Roggenkamp, D.} {The D-branes of 
%SU(n).} JHEP 0407, 015 (2004).


%\bibx{Ganmd} {Gannon, T.} {Modular data: the algebraic combinatorics of
%rational conformal field theory.} J. Alg. Combin. 22  (2005), 211--250.


%\bibx{G2} {Gannon, T.} {The classification of affine $\SUd$
%modular invariants.} {\CMP} {161}, 233--264 (1994).
%\bibx{Gannbook}{gannon, T.} {Moonshine Beyond the Monster:  The bridge connecting
% algebra, modular forms and physics,} Cambridge University Press, 2006.


%\bibx{GaVa} {Gannon, T., Vasudevan, M.} {Charges of exceptionally twisted branes.}
%{J. High Energy Phys.} 07, 035 (2005); hep-th/0504006v4.

\bibx{Ganun} {Gannon, T.} {Comments on nonunitary conformal field theories.} {Nucl. Phys. B670,} 335--358 (2003).

\bibx{vvmf} {Gannon, T.} {The theory of vector-valued modular forms for the modular group.}  In: Conformal Field
Theory, Automorphic Forms and Related Topics.  Eds: W Kohnen, R Weissauer.  (Springer Verlag, 2014) 247--286.

%\bibx{GaWa} {Gannon, T., Walton, M. A.} {On the classification of diagonal
%coset modular invariants.} {Commun. Math. Phys.} {173}, {175--197} (1995).  

%\bibx{GaWa2} {Gannon, T., Walton, M. A.} {On fusion algebras and modular
%matrices.} {Commun. Math. Phys.} 206, 1--22 (1999).

%\bibx{Gnt} {Ganter, N.} {Hecke operators in equivariant elliptic cohomology and
% generalized Moonshine.} In:  Groups and Symmetries. Proceedings, Montreal 2007.
% CRM Proc. 
%Lecture Notes 47. Providence: Amer. Math. Soc., 2009, pp. 173--209.

%\bibx{Gep} {Gepner, D.} {Fusion rings and geometry.} Commun. Math. Phys.
%{141}, 381--411 (1991).

%\bibx{GKO} {Goddard, P., Kent, A., Olive, D.} {Unitary representations of Virasoro and
%super-Virasoro algebras.} {Commun. Math. Phys.} {103}, {105--119} (1986).

%\bibx{GL2} {Guido, D., Longo, R.} {The conformal spin and
%statistics theorem.} {\CMP} {181}, {11--35} {(1996)}.

%\bibitem{haag:1992}
%Haag, R. {\it Local Quantum Physics}, Springer, Berlin, 1992.

%\bibx{UHaag} {Haag, U.} {On $\bbZ/2\bbZ$ graded KK-theory and its relation with the graded
%Ext-functor.} J. Operator Theory {42}, {3--36} (1999). 

%\bibx{HMT} {Hanaki, A., Miyamoto, M., Tambara, D.} {Quantum Galois theory for finite groups.} {Duke Math. J.} {97} (1999), 541--544.


%\bibx{HH} {Hanany, A., He, Y.-H.} {Non-abelian finite gauge theories.} JHEP 
%9902, 013 (1999); hep-th/9811183.

%\bibx{HV} {Hanany, A., Vegh, D.} {Quivers, tilings, branes and rhombi.} JHEP 
%2007, no.10, 029 (2007); hep-th/0511063.

%\bibx{Ha} {Hatcher, A.} {Algebraic Topology.} {Cambridge: Cambridge University 
%Press, 2002}.

%\bibx{Ht} {Hatcher, A.} {Spectral Sequences in Algebraic Topology.} (available
%electronically from Hatcher's homepage //www.math.cornell.edu/\~{}hatcher).

%\bibx{HS} {He, Y.-H., Song, J. S.} {Of McKay correspondence, non-linear sigma
%model and conformal field theory.} Adv. Theor. Math. Phys. 4, 747--790 (2000);
%hep-th/9903056.


%\bibx{HG} {Higson, N., Guentner, E.} {Group $C^{}$-algebras and $K$-theory.}
%In: Doplicher, S., Longo, R. (eds.) Noncommutative geometry. 
%Lecture Notes in Math. 1831. Berlin: Springer, 2004, pp. 137--251. 

%\bibx{HR} {Higson, N., Roe, J.} {Analytic K-homology.} Oxford: Oxford 
%University Press, 2000.

%\bibx{Hugh} {Hughes, N. J. S.} {The use of bilinear mappings in the classification of
%groups of class 2.} {Proc. AMS} 2 (1951), 742--747.

%\bibx{Isa} {Isaacs, I. M.} {Character theory of finite groups.} {New York: Dover, 1994.}


%\bibx{IN} {Ito, Y., Nakamura, I.} {Hilbert schemes and simple singularities.}
%In: Hulek, K., Catanese, F., Peters, C., Reid, M. (eds.)  New trends in 
%algebraic geometry. Proceedings, Warwick 1996. London Math. Soc.
% Lecture Note Ser. 264. Cambridge: Cambridge University Press, 1999, pp. 151--233.

%\bibx{IR} {Ito, Y., Reid, M.} {The McKay correspondence for finite subgroups
%of SL$(3,\bbC)$.} In: Andreatta, M., Peternell, T. (eds.) 
%Higher-dimensional complex varieties. Proceedings, Trento 1994.
%Berlin: de Gruyter, 1996, pp. 221--240.

%\bibx{Iz} {Izumi, M.}
%{The structure of sectors associated with Longo-Rehren inclusions, I.
%General theory.} \CMP{213}, {127--179} {(2000)}.

%\bibx{JW} {Jeffrey, L. C., Weitsman, J.} {Bohr-Sommerfeld orbits in the 
%moduli space of flat connections and the Verlinde dimension formula.}  
%Commun. Math. Phys.  150   no. 3,   593--630 (1992).

%\bibx{J1} {Jones, V. F. R.} 
%{An invariant for group actions.} In: de la Harpe, P. (ed.) 
%Alg\`ebres d'op\`erateurs. Proceedings, Les Plans-sur-Bex 1978.
%Lecture Notes in Math. 725.  Berlin: Springer, 1979, pp. 237--253. 

%\bibx{kac} {Kac, V.G.} {Infinite-dimensional Lie Algebras, 3rd edn.} Cambridge: Cambridge University
%Press, 1990.


%\bibx{Jon1} {Jones, V. F. R.}  {An invariant for group actions.} {Alg\`ebres d'op\'erateurs (S\'em., Les Plans-sur-Bex, 1978),} pp. 237--253, {Lecture Notes in Math., 725, Springer, Berlin, 1979.} 

%\bibx{Jon2} {Jones, V. F. R.}  {Actions of finite groups on the hyperfinite type II$_1$ factor.} {Mem. Amer. Math. Soc.} 28 (1980), no. 237.



%\bibx{Kar1} {Karoubi, M.} {Twisted K-theory, old and new.}  In:
 %$K$-theory and noncommutative geometry.  EMS Ser. Congr. Rep. Z\"urich: European Math. Soc.,  2008, pp. 117--149; arXiv: 
%math.KT/0701789. 
 
%\bibx{Kar2} {M. Karoubi} {K-theory. An elementary introduction,} math.KT/0602082.

%\bibx{KL} {M. Karoubi, T. Lambre} {K-theory of the norm functor,} 
%math.KT/0701628.

%\bibx{Karp} {Karpilovsky, G.}  {Group representations. Vol. 2.} North-Holland Mathematics Studies, 177. North-Holland Publishing Co., Amsterdam, 1993.

%\bibx{Kas} {Kasparov, G. G.} {The operator $K$-functor and extensions of $C^*$-algebras.} {Izv.\ Akad.\ Nauk.\ SSSR, Ser.\ Mat.} {44 (1980), 571--636.}

%\bibx{GK1} {Kasparov, G. G.} {Topological invariants of elliptic operators. {I}.
 % {$K$}-homology.} Izv. Akad. Nauk SSSR Ser. Mat. {39}, 
 % 796--838 (1975).

%\bibx{GK2} {Kasparov, G. G.}
%{Equivariant {$KK$}-theory and the Novikov conjecture.} Invent.
 % Math. {91}, 147--201 (1988).

%\bibx{GK3} {Kasparov, G. G.} {{$K$}-theory, group {$C\sp *$}-algebras, and higher signatures
  %(conspectus).} In: Ferry, S., Ranicki, A., Rosenberg, J. (eds.)
%Novikov conjectures, index theorems and rigidity, Vol.\ 1. Proceedings, 
%Oberwolfach 1993. London Math. Soc. Lecture Note Ser. 226. Cambridge: Cambridge
  %University Press, 1995, pp. 101--146.

%\bibx{Yasu} {Kawahigashi, Y.} {From operator algebras to superconformal field theory,}
%J. Math. Phys. 51 (2010), 015209. arXiv:1003.2925

%\bibx{Kir} {Kirillov, A. A.} {Lectures on the Orbit Method.} 
%Providence: American Mathematical Society, 2004.

%\bibx{KS}  {Klemm, A., Schmidt, M. G.} {Orbifolds by cyclic permutations of
%tensor product conformal field theories.} Phys. Lett. B245, 53--58 (1990).

%\bibx{KrSch} {Kreuzer, M., Schellekens, A. N.}
%{Simple currents versus orbifolds with discrete torsion -- a complete classification.} 
%{Nucl. Phys.} {B411} (1994),  97--121. 


%\bibx{L2} {Landweber, G. D.} {Harmonic spinors on homogeneous spaces.}
%  Represent. Theory 4, 466--473 (2000).

%\bibx{La} {Landweber, G. D.} {Twisted representation rings and Dirac
%induction.} J. Pure Appl. Algebra 206, 21--54 (2006); {math.RT/0403524}.

%\bibx{MMS} {Maldacena, J., Moore, G., Seiberg, N.} {Geometrical interpretation
%of D-branes in gauged WZW models.} JHEP 2001 no. 7, 046 (2001); hep-th/0105038.


%\bibx{Ln5} {Longo, R.} {Index
% of subfactors and statistics of quantum fields, I.} {\em Commun. Math. Phys.} {\bf 126}, 217--247 (1989).

%\bibx{LR}
%{Longo, R.,  Rehren, J.-H.} {Nets of subfactors, Rev. Math. Phys}. {7} (1995), 567--597.

%\bibx{McK} {McKay, J.} {Graphs, singularities, and finite groups.} In:
%Cooperstein, B., Mason, G. (eds.) The Santa Cruz
%Conference on Finite Groups. Proceedings, Santa Cruz 1979. 
%Proc. Sympos. Pure Math. 37. Providence:
%American Mathematical Society, 1980, pp. 183--186.

%\bibx{M} {Meinrenken, E.}
%{On the quantization of conjugacy classes.} L'Enseignement 
%Math. (2)  55, 33--75 (2009); math.DG/0707.3963.

%\bibx{Mi1} {J. Mickelsson} {Gerbes, (twisted) $K$-theory, and the 
%supersymmetric WZW model,} Infinite dimensional groups and manifolds, 
%93--107, IRMA Lect. Math. Theor. Phys. 5, de Gruyter, Berlin, 2004. 

%\bibx{Mi2} {J. Mickelsson} {Twisted $K$ theory invariants,}
% Lett. Math. Phys. 71 (2005), no. 2, 109--121.

%\bibx{Miy} {Miyamoto, M.} {$C_2$-cofiniteness of cyclic-orbifold models.} {arXiv: 1306.5031v1.}

%\bibx{Mon} {Monnier, S.} {Kondo flow invariants, twisted K-theory and 
%Ramond-Ramond charges.} JHEP 06, 022 (2008); hep-th/0803.1565

%\bibx{Mokt} {Moore, G.} {$K$-theory from a physical perspective.} In: 
%Topology, Geometry and Quantum Field Theory.  London Math. Soc. Lecture Note 
%Ser. 308. Cambridge: Cambridge Univ. Press, 2004, pp. 194--234.

%\bibx{MS2}
%{Moore, G., Seiberg, N.}
%{Naturality in conformal field theory.}
%\nup {B313}, 16-40 (1989).

%\bibx{MW} {Mohrdieck, S., Wendt, R.} {Integral conjugacy classes of compact
%Lie groups.} Manuscripta Math. 113, 531--547 (2004).

%\bibx{Nahm} {Nahm, W.} {Quasi-rational fusion products.} 
 %Internat. J. Modern Phys. B8, 3693--3702  (1994).

%\bibx{Na} {Nahm, W.} {Lie group exponents and $SU(2)$ current algebras.} 
%Commun. Math. Phys. 118, 171--176 (1988).

%\bibx{NeTu} {Neshveyev, S., Tuset, L.} {Hopf algebra equivariant cyclic
%cohomology, $K$-theory, and index formulas.} K-theory 31, 357--378 (2004); 
%math/0304001.

%\bibx{O}
%{Ocneanu, A.} {Paths on Coxeter diagrams: From Platonic solids and
%singularities to minimal models and subfactors.} (Notes recorded by S.\ Goto).
%In: Rajarama Bhat, B.V.\ et al.\ (eds.) Lectures on Operator Theory.
%Providence: American Mathematical Society, 2000, pp. 243--323. 


%\bibx{ocn4}{Ocneanu, A.}
%{The classification of subgroups of quantum ${\mathit{SU}}(N)$.} 
%In: Coquereaux, R., Garcia, A., Trinchero, R. (eds.)
%{{Quantum Symmetries in Theoretical Physics
%and Mathematics}. Proceedings, Bariloche 2000,  Contemp. Math. 294. 
%Providence: American Mathematical Society, 2002, pp. 133--159}.

%\bibx{Ost1} {Ostrik, V.} {Module categories, weak Hopf algebras and modular invariants.} Transform. Groups
%8 (2003), 177--206.

%\bibx{Ost} {Ostrik, V.} {Module categories for quantum doubles of finite groups.}
%Int. Math. Res. Notices 27 (2003), 1507--1520.%; math.QA/0202130.

%\bibx{EMP1} {Parker, E. M.} {The Brauer group of graded $C^{*}$-algebras.} 
%Trans. Amer. Math. Soc. 308, 115--132 (1988).

%\bibx{EMP2} {Parker, E. M.} {Graded continuous trace $C^{*}$-algebras and duality.}
%In: Arveson, W. B. et al (eds.)
%Operator algebras and topology. Proceedings, Craiova, 1989.
%Pitman Res. Notes Math. Ser. 270. Harlow:  
% Longman Sci. Tech., 1992, pp. 130--145. 

%\bibx{PZ} {Petkova, V. B., Zuber, J.-B.} {The many faces of Ocneanu cells.}
%Nucl. Phys. B603, 449--496 (2001); hep-th/0101151.

%\bibx{Pimsner} {Pimsner, M.} {Lecture at Odense Meeting,} April  2007.

%\bibx{Pit} {Pittie, H. V.} {Homogeneous vector bundles on homogeneous
%spaces.} Topology 11, 199--203 (1972).

%\bibx{RW} {Raeburn, I., Williams, D. P.} {Morita equivalence and continuous 
%trace $C^{*}$-algebras.}
% Mathematical Surveys and Monographs 60. Providence:
%American Mathematical Society, 1998.

%\bibitem{rehren:2000}
% Rehren, K.-H. Chiral observables and modular invariants, Comm. Math. Phys. {\bf 208} (2000), 689--712.

%\bibx{Ren} {Rennie, A.} {Smoothness and locality for nonunital spectral
%triples.} K-theory 28, 127--165 (2003).

%\bibx{Ros1} {Rosenberg, J.} {Homological invariants of extensions of $C^{*} 
%$-algebras.} In: Kadison, R. V. (ed.)
%Operator algebras and applications, Part 1. Proceedings, Kingston
%1980. Proc. Sympos. Pure Math. 38. Providence: American Mathematical Society, 
%Providence, 1982, pp. 35--75.

%\bibx{Ros2} {Rosenberg, J.} {Continuous-trace algebras from the bundle 
%theoretic point of view.}
% J. Austral. Math. Soc. Ser. A 47, 368--381 (1989).

%\bibx{RS} {Rosenberg, J., Schochet, C.} {The K\"unneth theorem and the
%universal coefficient theorem for equivariant K-theory and KK-theory.}
%Mem. Amer. Math. Soc. 62 (1986).

%\bibx{RFS} {Runkel, I., Fjelstad, J., Fuchs, J., Schweigert, C.} {Topological
%and conformal field theory as Frobenius algebras.} In: Davydov, A. et al (eds.)
%Categories in Algebra, Geometry and Mathematical Physics.
%Proceedings, Sydney 2005. Contemp. Math. 431. Providence: American 
%Mathematical Society, 2007, pp. 225--248; math.CT/0512076.

%\bibx{SN1} {Sch\"afer-Nameki, S.} {D-branes in $N=2$ coset models and twisted
% equivariant K-theory.} hep-th/0308058. 

%\bibx{SN2} {Sch\"afer-Nameki, S.} {$K$-theoretical boundary rings in $N=2$ 
%coset models.} Nuclear Phys. B706, 531--548 (2005); hep-th/0408060.



%\bibx{ScWa} {Schellekens, A. N., Warner, N. P.} {Conformal subalgebras of
%Kac-Moody algebras.} Phys. Rev. D34, 3092--3096 (1986).

%\bibx{SY} {Schellekens, A. N., Yankielowicz, S.} {Simple currents, modular invariants and 
  %fixed  points.} Int. J. Mod. Phys. 5A, 2903--2952 (1990).

%\bibx{Se2} {Segal, G.} {Classifying spaces and spectral sequences.} Inst.
%Hautes \'Etudes Sci. Publ. Math. No. 34, 105--112 (1968).

%\bibx{Se} {Segal, G.} {Equivariant K-theory.} Inst. Hautes \'Etudes Sci. Publ. 
%Math. No. 34 (1968), 129--151.

%\bibx{Seg} { {Segal, G.}} 
%{The definition of conformal field theory.}
%{Differential Geometric Methods in Theoretical Physics (Como, 1987)} 
%(Academic Press, Boston 1988) 165--171.


%\bibx{Sim} {Simon, B.} {Representations of Finite and Compact Groups.}
%Providence: American Mathematical Society, 1996.

%\bibx{Sta} {Stanciu, S.} {An illustrated guide to D-branes in SU(3).} hep-th/0111221.

%\bibx{Suth} 
%{Sutherland, C.} 
%{Cohomology and extensions of von Neumann algebras. I, 
%II.} {Publ.\ RIMS.\ Kyoto Univ.} {16}, 105--133, {135--174} {(1980)}. 

%\bibx{Tro} {Trout, J.} {On graded K-theory, elliptic operators and the
%functional calculus.} Illinois J. Math. 44, 294--309 (2000); {math.OA/9906106}.

%\bibx{Tu1} {Tu, J. L.} {Non-Hausdorff groupoids, proper actions and K-theory.} {Doc. Math.} 9 (2004), 565--597. 

%\bibx{Tu} {Tu, J. L.} {Twisted $K$-theory and Poincar\'e duality.}
%{Trans. Amer. Math. Soc.} { 361 (2009), 1269--1278;} math.KT/0609556.

%\bibx{TX} {Tu, J. L., Xu, P.} {The ring structure for equivariant twisted
%K-theory.} J. Reine Angew. Math. 635, 97--148 (2009); math.KT/0604160.


%\bibx{tur}
%{Turaev, V.G.}{Quantum Invariants of Knots and 3-manifolds.} {
%{de Gruyter Studies in Mathematics}, vol 18. Berlin: Walter de
%Gruyter, 1994}.


%\bibx{Vrr} {Verrill, R. W.} {Positive energy representations of $L^\sigma SU(2r)$ and orbifold
%fusions.} {Ph.D. thesis,} U Cambridge (2002). 

%\bibx{Vrst} {Verstegen, D.} {Conformal embeddings, rank-level duality and 
%exceptional modular invariants.}  Commun. Math. Phys.  137  (1991),  no. 3, 
%567--586. 

%\bibx{Wng} {Wang, W.} {Equivariant $K$-theory, wreath products, and Heisenberg 
%algebra.}  Duke Math. J.  103  (2000),  no. 1, 1--23. 



%\bibx{Wassthes} {Wassermann, A.} {PhD thesis,} U Pennsylvania (1981).

%\bibx{W} {Wassermann, A.} {Operator algebras and conformal field theory
%{\rm III}: fusion of positive representations of $\rm{LSU}(n)$
%using bounded operators.} \Inv{} {133}, {467--538} {(1998)}.

%\bibx{W3} {Wassermann, A.} {Subfactors and Connes fusion for twisted loop 
%groups.} arXiv:1003.2292.


%\bibx{xu}{Xu, F.}
%{New braided endomorphisms from conformal inclusions.} \CMP{} {192},
%{349--403} {(1998)}.


%\bibx{xu2}{Xu, F.}
%{Mirror  extensions of local nets.}  Comm. Math. Phys.  {270},  {835--847}  {(2007)}.



%\bibx{EGh} {Evans, D. E., Gannon, T.}
%{The exoticness and realisability of twisted Haagerup-Izumi modular data.} 
%\CMP{} \textbf{307},  {463--512} {(2011)}; arXiv:1006.1326.

%\bibx{EG1} {Evans, D. E., Gannon, T.}
%{Modular invariants and twisted equivariant K-theory.} {Commun. Number Theory 
%Phys.} \textbf{3}, 209--296 (2009).

%\bibx{EK1} {Evans, D.E., Kawahigashi, Y.}  { Orbifold subfactors from
%Hecke algebras II: quantum double and braiding.} \CMP{}  \textbf{196}, {331--361} {(1998)}; arXiv:
  %funct-an/9702018.


%\bibx{EK} {Evans, D.E., Kawahigashi, Y.} {Quantum Symmetries on
%Operator Algebras.}{Oxford University Press 1998}.




%\biba{EP} {Evans, D.E., Pinto, P.R.}{Subfactor realisation
%of modular invariants}\CMP{237} {309--363} {(2003)}.

%\bibx{EP2}
%{Evans, D.E., Pinto, P.R.} {\sl Modular invariants and their
%fusion rules.} {In: {Advances in Quantum Dynamics} (Mount
%Holyoke, 2002), pp.\ 119--130, Contemp.\ Math.\ {\bf 335}, Amer.\
%Math.\ Soc., Providence, RI, 2003}.

%\bibx{EP3}
%{Evans, D.E., Pinto, P.R.} {Modular invariants
%and the double of the Haagerup subfactor.} In: Advances in
%Operator Algebras and \hbox{Mathematical} Physics (Sinaia 2003). 
%\break  Boca, F.-P.,  Bratteli, O.,  Longo, R.,  Siedentop, H. (eds.). The Theta Foundation,
%Bucharest 2006, pp.67--88.

%\bibx{EP4} {Evans, D.E., Pinto, P.R.} {Subfactor realisation of modular invariants: II.} (in preparation).


%\bibx{FrZh} {Frenkel, I.B., Zhu, Y.} {Vertex operator algebras associated to representations of affine 
%and Virasoro algebras.}  Duke Math. J. \textbf{66}, 123--168 (1992). 

%\bibx{FFRS}  {Fr\"ohlich, J., Fuchs, J., Runkel, I., Schweigert, C.} {Defect lines, dualities,
%and generalised orbifolds.} arXiv: math-ph/0909.5013.

%\biba{frs} {Fuchs, J., Runkel, I., Schweigert, C.}
%{TFT construction of RCFT correlators I: Partition
%functions}\nupb{646}{354--497}{(2002)}.

%\biba{frs2} {Fuchs, J., Runkel, I., Schweigert, C.}
%{TFT construction of RCFT correlators II: Unoriented world
%sheets}\nupb{678}{511--637}{(2004)}.

%\biba{frs3}{Fuchs, J., Runkel, I., Schweigert, C.}
%{TFT construction of RCFT correlators III: Simple currents}
%\nupb{694}{277--353}{(2004)}.

% \biba{GabGan}{Gaberdiel, M.\ R.; Gannon, T.} {Boundary states for WZW models} \nupb{639} {471--501} {(2002)}.

%\bibx{Gan2} {Gannon, T.} {Modular data: the algebraic combinatorics of conformal field
%theory.} {J. Alg. Combin.} \textbf{22}, 211--250 (2005).

%\bibx{Gorth} {Gannon, T.} {The level 2 and 3 modular invariants for the orthogonal algebras.}
%{Canad. Math. J.} \textbf{52}, {503--521} {(2000)}.

%\bibx{GKO} {Goddard, P., Kent, A.,  Olive, D.}
%{Unitary representations of Virasoro and super-Virasoro algebras.}
%{Commun.\ Math.\ Phys.} {\bf 103}, 105--119 (1986).

%\biba{GHJ}  {Goodman, F.M., de la Harpe, P., Jones, V.F.R.} {Coxeter 
%graphs and towers of algebras} {vol.} {14} {Mathematical Sciences Research Institute 
%Publications} (Springer-Verlag, New York 1989).

\bibx{GrSn} {Grossman, P., Snyder, N.} {Quantum subgroups of the Haagerup fusion categories.}
{Commun. Math. Phys.} \textbf{311}, 617--643 (2012).


\bibx{Haag}{Haag, R.}   
%\bibitem{haag:1992}
{\it Local Quantum Physics.} Springer, Berlin, 1992. 

\bibx{Haag} {Haagerup, U.} {Principal graphs of subfactors in the index range $4<[M:N]<3+
\sqrt{3}$.} In: Subfactors. H. Araki et al (eds.). World Scientific, 1994, pp.1--38.

\bibx{HaYa} {Hayashi, T., Yamagami, S.}
{Amenable tensor categories and their realizations as AFD bimodules.}
J. Funct. Anal. {\bf 172}, (2000)  19--75. 

%\bibx{Hn} {Han, R.} {A Construction of the ``2221'' Planar Algebra.} {}{}{}Ph.D. Thesis (University of California, Riverside) (2010); {arXiv:1102.2052v1}.

%\bibx{HS} {Hong, J.H., Szymanski, W.} {Composition of subfactors and twisted bicrossed
%products.} {J. Operator Theory} \textbf{37},  {281--302} {(1997)}.

%\bibx{HRW} {Hong, S.-M., Rowell, E., Wang, Z.} {On exotic modular tensor 
%categories.} Commun. Contemp. Math. \textbf{10}, Suppl. 1, 1049--1074 (2008). 

%\bibx{HoJo} {Horn, R.A., Johnson, C.R.} {Matrix Analysis.} (Cambridge University Press, 1985).{}{}

%\bibx{Hua} {Huang, Y.-Z.} {Vertex operator algebras, the Verlinde conjecture, and modular tensor 
%categories.}  Proc. Natl. Acad. Sci. USA \textbf{102}, 5352--5356 (2005).

\biba{iz0} {Izumi, M.} {Subalgebras of infinite C*-algebras with finite Watatani indices,
 I. Cuntz Algebras} \CMP{155}{157--182}{(1993)}.

% \biba{iz}{Izumi, M.}
%{Subalgebras of infinite $C\sp *$-algebras with finite Watatani indices. II. Cuntz-Krieger algebras}
%\Duke{91} {409--461}{(1998)}.

\biba{iz1}{Izumi, M.}
{The structure of sectors associated with Longo-Rehren inclusions,
I. General theory}\CMP{213}{127--179}{(2000)}.

\biba{iz3}
{Izumi, M.} {The structure of sectors associated with Longo-Rehren
inclusions, II. Examples}\RMP{13}{603--674}{(2001)}.

\bibx{iz4}
{Izumi, M.} {(notes on the Haagerup series).} {January 2012.} {} {}  

%\bibx{iz5} {Izumi, M.} {(notes on near-group categories).} {} {} {}  

%\bibx{IK}
%{Izumi, M., Kawahigashi, Y.}  {Classification of subfactors with the principal graph $D_{n}^{(1)}$.}\JFA{} \textbf{112}, { 257--286} {(1993)}.

%\bibx{IK} {Izumi, M., Kosaki, H.} {Kac algebras arising from composition of subfactors: general theory and classification.} {Mem. Amer. Math. Soc.} \textbf{158} no. 750, (2002).
%{Finite-dimensional Kac algebras arising from certain
%group actions on a factor.} {Intern. Math. Res. Not.} \textbf{1996}, {No. 8}, {357--370} {(1996)}.

%\biba{KI}
%{Izumi, M., Kosaki, H.} {On the subfactor analogue of the second
%cohomology}\RMP{14} {733--757}{(2002)}.

%\bibx{J2} {Jones, V. F. R.} 
%{An invariant for group actions.} In: de la Harpe, P. (ed.) 
%Alg\`ebres d'op\`erateurs (Les Plans-sur-Bex 1978).
%Lecture Notes in Math. 725.  Springer, Berlin 1979, pp. 237--253. 


%\biba{J1}
%{Jones, V.F.R.} {Index for subfactors} \Inv{72} {1--25}{(1983)}.

%\bibx{kac} {Kac, V.G.} {Infinite-dimensional Lie algebras, 3rd edn.} (Cambridge University
%Press, Cambridge 1990).

%\bibx{KacT} {Kac, V.G., Todorov, I.T.} {Affine orbifolds and rational conformal field theory
%extensions of $W_{1+\infty}^*$.} \CMP\  {\textbf{190}}, {57--111} {(1997)}.

%\bibx{KaWa} {Kac, V.G., Wakimoto, M.} {Modular and conformal invariance constraints
%in representation theory of affine algebras.} {Adv. Math.} \textbf{70}, 156--236 (1988).

%\bibx{KLM} {Kawahigashi, Y., Longo, R., M\"uger, M.} {Multi-interval subfactors and modularity of representations in conformal field theory.} 
%{\em Commun. Math. Phys.} {\bf 219}, {631--669} {(2001)};
%arXiv: math.OA/9903104.

%\biba{KY1}
%{Kosaki, H., Munemasa, A., Yamagami, S.} {Irreducible bimodules
%associated with crossed product algebras} \IJM{3} {661--676}
%{(1992)}.


%\bibx{Leav1} {Leavitt, W. G.} {The module type of a ring.} {Trans. AMS} 103 (1962) 113--130.


\bibx{Leav2} {Leavitt, W. G.} {The module type of homomorphic images.} {Duke Math J.} 32 (1965) 305--311.


%\bibx{LeLi} {Lepowski, J., Li, H.} {Introduction to Vertex Operator Algebras and their 
%Representations.} {Boston, Birkh\"auser 2004}.

%\biba{longo2}{Longo, R.}{Index of subfactors and statistics of quantum
%fields, II}\CMP{130}{285--309} {(1990)}.

%\biba{longo}{Longo, R.}{Duality for Hopf algebras and for subfactors, I}
%\CMP{159}{133--150}{(1994)}.

%\bibx{Ln5} {Longo, R.} {Index
% of subfactors and statistics of quantum fields, I.} {\em Commun. Math. Phys.} {\bf 126}, 217--247 (1989).


%\biba{LR}
%{Longo, R., Rehren, K.-H.} {Nets of subfactors}\RMP{7} {567--597}
%{(1995)}.

%\bibx{maclane}{MacLane, S.}
%{Cohomology theory of Abelian groups.
%Proceedings of the International Congress of Mathematics, pp.\ 8--14, 1950.}

%\biba{masuda}
%{Masuda, T.} {An analogue of Longo's canonical endomorphism for
%bimodule theory and its application to asymptotic
%inclusions}\IJM{8} {249--265}{(1997)}.

 
%\bibx{MoSe} {Moore, G., Seiberg, N.} {Lectures on RCFT.} In: Physics, Geometry and Topology, 
%H.  C. Lee (ed.), Nato ASI Series, Vol. 238, 1990, pp.263--361. 


%\biba{MoSn}  {Morrison, S., Snyder, N.}
 %   {Subfactors of index less than 5, part 1: the principal graph odometer} \CMP{{312}} {1--35}
   % {(2012).} %   {arXiv:1007.1730}.

%\bibx{MoPe}  {Morrison, S., Peters, E.} {The little desert? Some subfactors with index in the
%interval $(5,3+\sqrt{5})$;} {arXiv: 1205.2742.} 

%\biba{m1} {Mueger, M.}
%{From subfactors to categories and topology I. Frobenius algebras
%in and Morita equivalence of tensor categories}\JPAA{180}
%{81--157} {(2003) or math.CT/0111204}.

\bibx{m1} {M\"uger, M.} {From subfactors to categories and topology II. The quantum double of 
tensor categories and subfactors.} {J. Pure Appl. Alg.}  {\bf 180},  {159--219}  {(2003)}.

\bibx{m3} {M\"uger, M.} 
{Tensor categories: a selective guided tour.}
{Rev. Un. Mat. Argentina} {\bf 51}, 95Ð163  (2010).

%\biba{NW} {Nobs, A., Wolfart, J.} {Die irreduziblen Darstellungen der Gruppen SL$_2(Z_p)$
%insbesondere SL$_2(Z_2)$ II} {Comment. Math. Helvetici} {51}  {491--526} {(1976)}.

%\bibx{Ost} {Ostrik, V.} {Fusion categories of rank 2.} {Math. Res. Lett.} \textbf{10}, 177--183 {(2003)}.

%\bibx{ocn3}{Ocneanu, A.}{Operator algebras, 3-manifolds and quantum field
%theory.}{Talk given in Kyoto, 1991}.

%\bibx{Ocriral}{Ocneanu, A.}{Chirality for operator algebras.}
%{{\it Subfactors} (Kyuzeso, 1993), World Sci.\ Publishing, River
%Edge, NJ, pp.\  39--63, 1994}.

%\biba{ost}{Ostrik, V.}
%{Module categories over the Drinfeld double of a finite group}
%{Inter.\ Math.\ Research Notices} {27} {1507--1520}{(2003)};
%{arXiv: math.QA/0202130}.

%\bibx{pet} {Peters, E.} {A planar algebra construction of the Haagerup subfactor.}  {Internat. J. Math.} \textbf{21}, {987--1045} {(2010)}; {arXiv: 
%0902.1294}. 

\bibx{Phil} {Phillips, N. C.} {Analogs of Cuntz algebras on $L^p$ spaces.} {arXiv:1201.4196.}

\bibx{Phil1} {Phillips, N. C.} {Private communication.}


\bibx{Rehren}
 {Rehren, K.-H.} {Chiral observables and modular invariants.} {Commun. Math. Phys.} {\bf 208} (2000), 689--712.


%\bibx{Sei} {Seitz, G.} {Finite groups having only one irreducible representation of degree
%greater than one.} {Proc.\ Amer.\ Math.\ Soc.} \textbf{19}, {459--461} {(1968)}. 

%\bibx{Sie} {Siehler, J.} {Near-group categories.} {Alg. Geom. Topol.} \textbf{3}, {719--775} {(2003)}.

%\biba{PZ}
%{Petkova, V.B., Zuber, J.B.} {The many faces of Ocneanu
%cells}\nupb{603}{449--496}{(2001)}.

%\bibx{PPP1}{Pinto, P.R.}{PhD Thesis.}{Cardiff University,\ 2003}.

%\bibx{R1}
%{Rehren, K.-H.} {Braid group statistics and their superselection
%rules.} {In: {The algebraic theory of superselection sectors}
%(Palermo 1989).  World Scientific, Singapore 1990, pp.\ 333--355}.

%\biba{rehren} {Rehren, K.-H.} {Space-time fields and exchange fields}\CMP{132}
%{461--483}{(1990)}.

%\bibx{SchWa} {Schellekens, A.N., Warner, N.P.} {Conformal subalgebras of Kac-Moody algebras.}
%Phys. Rev. \textbf{D34}, 3092--3096 (1986).
%\bibx{serre}
%{Serre, J.-P.} {Linear Representations of Finite Groups.}
%{Graduate Texts in Mathematics, Vol.\ 42. Springer-Verlag,
%New York-Heidelberg, 1977.}

%\biba{suth}{Sutherland, C.}
%{Cohomology and extensions of von Neumann algebras. I,
%II}\RIMS{16} {135--174} {(1980)}.

%\bibx{Tak} {Takagi, T.} {}  {Japanese J. Math.} {1} (1927) 83.

%\bibx{TY} {Tambara, D., Yamagami, S.} {Tensor categories with fusion rules of self-duality
%for finite abelian groups.} {J. Algebra} \textbf{209},  692--707 (1998).

%\bibx{Thor} {Thornton, J.} {On braided near-group categories.} arXiv: 1102.4640v1.

%\bibx{tur}
%{Turaev, V.G.}{Quantum Invariants of Knots and 3-manifolds.} 
%{de Gruyter Studies in Mathematics}, {vol 18.} (Walter de
%Gruyter, Berlin, 1994).

%\bibx{Ve} {Verlinde, E.}
%{Fusion rules and modular transformations in 2D conformal field
%theory.}  Nucl. Phys. \textbf{B300}, {360--376} {(1988)}.

%\bibx{Wal} {Walton, M.A.} {Conformal branching rules and modular invariants.} Nucl. Phys. 
%\textbf{B322}, {775--790} {(1989)}.

%\bibx{W1} {Wassermann, A.}
%{Coactions and Yang-Baxter equations for ergodic actions and subfactors.} {In: Operator algebras and applications, Vol. 2,} {203--236,} 
%{London Math. Soc. Lecture Note Ser., 136} (Cambridge Univ. Press, Cambridge, 1988).

%\bibx{W2}{Wassermann, A.}
%{Quantum subgroups and vertex algebras.}{MSRI talk, December
%2000,\par
%http://www.msri.org/publications/ln/msri/2000/subfactors/wassermann/1}.

%\bibx{Wi} {Witten, E.} {The search for higher symmetry in string theory. Physics and mathematics 
%of strings.} Philos. Trans. Roy. Soc. London Ser. \textbf{A329} (1989), no. 1605, 349--357. 


%\bibx{Xu} {Xu, F.} {Mirror extensions of local nets.}  Commun. Math. Phys. \textbf{270}, 835--847 (2007). 

\bibx{Xu} {Xu, F.}
{Algebraic orbifold conformal field theories.}  
Proc. Natl. Acad. Sci. USA 97, no. 26, 14069--14073 (2000). 

\bibx{Ya} {Yamagami, S.} {A categorical and diagrammatical approach to Temperley--Lieb algebras;} arXiv:0405267.
{}

%\bibx{zhu} {Zhu, Y.}  
%{Modular invariance of characters of vertex operator algebras.}
%J. Amer. Math. Soc. \textbf{9}, 237--302 (1996). 


\end{scriptsize}

\end{thebibliography}
\end{document}